\documentstyle[amscd]{amsart}
\numberwithin{equation}{section}
\theoremstyle{plain}
\newtheorem{thm}{Theorem}[section]
\newtheorem{cor}[thm]{Corollary}
\newtheorem{lem}[thm]{Lemma}
\newtheorem{prop}[thm]{Proposition}
\begin{document}
\title{$C^*$-algebras associated with Hilbert $C^*$-quad modules 
of $C^*$-textile dynamical systems
}
\author{Kengo Matsumoto}
\address{ 
Department of Mathematics, 
Joetsu University of Education,
Joetsu 943-8512, Japan}
\email{kengo@@juen.ac.jp}
\maketitle
\begin{abstract}
A $C^*$-textile dynamical system
$({\cal A}, \rho,\eta,\Sigma^\rho,\Sigma^\eta, \kappa)$
connsists of a unital $C^*$-algebra ${\cal A}$,
two families of endomorphisms
$\{ \rho_\alpha \}_{\alpha \in \Sigma^\rho}$
and
$\{ \eta_a \}_{a \in \Sigma^\eta}$
of ${\cal A}$ and
certain commutation relations $\kappa$
among them.
It yields a two-dimensional subshift
and 
multi structure of 
Hilbert $C^*$-bimodules,
which we call a Hilbert $C^*$-quad module.
We introduce a $C^*$-algebra from the Hilbert $C^*$-quad module
as a two-dimensional analogue of Pimsner's construction of $C^*$-algebras from Hilbert $C^*$-bimodules.
We study the $C^*$-algebras defined by  the Hilbert $C^*$-quad modules
and prove that they have  universal properties subject to certain operator relations.
We also present its examples arising from commuting matrices.
\end{abstract}


\def\Zp{{ {\Bbb Z}_+ }}
\def\CTDS{{ ({\cal A}, \rho, \eta, \Sigma^\rho, \Sigma^\eta, \kappa)}}
\def\OL{{{\cal O}_{{\frak L}}}}
\def\FL{{{\cal F}_{{\frak L}}}}
\def\FKL{{{\cal F}_k^{l}}}
\def\FHK{{{\cal F}_{{\cal H}_\kappa}}}
\def\FU{{{\cal F}_{{\cal H}_\kappa}^{uni}}}
\def\FHLI{{{\cal F}_{[{\frak L}]} }}
\def\FRHO{{ {\cal F}_\rho}}
\def\DRHO{{{{\cal D}_{\rho}}}}
\def\ORHO{{{\cal O}_\rho}}
\def\OETA{{{\cal O}_\eta}}
\def\ORE{{{\cal O}^{\kappa}_{\rho,\eta}}}
\def\OREK{{{\cal O}_{(\rho,\eta;\kappa)}}}
\def\OHK{{{\cal O}_{{\cal H}_\kappa}}}
\def\OU{{{\cal O}_{{\cal H}_\kappa}^{uni}}}
\def\OAK{{{\cal O}_{\A,\kappa}}}
\def\OHRE{{{\widehat{\cal O}}^{\kappa}_{\rho,\eta}}}
\def\FRE{{{\cal F}_{\rho,\eta}}}
\def\DRE{{{\cal D}_{\rho,\eta}}}
\def\M{{{\cal M}}}
\def\N{{{\cal N}}}
\def\HK{{{\cal H}_\kappa}}
\def\EK{{E_\kappa}}
\def\SK{{\Sigma_\kappa}}
\def\K{{{\cal K}}}
\def\A{{{\cal A}}}
\def\B{{{\cal B}}}
\def\BE{{{\cal B}_\eta}}
\def\BR{{{\cal B}_\rho}}
\def\BK{{{\cal B}_\kappa}}
\def\BU{{{\cal B}_\kappa^{uni}}}
\def\Aut{{{\operatorname{Aut}}}}
\def\End{{{\operatorname{End}}}}
\def\Hom{{{\operatorname{Hom}}}}
\def\Ker{{{\operatorname{Ker}}}}
\def\ker{{{\operatorname{ker}}}}
\def\Coker{{{\operatorname{Coker}}}}
\def\id{{{\operatorname{id}}}}
\def\dim{{{\operatorname{dim}}}}
\def\min{{{\operatorname{min}}}}
\def\exp{{{\operatorname{exp}}}}
\def\supp{{{\operatorname{supp}}}}
\def\Proj{{{\operatorname{Proj}}}}
\def\Im{{{\operatorname{Im}}}}
\def\span{{{span}}}
\def\OA{{{\cal O}_A}}
\def\OB{{{\cal O}_B}}
\def\T{{  {\cal T}_{{\K}^{\M}_{\N}} }}
\def\TK{{ {\cal T}_{{\cal H}_\kappa} }}

\section{Introduction}
In \cite{MaDoc99}, the author has introduced a notion of 
$\lambda$-graph system as 
a generalization of finite labeled graphs.
The $\lambda$-graph systems yield $C^*$-algebras so that its  K-theory groups  are
related to topological conjugacy invariants of the underlying symbolic dynamical systems.
He has extended the notion of $\lambda$-graph system to  $C^*$-symbolic dynamical system,
which 
is a generalization of both
a $\lambda$-graph system and an automorphism of a unital $C^*$-algebra.
It is denoted by $(\A, \rho,\Sigma)$
and consists of a finite family $\{ \rho_{\alpha} \}_{\alpha \in \Sigma}$ 
of endomorphisms of a unital $C^*$-algebra 
${\cal A}$ 
such that 
$\rho_\alpha(Z_\A) \subset Z_\A, \alpha \in \Sigma$
and
$\sum_{\alpha\in \Sigma}\rho_\alpha(1) \ge 1$
where
$Z_\A$ denotes the center of $\A$.
A $\lambda$-graph system ${\frak L}$ 
yields a $C^*$-symbolic dynamical system 
$({\cal A}_{\frak L},\rho^{\frak L}, \Sigma)$ 
such that ${\cal A}_{\frak L}$ 
is $C(\Omega_{\frak L})$ for some compact Hausdorff space 
$\Omega_{\frak L}$ with 
$\dim \Omega_{\frak L} = 0$.
A $C^*$-symbolic dynamical system $({\cal A},\rho, \Sigma)$
provides a subshift 
$\Lambda_\rho$  
over $\Sigma$
and a Hilbert $C^*$-bimodule 
${\cal H}_{\cal A}^{\rho}$
over $\A$
which gives rise to a  $C^*$-algebra
${\cal O}_\rho$ as a Cuntz-Pimsner algebra 
(\cite{MaCrelle}, cf. \cite{KPW}, \cite{Pim}).

G. Robertson--T. Steger \cite{RoSte}
have initiated a certain study of higher dimensional analogue of Cuntz--Krieger algebras
from the view point of tiling systems of 2-dimensional plane.
After their work, 
A. Kumjian--D. Pask \cite{KP}
have generalized their construction to introduce the notion of 
higher rank graphs and its $C^*$-algebras.
Since then, there have been 
many studies  on these $C^*$-algebras by many authors
(see for example \cite{Dea}, \cite{EGS}, \cite{KP}, \cite{PCHR}, \cite{PRW1}, \cite{RoSte}, etc.).

M. Nasu in \cite{NaMemoir} has introduced the notion of textile system 
which is useful in analyzing automorphisms and endomorphisms of topological Markov shifts.
A textile system also gives rise to a two-dimensional tiling called 
Wang tiling. 
Among textile systems,
 LR textile systems have specific properties 
that consist of two commuting symbolic matrices. 
In \cite{MaYMJ2008},
the author has extended the notion of textile systems to $\lambda$-graph systems
and has defined a notion of textile systems on $\lambda$-graph systems,
which are called textile $\lambda$-graph systems for short.
 $C^*$-algebras associated to  textile systems  
have been initiated by V. Deaconu (\cite{Dea}).

In \cite{MaPre2011}, 
the author has extended the notion of
$C^*$-symbolic dynamical system to 
$C^*$-textile  dynamical system
which is a higher dimensional analogue of $C^*$-symbolic dynamical system.
The $C^*$-textile  dynamical system
$({\cal A}, \rho, \eta, \Sigma^\rho, \Sigma^\eta, \kappa)$
 consists of two 
$C^*$-symbolic dynamical systems 
$({\cal A}, \rho,  {\Sigma^\rho})$
and 
$({\cal A}, \eta, {\Sigma^\eta})$
with a common unital $C^*$-algebra $\A$ and a  commutation relation
between $\rho$ and $\eta$ through a map $\kappa$ below.
Set 
$$
\Sigma^{\rho \eta} = \{ (\alpha, b ) \in \Sigma^\rho \times \Sigma^\eta \mid 
\eta_b \circ \rho_\alpha  \ne 0 \},\quad
\Sigma^{\eta \rho} = \{ (a, \beta) \in \Sigma^\eta \times \Sigma^\rho \mid 
\rho_\beta \circ \eta_a  \ne 0 \}.
$$
We require that there exists a bijection
$
\kappa : \Sigma^{\rho \eta} \longrightarrow \Sigma^{\eta\rho},
$
which we fix and call a specification.
Then the required commutation relations are 
\begin{equation}
\eta_b \circ \rho_\alpha = \rho_\beta \circ \eta_a
\qquad
\text{ if } 
\quad 
\kappa(\alpha, b) = (a,\beta). \label{eqn:kappa}
\end{equation}
The author has also
 introduced a $C^*$-algebra 
  $\ORE$
from 
$({\cal A}, \rho, \eta, \Sigma^\rho, \Sigma^\eta, \kappa)$
which is realized as the universal $C^*$-algebra 
$C^*(x, S_\alpha, T_a ; x \in \A, \alpha \in \Sigma^\rho, 
a \in \Sigma^\eta)$  
generated by 
$x \in \A$ 
and two families of partial isometries 
$S_{\alpha}, \alpha \in \Sigma^\rho$,
$T_a, a \in \Sigma^\eta$
subject to the following relations called $(\rho,\eta;\kappa)$:
\begin{align}
\sum_{\beta \in \Sigma^\rho}S_{\beta}S_{\beta}^*  =1,\qquad
x S_\alpha S_\alpha^* & =  S_\alpha S_\alpha^*  x,\qquad
S_\alpha^* x S_\alpha  = \rho_\alpha(x), \label{eqn:S} \\
\sum_{b \in \Sigma^\eta} T_{b} T_{b}^*  =1,\qquad
x T_a T_a^*  & =  T_a T_a^*  x,\qquad
T_a^* x T_a  = \eta_a(x), \label{eqn:T}\\  
S_\alpha T_b  = T_a S_\beta &
\qquad
\text{ if } 
\quad
\kappa(\alpha, b) = (a,\beta) \label{eqn:STTS}
\end{align}
for all $x \in {\cal A}$ and $\alpha \in \Sigma^\rho, a \in \Sigma^\eta$
(\cite{MaPre2011}).
The algebra is a generalization of some of higher rank graph algebras.

In the present paper,
the author will introduce another kind of $C^*$-algebras associated with the $C^*$-textile dynamical systems
from the view point of Hilbert $C^*$-modules.
The resulting $C^*$-algebras $\OHK$ are different from the above algebras $\ORE$.
A $C^*$-textile dynamical system
 provides  a two-dimensional subshift and 
multi structure of Hilbert $C^*$-bimodules that have  multi right actions and 
multi left actions and multi inner products.
We call it a 
 Hilbert $C^*$-quad module
 denoted  by
 $\HK$.
The  $C^*$-algebra $\OHK$, which we will introduce in the present paper, 
is constructed in a concrete way from 
the  structure of the Hilbert $C^*$-quad module
$\HK$ by a two-dimensional analogue of 
Pimsner's construction from Hilbert $C^*$-bimodules.
It is generated by the quotient images 
of creation operators on two-dimensional analogue of Fock Hilbert 
module by module maps of compact operators. 
As a result,
we will show the $C^*$-algebra has a universal property 
subject to certain operator relations 
of generators.

For a $C^*$-textile dynamical system
$\CTDS$, 
consider the set of quadruplet of symbols 
\begin{equation}
\SK =
\{\omega = (\alpha, b, a,\beta) 
\in \Sigma^\rho\times \Sigma^\eta\times \Sigma^\eta\times \Sigma^\rho
\mid \kappa(\alpha,b) = (a,\beta) \}.  \label{eqn:SK}
\end{equation}
Each element of $\SK$
is regarded as a tile \, \, 
{\footnotesize 
$
\begin{CD}
\cdot @>\alpha>> \cdot \\
@V{a}VV  @VV{b}V \\
\cdot @>>\beta> \cdot 
\end{CD}
$
} \,
of the associated two-dimensional subshift.
Denote by $\ORHO$ and by  $\OETA$
the $C^*$-algebras associated with the 
$C^*$-symbolic dynamical systems $(\A, \rho, \Sigma^\rho)$
and $(\A, \eta, \Sigma^\eta)$
respectively.
 Let
$S_\alpha, \alpha \in \Sigma^\rho$
and
$T_a, a \in \Sigma^\eta$
be the generating partial isometries 
of $\ORHO$ 
and those  
of $\OETA$, which satisfy 
\eqref{eqn:S} and \eqref{eqn:T} respectively.
Denote 
by $\BR$  
the $C^*$-subalgebra of $\ORHO$
generated by elements 
$
S_\alpha x S_\alpha^* , \alpha \in \Sigma^\rho,  x \in \A 
$
and 
by $\BE$ 
that of $\OETA$
generated by elements 
$
T_a x T_a^* , a \in \Sigma^\eta,  x \in \A
$
respectively.
The endomorphisms
$\rho_\alpha, \alpha \in \Sigma^\rho$
and
$\eta_a, a \in \Sigma^\eta
$
on $\A$
extend to
$\BR$ and to $\BE$
as $*$-homomorphisms
$\widehat{\rho}_\alpha: \BR \longrightarrow \A$
and
$\widehat{\eta}_a: \BE \longrightarrow \A$
by 
\begin{equation}
\widehat{\rho}_\alpha(w)  = S_\alpha^* w S_\alpha \in \A,
\quad
 w \in \BR
 \quad 
 \text{and}
 \quad
\widehat{\eta}_a(z)  = T_a^* z T_a \in \A,
\quad
 z \in \BE. \label{eqn:homo}
\end{equation}
They also extend to
$\BE$ and to $\BR$ as 
endomorphisms 
$\widehat{\rho}_\alpha^\eta
$
on $\BE$ 
and
$\widehat{\eta}_a^\rho
$
on 
$\BR$
by 
\begin{align}
\widehat{\rho}_\alpha^\eta(z) & =
\sum
\begin{Sb}
b,a,\beta\\
(\alpha,b,a,\beta)\in\SK
\end{Sb}
T_b \rho_\beta( \widehat{\eta}_a(z)) T_b^* \in \BE,
\qquad z \in \BE, \\
\widehat{\eta}_a^\rho(w) &  =
\sum
\begin{Sb}
\alpha, b,\beta\\
(\alpha,b,a,\beta)\in\SK
\end{Sb}
S_\beta \eta_b( \widehat{\rho}_\alpha(w)) S_\beta^* \in \BR,
\qquad w \in \BR.
\end{align}
For 
$\omega \in \SK$, 
put a projection
$E_\omega = \eta_b(\rho_\alpha(1)) ( = \rho_\beta(\eta_a(1)) \in \A$.
The  vector space
\begin{equation}
\HK = \sum_{\omega \in \Sigma_\kappa} {\Bbb C}e_\omega \otimes E_\omega \A
\end{equation}
has a natural structure of  a Hilbert $C^*$-right $\A$-module. 
In addition to the $\A$-module structure,
$\HK$ has multi structure of Hilbert $C^*$-bimodules,
a Hilbert $C^*$-bimodule structure over $\BR$ and 
a Hilbert $C^*$-bimodule structure over $\BE$.
We call it Hilbert $C^*$-quad module over $(\A;\BR,\BE)$.
We will construct a $C^*$-algebra
$\OHK$ in a concrete way from 
the Hilbert $C^*$-quad module
$\HK$ by a two-dimensional analogue of 
Pimsner's construction of $C^*$-algebras from Hilbert $C^*$-bimodules.
It is generated by two kinds of creation operators, 
the horizontal creation operators and the vertical creation operators,  on two-dimensional analogue of Fock Hilbert module. 
Denote by
$\iota_\rho: \A \hookrightarrow \BR$
and 
$\iota_\eta: \A \hookrightarrow \BE$
natural embeddings.
We assume that the algebra $\A$ is commutative.
The main result of the paper is 
the following theorem, which states that the algebraic structure of the algebra 
$\OHK$ is determined by the behavior of the $*$-homomorphisms
$\widehat{\rho}_\alpha, \widehat{\eta}_a, \widehat{\rho}^\eta_\alpha$ and
$\widehat{\eta}^\rho_a$.  
\begin{thm}[Theorem \ref{thm:main}]
For a $C^*$-textile dynamical system 
$\CTDS$,
the $C^*$-algebra 
$\OHK$ associated with the 
Hilbert $C^*$-quad module $\HK$
is realized as  
the universal concrete $C^*$-algebra
generated by the operators
$z\in \BE, w\in \BR $
and partial isometries
${\tt u}_\alpha, \alpha \in \Sigma^\rho,  
{\tt v}_a, a \in \Sigma^\eta
$
subject to the relations: 
\begin{align*}
\sum_{\beta \in \Sigma^\rho} {\tt u}_\beta {\tt u}_\beta^* +
& \sum_{b \in \Sigma^\eta} {\tt v}_b {\tt v}_b^* =1, \\ 
{\tt u}_\alpha {\tt u}_\alpha^* w  = w {\tt u}_\alpha {\tt u}_\alpha^*, 
& \qquad
{\tt v}_a {\tt v}_a^* w  = w {\tt v}_a {\tt v}_a^*, \\
{\tt u}_\alpha {\tt u}_\alpha^* z  = z {\tt u}_\alpha {\tt u}_\alpha^*, 
& \qquad
{\tt v}_a {\tt v}_a^* z = z {\tt v}_a {\tt v}_a^*, \\
\widehat{\rho}_\alpha(w)
=
{\tt u}_\alpha^* w {\tt u}_\alpha, 
&  \qquad
\widehat{\eta}_a(z)
=
{\tt v}_a^* z {\tt v}_a,\\ 
\widehat{\rho}^\eta_\alpha(z)
=
{\tt u}_\alpha^* z {\tt u}_\alpha, 
&  \qquad
\widehat{\eta}^\rho_a(w)
=
{\tt v}_a^* w {\tt v}_a,\\ 
\iota_\eta(y) & = \iota_\rho(y) 
\end{align*}
for
$
w \in \BR, z \in \BE, \alpha \in \Sigma^\rho, a \in \Sigma^\eta,
y \in \A$.
\end{thm}
Thanks to the above theorem,
simplicity condition of the $C^*$-algebra 
$\OHK$ will be presented
(Theorem \ref{cor:simplicity}).

Let
$A,B$ be two $N \times N$ matrices with entries in 
nonnegative integers.
They yield directed graphs
$G_A = (V, E_A)$
and 
$G_B = (V, E_B)$
with a common vertex set
$V = \{ v_1,\dots, v_N \}$ and 
edge sets $E_A$ and 
$E_B$ respectively,
where
the edge set 
$E_A$ consist of $A(i,j)$-edges from the vertex 
$v_i$ to the vertex $v_j$
and 
$E_B$ consist of $B(i,j)$-edges from the vertex 
$v_i$ to the vertex $v_j$.
We then have two
$C^*$-symbolic dynamical systems
$(\A_N, \rho^A, E_A)$
and
$(\A_N, \rho^B, E_B)$
with
$\A_N = {\Bbb C}^N$.
Denote by $s(e), r(e)$
the source vertex and the range vertex of an edge $e$.
Put
\begin{align*}
\Sigma^{AB}
& =\{(\alpha,b)\in E_A \times E_B 
\mid  r(\alpha) = s(b)  \}, \\ 
\Sigma^{BA}
& =\{(a,\beta)\in E_B \times E_A 
\mid  r(a) = s(\beta) \}.
\end{align*}
Assume that the commutation relation
\begin{equation}
AB = BA  \label{eqn:abba}
\end{equation} 
holds.
We may take a bijection 
$\kappa: \Sigma^{AB} \longrightarrow \Sigma^{BA}$
such that $s(\alpha) = s(a), r(b) = r(\beta)$
for $\kappa(\alpha,b) = (a,\beta)$
which we fix.
This situation is called an LR-textile system introduced by Nasu
(\cite{NaMemoir}).
We then have a $C^*$-textile dynamical system
$(\A_N, \rho^A,\rho^B, E_A, E_B, \kappa)$.
We set
\begin{equation*}
\Omega_\kappa
=\{  (\alpha,a) \in   E_A\times E_B |  s(\alpha) = s(a),   
\kappa(\alpha,b) =(a,\beta) \text{ for some } \beta \in E_A, b\in E_B \} 
\end{equation*}
and define two $|\Omega_\kappa| \times |\Omega_\kappa|$-matrcies
$A_\kappa$ and $B_\kappa$  with entries in $\{0,1\}$
by 
\begin{align*}
A_\kappa((\alpha,a),(\delta,b))
& = 
{\begin{cases}
1 & \text{ if there exists } \beta\in E_A  \text{ such that } \kappa(\alpha,b) = (a,\beta),\\  
0 & \text{ otherwise}
\end{cases}
}
\\
\intertext{ for
$(\alpha,a),(\delta,b) \in \Omega_\kappa$, and} 
B_\kappa((\alpha,a), (\beta,d))           
& = 
{\begin{cases}
1 & \text{ if there exists } b \in E_B  \text{ such that } \kappa(\alpha,b) = (a,\beta),\\  
0 & \text{ otherwise}
\end{cases}
}
 \end{align*}
for 
$(\alpha,a), (\beta,d) \in \Omega_\kappa$
respectively.
Denote by
${\cal H}_\kappa^{A,B}$ the associated Hilbert $C^*$-quad module.
\begin{thm}[Theorem \ref{thm:CK}]
The $C^*$-algebra ${\cal O}_{{\cal H}_\kappa^{A,B}}$
associated with the 
Hilbert $C^*$-quad module ${\cal H}_\kappa^{A,B}$
defined by commuting matrices
$A, B$
and a specification $\kappa$
is 
generated by 
two families of  partial isometries
$
S_{(\alpha,a)}, T_{(\alpha,a)}
$ 
for
$(\alpha,a)\in \Omega_\kappa$
satisfying the  relations: 
\begin{align*}
\sum_{(\delta,b) \in\Omega_\kappa} 
&S_{(\delta,b)}S_{(\delta,b)}^*
+
\sum_{(\beta,d) \in\Omega_\kappa} 
T_{(\beta,d)}T_{(\beta,d)}^* = 1, \\
S_{(\alpha,a)}^*S_{(\alpha,a)}
=
&
\sum_{(\delta,b) \in\Omega_\kappa}
A_\kappa((\alpha,a),(\delta,b)) (
S_{(\delta,b)}S_{(\delta,b)}^* + T_{(\delta,b)} T_{(\delta,b)}^*),\\
T_{(\alpha,a)}^* T_{(\alpha,a)}
=
&
\sum_{(\beta,d) \in\Omega_\kappa} B_\kappa((\alpha,a),(\beta,d)) 
(
S_{(\beta,d)}S_{(\beta, d)}^* + T_{(\beta,d)} T_{(\beta,d)}^*)
\end{align*}
for $(\alpha,a) \in \Omega_\kappa$.
Hence the $C^*$-algebra ${\cal O}_{{\cal H}_\kappa^{A,B}}$
is $*$-isomorphic to the Cuntz-Krieger algebra
${\cal O}_{H_\kappa}$
for the matrix 
$
H_\kappa
=
\begin{bmatrix}
A_\kappa & A_\kappa \\
B_\kappa & B_\kappa 
\end{bmatrix}.
$
\end{thm}
The paper is organized as in the following way:
In Section 2, we will state basic facts on the $C^*$-symbolic dynamical systems and
the $C^*$-textile dynamical systems.
In Section 3, we will introduce Hilbert $C^*$-quad modules from 
$C^*$-textile dynamical systems.
In Section 4, we will introduce Fock Hilbert $C^*$-quad modules
which are two-dimensional analogue of Fock Hilbert $C^*$-bimodules,
and study creation operators on 
the Fock Hilbert $C^*$-quad modules.
In Section 5, we will prove the main result stated as Theorem 1.1.
In Section 6, we will state a relationship between the $C^*$-algebras 
$\OHK$ and $\ORE$
so that 
the algebra  
$\OHK$ is realized as a $C^*$-subalgebra of the tensor product
$\ORE \otimes {\cal O}_2$
in a natural way.
In Section 7, we will  study the $C^*$-algebras arising from the 
Hilbert $C^*$-quad modules of the $C^*$-textile dynamical systems defined by commuting matrices
and will prove Theorem 1.2. 

Throughout the paper,
we will denote by $\Zp$ 
the set of nonnegative integers 
and
by ${\Bbb N}$
the set of positive integers.

\section{ $C^*$-symbolic dynamical systems
and  $C^*$-textile dynamical systems}
In this section, we will briefly state basic facts on  
 $C^*$-symbolic dynamical systems and $C^*$-textile dynamical systems.
Throughout the section,
$\Sigma$ denotes a finite set with its discrete topology, 
that is called an alphabet.
Each element of $\Sigma$ is called a symbol.
Let $\Sigma^{\Bbb Z}$ 
be the infinite product space 
$\prod_{i\in {\Bbb Z}}\Sigma_{i},$ 
where 
$\Sigma_{i} = \Sigma$,
 endowed with the product topology.
 The transformation $\sigma$ on $\Sigma^{\Bbb Z}$ 
given by 
$
\sigma( (x_i)_{i \in {\Bbb Z}})
 = (x_{i+1})_{i \in {\Bbb Z}}
$ 
is called the full shift over $\Sigma$.
 Let $\Lambda$ be a shift invariant closed subset of $\Sigma^{\Bbb Z}$ i.e. 
 $\sigma(\Lambda) = \Lambda$.
  The topological dynamical system 
  $(\Lambda, \sigma\vert_{\Lambda})$
   is called a two-sided subshift, written as $\Lambda$ for brevity.
Finite directed graphs present 
a class of subshifts called shifts of finite type. 
More generally, 
finite directed labeled graphs present 
a class of subshifts called sofic shifts.
The author has introduced a notion of $\lambda$-graph system
as a generalization of finite labeled graphs.
The $\lambda$-graph systems present all the  subshifts.
Furthermore,
the author has introduced a notion of $C^*$-symbolic dynamical system
which generalize $\lambda$-graph systems and automorphisms of unital $C^*$-algebras.
$C^*$-symbolic dynamical systems may be presentation of subshifts
to $C^*$-algebras.

Let ${\cal A}$ be a unital $C^*$-algebra.
In what follows,
an endomorphism of $\cal A$ means 
a $*$-endomorphism of $\cal A$ that does not necessarily preserve the unit
$1_\A$ of 
$\cal A$.
The unit $1_\A$ is denoted by $1$ unless we specify.
Denote by 
$Z_\A$ the center of $\A$.
Let $\rho_\alpha, \alpha\in \Sigma$
be a finite family of endomorphisms
of $\A$ indexed by symbols of a finite set $\Sigma$. 
We assume that
$\rho_\alpha(Z_\A) \subset Z_\A, \alpha \in \Sigma$.
The family 
$\rho_\alpha, \alpha\in \Sigma$
of endomorphisms of $\A$ 
is said to be essential if
$\rho_{\alpha}(1) \ne 0$ for all $\alpha \in \Sigma$
 and   
$\sum_{\alpha}\rho_\alpha(1) \ge 1$.
It is said to be faithful if for any nonzero $x \in \cal A$ 
there exists a symbol $\alpha\in \Sigma$ such that $\rho_{\alpha}(x) \ne 0$.
A $C^*$-{\it symbolic dynamical system}\ 
is a triplet $({\cal A}, \rho, \Sigma)$ 
consisting of a unital $C^*$-algebra $\A$ 
and 
an essential and faithful finite family 
$\{ \rho_{\alpha} \}_{\alpha \in \Sigma}$  of endomorphisms
of ${\cal A}$.
In \cite{MaCrelle}, \cite{MaContem}, \cite{MaMZ2010},
we have defined a $C^*$-symbolic dynamical system in a less restrictive way than the above definition.
In stead of the above condition 
$\sum_{\alpha \in \Sigma}\rho_\alpha(1) \ge 1$
with
$\rho_\alpha(Z_\A) \subset Z_\A, \alpha \in \Sigma$,
we have used the condition in the  papers  that 
the closed ideal generated by $\rho_\alpha(1), \alpha \in \Sigma$ 
coincides with $\A$.
 All of the examples appeared in the papers
\cite{MaCrelle}, \cite{MaContem}, \cite{MaMZ2010}
satisfy the condition
$\sum_{\alpha \in \Sigma}\rho_\alpha(1) \ge 1$
with
$\rho_\alpha(Z_\A) \subset Z_\A, \alpha \in \Sigma$,
and all discussions in the papers well work under the new defition.

A $C^*$-symbolic dynamical system $({\cal A}, \rho, \Sigma)$
 yields  a subshift $\Lambda_\rho$
over $\Sigma$ such that a word 
$\alpha_1\cdots\alpha_k$ of $\Sigma$ is admissible for $\Lambda_\rho$
 if and only if 
$(\rho_{\alpha_k}\circ \cdots\circ \rho_{\alpha_1})(1) \ne 0$
(\cite[Proposition 2.1]{MaCrelle}).
Denote by 
$B_k(\Lambda_\rho)$
 the set of admissible words of
$\Lambda_\rho$ with length $k$.
Put
$
B_*(\Lambda_\rho) = \cup_{k=0}^{\infty}B_k(\Lambda_\rho),
$
where $B_0(\Lambda_\rho) $ 
consists of the empty word. 
The $C^*$-algebra $\ORHO$ 
associated with 
$(\A,\rho,\Sigma)$
has been originally constructed in \cite{MaCrelle} 
from an associated  Hilbert $C^*$-bimodule
(cf. \cite{Pim}, \cite{KPW} etc.).
It is realized as a universal $C^*$-algebra
$C^*(x, S_{\alpha}; x \in \A, \alpha \in \Sigma)$  
generated by 
$x \in \A$ 
and partial isometries $S_{\alpha}, \alpha \in \Sigma$
subject to the following relations called $(\rho)$:
$$
\sum_{\beta \in \Sigma}S_{\beta}S_{\beta}^* =1,\qquad
x S_\alpha S_\alpha^* =  S_\alpha S_\alpha^*  x,\qquad
S_\alpha^* x S_\alpha = \rho_\alpha(x)  
$$
for all $x \in \cal A$ and $\alpha \in \Sigma.$
The $C^*$-algebra $\ORHO$ is a generalization of 
the $C^*$-algebra ${\cal O}_{\frak L}$
associated with the $\lambda$-graph system $\frak L$
(cf. \cite{DocMath02}).


Let 
$\CTDS$ be a 
$C^*$-{\it textile dynamical system}.
It  consists of two 
$C^*$-symbolic dynamical systems 
$({\cal A}, \rho,  {\Sigma^\rho})$
and 
$({\cal A}, \eta, {\Sigma^\eta})$
with common unital $C^*$-algebra $\A$
and commutation relations
between their endomorphisms 
$\rho_\alpha, \alpha\in \Sigma^\rho$
and 
$\eta_a, a\in \Sigma^\eta$
through a bijection $\kappa$
satisfying 
\eqref{eqn:kappa}.
Let
$S_\alpha, \alpha \in \Sigma^\rho$
and
$T_a, a \in \Sigma^\eta$
be the generating partial isometries 
of $\ORHO$ and 
of $\OETA$, 
which satisy
\eqref{eqn:S} and \eqref{eqn:T}
respectively. 
We set two $C^*$-algebras
\begin{equation*}
\BR  = C^*(S_\alpha x S_\alpha^* : \alpha \in \Sigma^\rho,  x \in \A ) ,
\qquad
\BE  = C^*(T_a x T_a^* : a \in \Sigma^\eta,  x \in \A ). 
\end{equation*}
They are realized concretely as subalgebras of $\ORHO$ and of $\OETA$
respectively. 
Both the algebras
$\BR$ and $\BE$ contain the algebra
$\A$ through the identification
\begin{equation}
x =
\sum_{\alpha \in \Sigma^\rho} S_\alpha \rho_\alpha(x)S_\alpha^*
=
\sum_{a \in \Sigma^\eta} T_a \eta_a(x)T_a^*,
\qquad 
x \in \A.
\end{equation}
We  put the projections
$$
P_\alpha =\rho_\alpha(1)
\text{ for } \alpha \in \Sigma^\rho,
\qquad 
Q_a =\eta_a(1)
\text{ for } a \in \Sigma^\eta.
$$
Elements  $w \in \BR$ and 
$z \in \BE$
are uniquely written as in the following way:
\begin{align}
w=
& \sum_{\alpha \in \Sigma^\rho} S_\alpha w_\alpha S_\alpha^*
 \text{ with } w_\alpha = P_\alpha w_\alpha P_\alpha \in \A, \, \alpha \in \Sigma^\rho,
\label{eqn:w} \\
z =
& \sum_{a \in \Sigma^\eta}T_a z_a T_a^*
\text{ with } z_a = Q_a z_a Q_a \in \A, \, a \in \Sigma^\eta. \label{eqn:z}
\end{align}
Define an alphabet set $\SK$ as in \eqref{eqn:SK}.
For
$\omega = (\alpha, b, a,\beta) \in \SK,
$
we set 
$$
\alpha = t(\omega) \in \Sigma^\rho,\quad
 b=r(\omega)\in \Sigma^\eta,\quad
  a = l(\omega)\in \Sigma^\eta,\quad 
  \beta= b(\omega)\in \Sigma^\rho,\quad
$$
which stand for:
top, right, left, bottom respectively
as in the following figure.
$$
\begin{CD}
\cdot @>\alpha=t(\omega)>> \cdot \\
@V{a=l(\omega)}VV  @VV{b=r(\omega)}V \\
\cdot @>>\beta= b(\omega)> \cdot 
\end{CD}
$$
Define $*$-homomorphisms 
$\widehat{\rho}_\alpha : \BR \longrightarrow \A
$   
for $\alpha \in \Sigma^\rho$
and
$
\widehat{\eta}_a
: \BE \longrightarrow \A
$
for
$ a \in \Sigma^\eta$
by  \eqref{eqn:homo}
which satisfy  the equalities
\begin{equation*}
\widehat{\rho}_\alpha (w)
= P_\alpha w_\alpha  P_\alpha
\quad
\text{ and }
\quad
\widehat{\eta}_a (z)
= Q_a z_a  Q_a
\end{equation*}
for $w=\sum_{\beta \in \Sigma^\rho} S_\beta w_\beta S_\beta^* \in \BR$ 
as in \eqref{eqn:w} 
and $z=\sum_{b \in \Sigma^\eta} T_b z_b T_b^* \in \BE$ 
as in \eqref{eqn:z}. 
Their restrictions to $\A$
coincide with $\rho_\alpha$ and $\eta_a$
respectively. 
\begin{lem} 
Keep the above notations.
 \begin{enumerate}
\renewcommand{\labelenumi}{(\roman{enumi})}
\item 
For $\alpha \in \Sigma^\rho$
and
$z = \sum_{b \in \Sigma^\eta} T_b z_b T_b^* \in \BE$ 
as in \eqref{eqn:z},
put
\begin{equation}
\widehat{\rho}_\alpha^\eta(z) =
\sum
\begin{Sb}
b,a,\beta\\
(\alpha,b,a,\beta)\in\SK
\end{Sb}
T_b \rho_\beta(z_a) T_b^* \in \BE.
\end{equation}
Then 
$\widehat{\rho}_\alpha^\eta:\BE \longrightarrow \BE$
is a $*$-homomorphism such that 
$\widehat{\rho}_\alpha^\eta(y) = \rho_\alpha(y)$ for $y \in \A$.
\item 
For $a \in \Sigma^\eta$
and
$w = \sum_{\beta \in \Sigma^\rho} S_\beta w_\beta S_\beta^* \in \BR$ 
as in \eqref{eqn:w},
put
\begin{equation}
\widehat{\eta}_a^\rho(w) =
\sum
\begin{Sb}
\alpha, b,\beta\\
(\alpha,b,a,\beta)\in\SK
\end{Sb}
S_\beta \eta_b(w_\alpha) S_\beta^* \in \BR.
\end{equation}
Then 
$\widehat{\eta}_a^\rho:\BR \longrightarrow \BR$
is a $*$-homomorphism such that 
$\widehat{\eta}_a^\rho(y) = \eta_a(y)$ for $y \in \A$.
\end{enumerate}
\end{lem}
\begin{pf}
(i)
Since $z_a = \widehat{\eta}_a(z)$, the equality (2.4) goes to
\begin{equation*}
\widehat{\rho}_\alpha^\eta(z) =
\sum
\begin{Sb}
b,a,\beta\\
(\alpha,b,a,\beta)\in\SK
\end{Sb}
T_b \rho_\beta(\widehat{\eta}_a(z)) T_b^*
\end{equation*}
as in (1.7).
It is easy to see that
$\widehat{\rho}_\alpha^\eta: \BE \longrightarrow \BE$
yields a $*$-homomorphism.
If in particular
$z = y \in \A$,
we have
$y = \sum_{b \in \Sigma^\eta}
T_b \eta_b(y) T_b^*
$ 
so that
\begin{equation*}
\widehat{\rho}_\alpha^\eta(y)
 =
\sum
\begin{Sb}
b, a,\beta\\
(\alpha,b,a,\beta)\in\SK
\end{Sb}
T_b \rho_\beta(\eta_a(y) )T_b^* 
=
\sum
\begin{Sb}
a,b,\beta\\
(\alpha,b,a,\beta)\in\SK
\end{Sb}
T_b \eta_b(\rho_\alpha(y) )T_b^* 
= \rho_\alpha(y).
\end{equation*}

(ii) is similar to (i).
\end{pf}
The commutation relations \eqref{eqn:kappa} on $\A$ 
extend to $\BR$ and to $\BE$ as in the following lemma.
\begin{lem}
For $\omega = (\alpha, b, a,\beta) \in \Sigma_\kappa$,
we have
 \begin{enumerate}
\renewcommand{\labelenumi}{(\roman{enumi})}
\item 
$\eta_b \circ \widehat{\rho}_\alpha(w) 
=\widehat{\rho}_\beta\circ \widehat{\eta}^\rho_a(w)$ for $w \in \BR$.
\item 
$
\rho_\beta\circ \widehat{\eta}_a(z)
=
\widehat{\eta}_b \circ \widehat{\rho}^\eta_\alpha(z) $ for $z \in \BE$.
\end{enumerate}
\end{lem}
\begin{pf}
(i)
For 
$w = \sum_{\alpha' \in \Sigma^\rho} S_{\alpha'}w_{\alpha'} S_{\alpha'}^*$
as in \eqref{eqn:w},
we have
$S_\beta^* \widehat{\eta}^\rho_a(w) S_\beta 
= S_\beta^*S_\beta \eta_b(w_\alpha)  S_\beta^*S_\beta
$ 
so that by \eqref{eqn:kappa}
\begin{align*}
\widehat{\rho}_\beta\circ \widehat{\eta}^\rho_a(w)
& = P_\beta \eta_b(w_\alpha) P_\beta \\
& = \rho_\beta(1)  \eta_b(\rho_\alpha(1))
 \eta_b(w_\alpha)  \eta_b(\rho_\alpha(1)) \rho_\beta(1) \\
& = \rho_\beta(\eta_a(1)) \eta_b(w_\alpha) \rho_\beta(\eta_a(1)) \\
& = \eta_b(\rho_\alpha(1) w_\alpha \rho_\alpha(1)) 
 = \eta_b(\widehat{\rho}_\alpha(w) ).
\end{align*}

(ii) is similar to (i).
\end{pf}

\section{Hilbert $C^*$-quad modules from $C^*$-textile dynamical systems}
We fix a $C^*$-textile dynamical system $\CTDS$.
For $\omega = (\alpha, b, a,\beta) \in \Sigma_\kappa$,
we put a projection
\begin{equation}
E_\omega = \eta_b(\rho_\alpha(1)) (= \rho_\beta(\eta_a(1)) \in Z_\A.
\end{equation}
Let
$e_\omega, \omega \in \Sigma_\kappa$
denote the orthogonal basis of the vector space 
${\Bbb C}^{| \Sigma_\kappa |}$
where $|\SK|$ means the cardinal number of the finite set 
$\SK$. 
Define the vector space
\begin{equation}
\HK = \sum_{\omega \in \Sigma_\kappa} {\Bbb C}e_\omega \otimes E_\omega \A
\end{equation}
which is naturally isomorphic to the vector space
$
\sum_{\omega \in \Sigma_\kappa}{}^\oplus E_\omega \A.
$

We first endow $\HK$ with right $\A$-module structure and 
$\A$-valued inner product as follows:
For 
$\xi = \sum_{\omega\in \SK} e_\omega \otimes E_\omega x_\omega,
 \xi' =\sum_{\omega'\in \SK} e_{\omega'} \otimes E_{\omega'} x'_{\omega'}$
with $x_\omega, x'_{\omega'} \in \A$ 
and $y \in \A$,
set
\begin{align*}
\xi \varphi_\A(y) & := \sum_{\omega \in \SK}e_\omega \otimes E_\omega x_\omega y \in \HK, \\ 
\langle \xi \mid \xi' \rangle_{\A} 
& := \sum_{\omega \in \SK}
x_\omega^* E_\omega x_\omega'\in \A
\end{align*}
which satisfy the relations:
\begin{equation*}
\langle \xi \mid \xi'\varphi_\A(y) \rangle_{\A}
= \langle \xi \mid \xi' \rangle_{\A}\cdot y,
\qquad
\langle \xi \mid \xi' \rangle_{\A}^*
=
\langle \xi' \mid \xi \rangle_{\A}.
\end{equation*}
We will further endow $\HK$  with two other Hilbert $C^*$-bimodule 
structure.
Such a system will be called a Hilbert $C^*$-quad module.
Let
$\xi = \sum_{\omega \in \SK} e_\omega \otimes E_\omega x_\omega, 
\xi' = \sum_{\omega' \in \SK}
e_{\omega'} \otimes E_{\omega'} x'_{\omega'} \in \HK$ 
with 
$x_\omega, x'_{\omega'} \in \A$,
and
$w =\sum_{\alpha'\in \Sigma^\rho} 
S_{\alpha'} w_{\alpha'} S_{\alpha'}^* \in \BR
$
as in \eqref{eqn:w},
$z =\sum_{a'\in \Sigma^\eta} T_{a'} z_{a'} T_{a'}^* \in \BE$
as in \eqref{eqn:z}.
We define

1. The right $\BR$-action $\varphi_\rho$ and the right $\BE$-action $\varphi_\eta$:
\begin{equation*}
\xi 
\varphi_\rho(w)
  := \sum_{\omega \in \SK}
 e_\omega \otimes E_\omega x_\omega w_{b(\omega)}, 
\qquad
\xi
\varphi_\eta(z)
  :=  \sum_{\omega \in \SK}
 e_\omega \otimes E_\omega x_\omega z_{r(\omega)}.
\end{equation*}

2. The left $\BR$-action $\phi_\rho$ and the left $\BE$-action $\phi_\eta$:
\begin{equation*}
\phi_\rho(w) \xi 
 :=  \sum_{\omega \in \SK}
e_\omega \otimes E_\omega \eta_{r(\omega)}(w_{t(\omega)}) x_\omega,
\qquad
\phi_\eta( z) \xi 
 :=  \sum_{\omega \in \SK}
e_\omega \otimes E_\omega \rho_{b(\omega)}(z_{l(\omega)}) x_\omega.
\end{equation*}

3. The  
 right $\BR$-valued inner product 
$\langle \cdot | \cdot \rangle_{\rho}$
and the
right $\BE$-valued inner product 
$\langle \cdot | \cdot \rangle_{\eta}$:
\begin{equation*}
\langle \xi | \xi' \rangle_{\rho}
: = \sum_{\omega \in \SK}
S_{b(\omega)} x_\omega^* E_\omega x'_\omega S_{b(\omega)}^* , 
\qquad
\langle \xi | \xi' \rangle_{\eta}
: = 
\sum_{\omega \in \SK}
T_{r(\omega)} x_\omega^* E_\omega x'_\omega T_{r(\omega)}^*.
\end{equation*}
The following lemma is straightforward. 
\begin{lem}
For $\xi \in \HK$ and $w, w' \in \BR$, $z, z' \in \BE$, we have
\begin{equation*}
\begin{aligned}
(\xi \varphi_\rho(w)) \varphi_\rho(w')
& = \xi \varphi_\rho(w w'), \\ 
 \phi_\rho(w)( \phi_\rho(w') \xi)
& = \phi_\rho(w w') \xi, \\
 \phi_\rho(w)( \xi \varphi_\rho(w')) 
& = (\phi_\rho(w) \xi) \varphi_\rho(w'),
\end{aligned}
\qquad
\begin{aligned}
(\xi \varphi_\eta(z)) \varphi_\eta(z')
& = \xi \varphi_\eta(z z'), \\ 
 \phi_\eta(z)( \phi_\eta(z') \xi)
 & = \phi_\eta(z z') \xi, \\ 
 \phi_\eta(z)( \xi \varphi_\eta(z')) 
& = (\phi_\eta(z) \xi) \varphi_\eta(z').
\end{aligned}
\end{equation*}
\end{lem}
\begin{lem}
For $\xi, \xi' \in \HK$ and $w \in \BR$, $z \in \BE$, we have
\begin{equation*}
\begin{aligned}
\langle \xi \mid \xi' \varphi_\rho(w) \rangle_\rho 
& = \langle \xi \mid \xi'  \rangle_\rho \cdot w, \\ 
\langle \xi \varphi_\rho(w) \mid \xi'  \rangle_\rho 
& = w^* \cdot \langle \xi \mid \xi'  \rangle_\rho, \\
\langle \phi_\rho(w)\xi \mid \xi' \rangle_\rho
& = 
\langle \xi \mid \phi_\rho(w^*) \xi' \rangle_\rho,
\end{aligned}
\qquad
\begin{aligned}
\langle \xi \mid \xi' \varphi_\eta(z) \rangle_\eta 
& = \langle \xi \mid \xi'  \rangle_\eta \cdot z, \\ 
\langle \xi \varphi_\eta(z) \mid \xi'  \rangle_\eta 
& = z^* \cdot \langle \xi \mid \xi'  \rangle_\eta, \\
\langle \phi_\eta(z)\xi \mid \xi' \rangle_\eta
& = 
\langle \xi \mid \phi_\eta(z^*) \xi' \rangle_\eta.
\end{aligned}
\end{equation*}
\end{lem}
\begin{pf}
We will show the equalities
\begin{equation*}
\langle \xi \mid \xi' \varphi_\rho(w) \rangle_\rho 
= \langle \xi \mid \xi'  \rangle_\rho \cdot w, 
\qquad
\langle \phi_\rho(w)\xi \mid \xi' \rangle_\rho
= 
\langle \xi \mid \phi_\rho(w^*) \xi' \rangle_\rho.
\end{equation*}
For 
$
\xi= \sum_{\omega \in \SK} e_\omega \otimes E_\omega x_\omega,
\xi' =\sum_{\omega' \in \SK} e_{\omega'} \otimes E_{\omega'} x'_{\omega'} \in \HK 
$ 
with 
$x_\omega, x'_{\omega'} \in \A$,
and
$w =\sum_{\gamma\in \Sigma^\eta} S_{\gamma} w_{\gamma} S_{\gamma}^* \in \BR$
as in \eqref{eqn:w},
we have
\begin{align*}
\langle \xi \mid \xi' \varphi_\rho(w) \rangle_\rho 
& = \sum_{\omega, \omega' \in \SK} \langle  e_\omega \otimes E_\omega x_\omega \mid 
    e_{\omega'} \otimes E_{\omega'} x'_{\omega'}  w_{b(\omega')} \rangle_\rho \\
& =
\sum_{\omega \in \SK}
S_{b(\omega)} x_\omega^* E_\omega x'_{\omega} w_{b(\omega)} S_{b(\omega)}^* \\ 
&
=
(\sum_{\omega \in \SK} S_{b(\omega)} x_\omega^* E_\omega x'_{\omega}  S_{b(\omega)} ^* )\cdot 
(\sum_{\gamma\in \Sigma^\rho} S_{\gamma} w_{\gamma} S_{\gamma}^*) \\
& 
=
\langle \xi \mid \xi'  \rangle_\rho \cdot w.
\end{align*}
We also have
\begin{align*}
\langle \phi_\rho(w) \xi \mid \xi' \rangle_\rho
& = \sum_{\omega, \omega' \in \SK} \sum_{\gamma \in \Sigma^\rho}
\langle \phi_\rho(S_{\gamma} w_{\gamma} S_{\gamma}^*)
(e_\omega \otimes E_\omega x_\omega)
 \mid e_{\omega'} \otimes E_{\omega'}x'_{\omega'} \rangle_\rho \\
& =  \sum_{\omega, \omega' \in \SK}
\langle  e_\omega \otimes E_\omega \eta_{r(\omega)}( w_{t(\omega)}) x_\omega
 \mid e_{\omega'} \otimes E_{\omega'}x'_{\omega'} \rangle_\rho \\
& = 
\sum_{\omega \in \SK}
S_{b(\omega)}  x_\omega^* E_\omega 
\eta_{r(\omega)} (w_{t(\omega)}^*)x'_{\omega}S_{b(\omega)}^* \\
& =\sum_{\omega, \omega' \in \SK} 
\langle  e_\omega \otimes E_\omega  x_\omega
 \mid e_{\omega'} \otimes E_{\omega'} \eta_{r(\omega')}(w_{t(\omega')}^*)
x'_{\omega'} \rangle_\rho \\
& = \sum_{\omega, \omega' \in \SK} 
\sum_{\gamma \in \Sigma^\rho}
\langle  e_\omega \otimes E_\omega x_\omega
 \mid \phi_\rho(S_{\gamma} w_{\gamma}^* S_{\gamma}^*)
 (e_{\omega'} \otimes E_{\omega'}x'_{\omega'})
 \rangle_\rho \\
& =
\langle  \xi \mid \phi_\rho(w^*) \xi' \rangle_\rho.
\end{align*}
The  three equalities of the right hand side are similarly shown to 
the above equalities. 
\end{pf}
Hence we have
\begin{align*}
\phi_\rho(w^*) =\phi_\rho(w)^* : & \text{ the adjoint with respect to the inner product } \langle \cdot \mid \cdot \rangle_\rho, \\
\phi_\eta(z^*) =\phi_\eta(z)^* : & \text{ the adjoint with respect to the inner product } \langle \cdot \mid \cdot \rangle_\eta. 
\end{align*}
The following lemma is direct 
and shows that the two module structure are compatible to each other. 
\begin{lem} 
For $w \in \BR, z \in \BE$ and $\xi \in \HK$, we have
\begin{enumerate}
\renewcommand{\labelenumi}{(\roman{enumi})}
\item $(\phi_\rho(w)\xi)\varphi_\eta(z) = \phi_\rho(w)(\xi \varphi_\eta(z)). $  
\item $(\phi_\eta(z)\xi)\varphi_\rho(w) = \phi_\eta(z)(\xi \varphi_\rho(w)). $
\end{enumerate}
\end{lem}
Then we have the following proposition
\begin{prop}
Keep the above notations.
\begin{enumerate}
\renewcommand{\labelenumi}{(\roman{enumi})}
\item
$(\HK,\varphi_\rho)$ is a right  $\BR$-module
with right $\BR$-valued inner product
$\langle \cdot \mid \cdot \rangle_\rho$ and
left $\BR$-action by $\phi_\rho$.
Hence 
$\HK$ is a Hilbert $C^*$-bimodule over $\BR$.
\item
$(\HK,\varphi_\eta)$ is a right  $\BE$-module
with right $\BE$-valued inner product
$\langle \cdot \mid \cdot \rangle_\eta$ and
left $\BE$-action by $\phi_\eta$.
Hence 
$\HK$ is a Hilbert $C^*$-bimodule over $\BE$.
\end{enumerate}
Therefore $\HK$ has multi structure of Hilbert $C^*$-bimodules,
which are compatible to each other.
\end{prop}
$\HK$ is originally a Hilbert $C^*$-right module over $\A$,
which is also compatible to the two left actions $\phi_\rho$ of $\BR$ and $\phi_\eta$ of $\BE$ 
as in the following lemma.
Its proof is straightforward.
\begin{lem}\label{lem:4.3}
For $\xi \in \HK$ and  $ y \in \A$, we have
 \begin{enumerate}
\renewcommand{\labelenumi}{(\roman{enumi})}
\item 
$\phi_\rho(w) (\xi \varphi_\A(y)) = (\phi_\rho(w) \xi) \varphi_\A(y)  $ 
for $w \in \BR$.
\item 
$\phi_\eta(z) (\xi \varphi_\A(y)) = (\phi_\eta(z) \xi) \varphi_\A(y)  $ 
for $z \in \BE$.
\end{enumerate}
Hence both
$\phi_\eta(z)$ and $\phi_\rho(w)$ are right $\A$-module maps.
\end{lem}
%
Define positive maps
$\psi_\rho : \A \longrightarrow \BR$  
and
$\psi_\eta : \A \longrightarrow \BE$ 
by
\begin{equation}
\psi_\rho(y)  = \sum_{\alpha\in \Sigma^\rho} 
S_{\alpha} y S_{\alpha}^* \in \BR ,\qquad 
\psi_\eta(y)  = \sum_{a\in \Sigma^\eta} 
T_{a} y T_{a}^* \in \BE 
\end{equation}
for $y \in \A.$
 Then we have
\begin{lem} \label{lem:4.4}
For $\xi \in \HK$ and  $ y \in \A$, we have
 \begin{enumerate}
\renewcommand{\labelenumi}{(\roman{enumi})}
\item 
$\xi\varphi_\rho(w \psi_\rho(y)) = ( \xi\varphi_\rho(w)) \varphi_\A(y)$ 
for $w \in \BR$.
\item 
$\xi \varphi_\eta(z \psi_\eta(y)) = ( \xi \varphi_\eta(z)) \varphi_\A(y)$ 
for $z \in \BE$.
\end{enumerate}
Hence we have
\begin{equation*}
\xi \varphi_\rho(\psi_\rho(y)) =\xi \varphi_\eta( \psi_\eta(y)) =  \xi  \varphi_\A(y).
\end{equation*}
\end{lem}
\begin{pf}
(i)
For 
$\xi = \sum_{\omega \in \SK} e_\omega \otimes E_\omega x_\omega \in \HK$ 
with 
$x_\omega \in \A$,
and
$w =\sum_{\alpha'\in \Sigma^\rho}
S_{\alpha'} w_{\alpha'} S_{\alpha'}^* \in \BR$
as in \eqref{eqn:w},
we have
\begin{equation*}
w \psi_\rho(y) 
= \sum_{\alpha'\in \Sigma^\rho}
S_{\alpha'} w_{\alpha'} S_{\alpha'}^* 
\sum_{\beta'\in \Sigma^\rho} 
S_{\beta'} y S_{\beta'}^*
= \sum_{\alpha'\in \Sigma^\rho} 
S_{\alpha'} w_{\alpha'} y S_{\alpha'}^* 
\end{equation*}
so that
\begin{align*}
\xi \varphi_\rho( w\psi_\rho(y)) 
& = \sum_{\omega \in \SK} e_\omega \otimes E_\omega x_\omega  w_{b(\omega)} y \\
& = \sum_{\omega \in \SK} [(e_\omega \otimes E_\omega x_\omega ) \varphi_\rho(w)]
      \varphi_\A(y) 
 = [\xi \varphi_\rho(w) ] \varphi_\A(y).
\end{align*}

(ii) is similar to (i).
\end{pf}
\begin{lem} \label{lem:4.7}
For 
$y \in \A$ and $\xi \in \HK$,
we have 
$\phi_\rho(y) \xi = \phi_{l}(y) \xi$. 
\end{lem}
\begin{pf}
Since the identities 
$y 
 = \sum_{\alpha \in \Sigma^\rho} S_\alpha \rho_\alpha(y) S_\alpha^*
 = \sum_{a \in \Sigma^\eta} T_a\eta_a(y)T_a^*
$ 
hold, 
we have
for $\xi = \sum_{\omega \in \SK} e_\omega \otimes E_\omega x_\omega$ 
with $x_\omega \in \A$,
\begin{equation*}
\phi_\rho(y) \xi 
 = \sum_{\omega \in \SK} e_\omega \otimes E_\omega 
    \eta_{r(\omega)}(\rho_{t(\omega)}(y)) x_\omega. 
\end{equation*}
On the other hand,
we have
\begin{equation*}
\phi_\eta(y) \xi 
 = \sum_{\omega \in \SK} e_\omega \otimes E_\omega 
    \rho_{b(\omega)}(\eta_{l(\omega)}(y)) x_\omega. 
\end{equation*}
As
$\eta_{r(\omega)}(\rho_{t(\omega)}(y))
=\rho_{b(\omega)}(\eta_{l(\omega)}(y))$,
we obtain the desired equality.
\end{pf}
By the above lemma, 
we may define the left action $\phi$ of $\A$ on $\HK$ by
\begin{equation*}
\phi(y) \xi :=\phi_\rho(y) \xi = \phi_{l}(y) \xi, \qquad
y \in \A, \, \xi \in \HK
\end{equation*} 
so that 
$\HK$ has a structure of a Hilbert $C^*$-bimodule over $\A$.
We note the following lemma.
\begin{lem}
If the algebra $\A$ is commutative, we have
\begin{equation*}
 \phi_\rho(w)\phi_\eta(z) = \phi_\eta(z) \phi_\rho(w),
\qquad w \in \BR, z \in \BE.
\end{equation*}
\end{lem}
\begin{pf}
For 
$
w = \sum_{\alpha \in \Sigma^\rho}S_\alpha w_\alpha S_\alpha^*$ 
as in \eqref{eqn:w},
$
z = \sum_{a' \in \Sigma^\eta} T_a z_a T_a^*
$ 
as in \eqref{eqn:z}
and
$\xi = \sum_{\omega \in \SK} e_\omega \otimes E_\omega x_\omega$
with $x_\omega \in \A$,
we have
\begin{align*}
\phi_\rho(w)\phi_\eta(z) \xi
& = \sum_{\omega \in \SK}
     e_\omega \otimes E_\omega 
     \eta_{r(\omega)}(w_{t(\omega)})\rho_{b(\omega)}(z_{l(\omega)}) x_\omega \\
& = \sum_{\omega \in \SK} 
      e_\omega \otimes E_\omega 
      \rho_{b(\omega)}(z_{l(\omega)})\eta_{r(\omega)}(w_{t(\omega)}) x_\omega \\
& = \phi_\eta(z) \phi_\rho(w) \xi.
\end{align*}
\end{pf}
Put
for $\alpha \in \Sigma^\rho$ and $a \in \Sigma^\eta$
\begin{equation*}
u_\alpha = \sum_{\omega \in \SK, \alpha =t(\omega)} e_\omega\otimes E_\omega \in \HK,
\qquad 
v_a = \sum_{\omega \in \SK, a =l(\omega)} e_\omega\otimes E_\omega \in \HK.
\end{equation*}
\begin{lem} \label{lem:4.8}
Keep the above notations.
 \begin{enumerate}
\renewcommand{\labelenumi}{(\roman{enumi})}
\item $\{ u_\alpha \}_{\alpha \in \Sigma^\rho}$ 
forms an essential orthogonal finite basis of $\HK$ 
with respect to the $\BE$-valued inner product 
$\langle \cdot \mid \cdot \rangle_\eta$
 as right $\BE$-module through $\varphi_\eta$.
\item $\{ v_a \}_{a \in \Sigma^\eta}$ forms an essential orthogonal finite basis of $\HK$ 
with respect to the $\BR$-valued inner product 
$\langle \cdot \mid \cdot \rangle_\rho$
as right $\BR$-module through $\varphi_\rho$.
\end{enumerate}
\end{lem}
\begin{pf}
(i)
For $\alpha, \beta \in \Sigma^\rho$,
we have
\begin{align*}
\langle u_\alpha \mid u_\beta \rangle_\eta
& =  \langle \sum_{\omega \in \SK, \alpha =t(\omega)} e_\omega\otimes E_\omega
\mid         \sum_{\omega' \in \SK, \beta =t(\omega')} e_{\omega'} \otimes E_{\omega'} \rangle_\eta \\
& =
{\begin{cases}
\sum_{\omega \in \SK, \alpha =t(\omega)} 
\langle e_\omega \otimes E_\omega \mid e_\omega \otimes E_\omega \rangle_\eta
  & \text{ if } \alpha = \beta, \\
0 & \text{ if } \alpha \ne \beta.
\end{cases}} 
\end{align*}
Since 
\begin{align*}
 \sum_{\omega \in \SK, \alpha =t(\omega)} 
  \langle e_\omega \otimes E_\omega \mid e_\omega \otimes E_\omega \rangle_\eta 
=&
\sum_{\omega \in \SK, \alpha =t(\omega)} 
T_{r(\omega)} E_\omega T_{r(\omega)}^* \\
=& 
\sum_{\omega \in \SK, \alpha =t(\omega)} 
T_{r(\omega)} T_{r(\omega)}^*  \rho_\alpha(1) T_{r(\omega)} T_{r(\omega)}^*
= P_\alpha,
\end{align*}
we see 
\begin{equation*}
\langle u_\alpha \mid u_\beta \rangle_\eta
= 
\begin{cases}
P_\alpha & \text{ if } \alpha = \beta, \\
0 & \text{ if } \alpha \ne \beta. \\
\end{cases}
\end{equation*}
Hence we have
\begin{equation*}
\sum_{\alpha \in \Sigma^\rho}
\langle u_\alpha \mid u_\alpha \rangle_\eta
= \sum_{\alpha \in \Sigma^\rho} P_\alpha \ge 1.
\end{equation*}
For 
$\xi = \sum_{\omega \in \SK}e_\omega\otimes E_\omega x_\omega\in \HK
$
with
$x_\omega \in \A$,
we have
\begin{equation*}
\langle u_\alpha \mid \xi \rangle_\eta
 =
\sum_{\omega \in \SK, \alpha =t(\omega)} 
T_{r(\omega)} E_\omega x_\omega T_{r(\omega)}^*.
\end{equation*}
It then follows that
\begin{align*}
u_\alpha
\varphi_\eta(\langle u_\alpha \mid \xi \rangle_\eta)
& =
(\sum_{\omega \in \SK, \alpha =t(\omega)} e_\omega\otimes E_\omega)
\varphi_\eta(
\sum_{\omega' \in \SK, \alpha =t(\omega')} 
T_{r(\omega')} E_{\omega'} x_{\omega'} T_{r(\omega')}^*)\\
& =
\sum_{\omega, \omega' \in \SK, \alpha =t(\omega)=t(\omega')} 
(e_\omega\otimes E_\omega)
\varphi_\eta(
T_{r(\omega')} E_{\omega'} x_{\omega'} T_{r(\omega')}^*)\\
& =
\sum_{\omega \in \SK, \alpha =t(\omega)} e_\omega\otimes E_\omega x_\omega
\end{align*}
so that
\begin{equation*}
\sum_{\alpha \in \Sigma^\rho}
u_\alpha
\varphi_\eta(\langle u_\alpha \mid \xi \rangle_\eta)
=\sum_{\alpha \in \Sigma^\rho}
\sum_{\omega \in \SK, \alpha =t(\omega)} e_\omega\otimes E_\omega x_\omega 
 =\sum_{\omega \in \HK} e_\omega \otimes E_\omega x_\omega = \xi.
\end{equation*}

(ii) is similar to (i).
\end{pf}
As 
$\langle u_\alpha \mid u_\alpha \rangle_\eta = P_\alpha$
and
$\langle v_a \mid v_a \rangle_\rho = Q_a$,
we note that the equality
\begin{equation*}
\eta_b(\langle u_\alpha \mid u_\alpha \rangle_\eta)
=
\rho_\beta(\langle v_a \mid v_a \rangle_\rho)
= E_\omega
\end{equation*}
holds for $\omega =(\alpha, b, a, \beta)\in \SK$.
\begin{lem}
For $\alpha \in \Sigma^\rho, a \in \Sigma^\eta$ and $y \in \A$, 
we have
\begin{enumerate}
\renewcommand{\labelenumi}{(\roman{enumi})}
\item $\phi(y) u_\alpha = u_\alpha \varphi_\eta(\rho_\alpha(y))$
and hence
$\rho_\alpha(y) = \langle u_\alpha \mid \phi(y) u_\alpha\rangle_\eta$.
\item $\phi(y) v_a = v_a \varphi_\rho(\eta_a(y))$
and hence
$\eta_a(y) = \langle v_a \mid \phi(y) v_a\rangle_\rho$.
\end{enumerate}
Therefore the commutation relations \eqref{eqn:kappa} are rephrased as the equality
\begin{equation}
\langle v_b \mid 
\phi(\langle u_\alpha \mid \phi(y)u_\alpha \rangle_\eta )
v_b \rangle_\rho
=
\langle
u_\beta\mid
\phi(\langle v_a \mid \phi(y)v_a \rangle_\rho)
u_\beta \rangle_\eta 
\end{equation}
for $\omega = (\alpha, b, a, \beta) \in \SK$ and $y \in \A$.
\end{lem}
\begin{pf}
(i)  It follows that
\begin{align*}
\phi(y) u_\alpha
 =\phi_\eta(y) u_\alpha
&= \phi_\eta(\sum_{a' \in \Sigma^\eta}T_{a'}\eta_a(y)T_{a'}^* )
   (\sum_{\omega \in \SK, \alpha =t(\omega)} e_\omega\otimes E_\omega)\\
& = \sum_{\omega \in \SK, \alpha =t(\omega)}
        e_\omega \otimes E_\omega \rho_{b(\omega)}(\eta_{l(\omega)}(y))  \\
& = \sum_{\omega \in \SK, \alpha =t(\omega)}
        e_\omega \otimes E_\omega \eta_{r(\omega)} (\rho_\alpha(y))  \\
& = \sum_{\omega \in \SK, \alpha =t(\omega)}
    [(e_\omega \otimes E_\omega)
     \varphi_\eta(\sum_{b' \in \Sigma^\eta}T_{b'}\eta_{b'}(\rho_\alpha(y))T_{b'}^* )]  \\
& = u_\alpha \varphi_\eta(\rho_\alpha(y)). 
\end{align*}
It then follows that
\begin{equation*}
\langle u_\alpha \mid \phi(y) u_\alpha\rangle_\eta
 = \langle u_\alpha \mid u_\alpha\varphi_\eta(\rho_\alpha(y)) \rangle_\eta 
  = \langle u_\alpha \mid u_\alpha  \rangle_\eta\cdot \rho_\alpha(y) 
   = \rho_\alpha(y). 
\end{equation*}

(ii) is similar to (i).
\end{pf}
More generally we have
 \begin{lem}
For $\alpha \in \Sigma^\rho$,
$a \in \Sigma^\eta$
and
$w \in \BR, z \in \BE$,
we have
\begin{enumerate}
\renewcommand{\labelenumi}{(\roman{enumi})}
\item 
$\phi_\rho(w) u_\alpha
 = u_\alpha \varphi_\eta(\langle u_\alpha \mid \phi_\rho(w) u_\alpha \rangle_\eta),
 \quad
 \phi_\rho(w) v_a
 = v_a \varphi_\rho(\langle v_a \mid \phi_\rho(w)v_a \rangle_\rho).$
\item 
$\phi_\eta(z) u_\alpha
 =  u_\alpha \varphi_\eta(\langle u_\alpha \mid \phi_\eta(z) u_\alpha \rangle_\eta),
 \quad
 \phi_\eta(z) v_a
 =  v_a \varphi_\rho(\langle v_a \mid \phi_\eta(z) v_a \rangle_\rho).$
\end{enumerate}
\end{lem}
\begin{pf}
(i)
For 
$\xi\in \HK$,
we have
$
\xi = \sum_{\alpha' \in \Sigma^\rho} u_{\alpha'} 
\varphi_\eta(\langle u_{\alpha'} \mid \xi \rangle_\eta).
$
For $\alpha \ne \alpha'$, 
we have
$\langle u_\alpha \mid \phi_\rho(w)u_{\alpha'} \rangle_\eta = 0$
so that
\begin{equation*}
 \phi_\rho(w) u_\alpha 
=
\sum_{\alpha' \in \Sigma^\rho}
 u_{\alpha'} \varphi_\eta(\langle u_{\alpha'} 
\mid \phi_\rho(w) u_\alpha \rangle_\eta)
 =  u_\alpha \varphi_\eta(\langle u_\alpha \mid \phi_\rho(w) u_\alpha \rangle_\eta).
\end{equation*}
Similarly 
for $a \ne a'$,
we have
$\langle v_{a'} \mid \phi_\rho(w)v_a \rangle_\rho = 0$
so that
\begin{equation*}
 \phi_\rho(w) v_a 
=
\sum_{a' \in \Sigma^\eta}
 v_{a'} \varphi_\rho(\langle v_a'  
\mid \phi_\rho(w) v_a \rangle_\rho)
 =  v_a \varphi_\rho(\langle v_a \mid \phi_\rho(w) v_a \rangle_\rho).
\end{equation*}

 (ii) is similar to (i).
\end{pf}
The following lemma states that the $*$-homomorphisms
$
\hat{\rho}_\alpha,
 \hat{\eta}_a^\rho$ 
on $\BR$
and
$\hat{\eta}_a,
\hat{\rho}_\alpha^\eta
$ on $\BE$ are given by inner products.
\begin{lem}
For $\alpha \in \Sigma^\rho$,
$a \in \Sigma^\eta$
and
$w \in \BR, z \in \BE$,
we have
\begin{enumerate}
\renewcommand{\labelenumi}{(\roman{enumi})}
\item 
$\hat{\rho}_\alpha(w) = \langle u_\alpha \mid \phi_\rho(w)u_\alpha\rangle_\eta, \quad
\hat{\eta}_a(z) = \langle v_a \mid \phi_\eta(z)v_a\rangle_\rho$.
\item
$\hat{\rho}^\eta_\alpha(z)
 = \langle u_\alpha \mid \phi_\eta(z)u_\alpha\rangle_\eta, \quad
\hat{\eta}^\rho_a(w)
= \langle v_a \mid \phi_\rho(w)v_a\rangle_\rho$.
\end{enumerate}
\end{lem}
\begin{pf}
(i)
For $w = \sum_{\alpha' \in \Sigma^\rho}S_{\alpha'}w_{\alpha'} S_{\alpha'}^* \in \BR$ 
as in \eqref{eqn:w},
we have
\begin{equation*}
\phi_\rho(w) u_\alpha 
=  
\sum_{\omega\in \SK, t(\omega) = \alpha} e_\omega \otimes E_\omega \eta_{r(\omega)}(w_\alpha)
\end{equation*}
so that
\begin{align*}
\langle u_\alpha \mid \phi_\rho(w)u_\alpha\rangle_\eta
& = \langle \sum_{\omega'\in \SK, t(\omega') = \alpha} e_{\omega'}\otimes E_{\omega'} 
\mid 
\sum_{\omega\in \SK, t(\omega) = \alpha} e_\omega \otimes E_\omega \eta_{r(\omega)}(w_\alpha) \rangle_\eta \\
& =  
\sum_{\omega\in \SK, t(\omega) = \alpha} T_{r(\omega)} E_\omega \eta_{r(\omega)}(w_\alpha)  T_{r(\omega)}^* \\
& =  
\sum \begin{Sb} b \\ (\alpha,b) \in \Sigma^{\rho\eta} \end{Sb} 
T_b T_b^* w_\alpha T_bT_b^* \\
& =  
\sum_{b \in \Sigma^{\eta}} 
T_b T_b^* S_\alpha^*S_\alpha w_\alpha S_\alpha^*S_\alpha T_bT_b^* 
=  
P_\alpha w_\alpha P_\alpha
= \widehat{\rho}_\alpha(w).
\end{align*}
The other equality for
$\widehat{\eta}_a(z)$
is similarly shown to the above equalities.

(ii)
For $z = \sum_{a \in \Sigma^\eta}T_a z_a  T_a^* \in \BE$ 
as in \eqref{eqn:z},
we have
\begin{equation*}
\phi_\eta(z) u_\alpha 
=  
\sum_{\omega\in \SK, t(\omega) = \alpha}
e_\omega \otimes E_\omega \rho_{b(\omega)}(z_{l(\omega)}) 
\end{equation*}
so that
\begin{equation*}
\langle u_\alpha \mid \phi_\eta(z) u_\alpha \rangle_\eta
 =
\sum_{\omega\in \SK, t(\omega) = \alpha}
T_{r(\omega)} \rho_{b(\omega)}(z_{l(\omega ) } ) T_{r(\omega)}^* 
=\hat{\rho}^\eta_\alpha(z).
\end{equation*}
The other equality for
$\widehat{\eta}_a^\rho(w)$
is similarly shown to the above equalities.
\end{pf}

We will next study the norms on $\HK$ induced by the 
two inner products
$\langle \cdot \mid \cdot \rangle_\eta$
and
$\langle \cdot \mid \cdot \rangle_\rho$
\begin{lem}
For 
$
\xi =\sum_{\omega \in \SK}e_\omega\otimes E_\omega x_\omega\in \HK
$
with
$x_\omega \in \A$, we have
\begin{enumerate}
\renewcommand{\labelenumi}{(\roman{enumi})}
\item
$\| \langle \xi \mid \xi \rangle_\rho \|
= \max_{\beta \in \Sigma^\rho} 
\| \sum_{a \in \Sigma^\eta} x_{a,\beta}^* x_{a,\beta}\|,$
\item
$\| \langle \xi \mid \xi \rangle_\eta \|
= \max_{b \in \Sigma^\eta} 
\| \sum_{\alpha \in \Sigma^\rho} x_{\alpha,b}^* x_{\alpha,b}\|,$
\end{enumerate}
where
$E_\omega x_\omega  = x_{\alpha,b} = x_{a,\beta}$
for
$\omega = (\alpha, b, a, \beta) \in \SK.$
\end{lem}
\begin{pf}
(i) We have
\begin{align*}
\| \langle \xi  \mid \xi \rangle_\rho \|
& = \| \sum_{\omega \in \SK} 
S_{b(\omega)} x_\omega^* E_\omega x_\omega S_{b(\omega)}^* \| \\
& = \| \sum_{(a,\beta) \in \Sigma^{\eta\rho}} 
S_\beta x_{a,\beta}^* x_{a,\beta} S_\beta^* \| \\
& = \| \sum_{\beta\in \Sigma^\rho} 
S_\beta (\sum_{a\in \Sigma^\eta} x_{a,\beta}^* x_{a,\beta} )S_\beta^* \| \\
& =  \max_{\beta\in \Sigma^\rho} 
\| S_\beta (\sum_{a\in \Sigma^\eta} x_{a,\beta}^* x_{a,\beta} )S_\beta^* \|.
\end{align*}
Since 
$x_{a,\beta} = x_{a,\beta} P_\beta$, 
we have 
for $\beta \in \Sigma^\rho$
$$
\| S_\beta (\sum_{a\in \Sigma^\eta} x_{a,\beta}^* x_{a,\beta} )S_\beta^* \|
\le 
\| \sum_{a\in \Sigma^\eta} x_{a,\beta}^* x_{a,\beta} \|
\le
\| \sum_{a\in \Sigma^\eta} 
   P_\beta x_{a,\beta}^* x_{a,\beta} P_\beta \|
$$
one has
$$
\| S_\beta (\sum_{a\in \Sigma^\eta} x_{a,\beta}^* x_{a,\beta} )S_\beta^* \|
= 
\| \sum_{a\in \Sigma^\eta} x_{a,\beta}^* x_{a,\beta} \|.
$$
Therefore we have
$$
\| \langle \xi  \mid \xi \rangle_\rho \|
=  \max_{\beta\in \Sigma^\rho} 
\| \sum_{a\in \Sigma^\eta} x_{a,\beta}^* x_{a,\beta} \|.
$$
(ii) is similar to (i).
\end{pf}
Define positive maps
$\lambda_\rho : \BR \longrightarrow \A$ 
and 
$\lambda_\eta : \BE \longrightarrow \A$  by
\begin{equation}
\lambda_\rho(w)  = \sum_{\alpha\in \Sigma^\rho} 
\widehat{\rho}_\alpha(w), \quad w \in \BR,
\qquad
\lambda_\eta(z)  = \sum_{a\in \Sigma^\eta} 
\widehat{\eta}_a(z), \quad z \in \BE.
\end{equation}
Then we have for $\xi, \xi' \in \HK$
\begin{equation}
\langle \xi | \xi' \rangle_{\A}
= \lambda_\rho(
\langle \xi | \xi' \rangle_{\rho})
= \lambda_\eta(
\langle \xi | \xi' \rangle_{\eta}). \label{eqn:innerprod}
\end{equation}
Put
$
C_\rho  = \| \lambda_\rho(1) \|,
C_\eta  = \| \lambda_\eta(1) \|.
$
As 
$\lambda_\rho(1) = \sum_{\alpha \in \Sigma^\rho} \rho_\alpha(1) \ge 1$,
one sees $C_\rho \ge 1$ and similarly
$C_\eta \ge 1$.
Define the three norms for $\xi \in \HK$ 
\begin{equation}
\| \xi \|_{\A} = \| \langle \xi \mid \xi \rangle_\A \|^{\frac{1}{2}},
\quad
\| \xi \|_{\rho} = \| \langle \xi \mid \xi \rangle_\rho \|^{\frac{1}{2}},
\quad
\| \xi \|_{\eta} = \| \langle \xi \mid \xi \rangle_\eta \|^{\frac{1}{2}}.
\end{equation}
\begin{lem}\label{lem:norm}
The following inequalities hold for $\xi \in \HK$:
$$
\|  \xi \|_\rho 
\le \| \xi \|_\A 
\le
C_\rho^{\frac{1}{2}}
\|  \xi  \|_\rho
\quad 
\text{ and }
\quad
\|  \xi \|_\eta 
\le \| \xi \|_\A 
\le
C_\eta^{\frac{1}{2}}
\|  \xi  \|_\eta.
$$
Hence the three norms
$
\| \xi \|_\A, \,
\|  \xi  \|_\eta, \,
\|  \xi  \|_\eta
$
are equivalent to each other.
\end{lem}
\begin{pf}
For 
$
\xi =\sum_{\omega \in \SK}e_\omega\otimes E_\omega x_\omega\in \HK
$
with
$x_\omega \in \A$,
where
$E_\omega x_\omega = x_{\alpha,b} = x_{a,\beta}$
for
$\omega = (\alpha, b, a, \beta) \in \SK,$
we have
$$
\| \langle \xi \mid \xi \rangle_\rho \|
= \max_{\beta \in \Sigma^\rho} 
\| \sum_{a \in \Sigma^\eta} x_{a,\beta}^* x_{a,\beta}\|.
$$
We then have
\begin{equation*}
\| \xi \|_{\A} = \| \langle \xi \mid \xi \rangle_\A \|^{\frac{1}{2}}
=
\| \sum_{\omega \in \SK} x_{\omega}^* x_{\omega} \|^{\frac{1}{2}}
=\| \sum_{\beta \in \Sigma^\rho} \sum_{a \in \Sigma^\eta} 
x_{a,\beta}^* x_{a,\beta}\|^{\frac{1}{2}}.
\end{equation*}
Since
$$
\| \sum_{a \in \Sigma^\eta} x_{a,\beta}^* x_{a,\beta}\|
\le
\| \sum_{\beta \in \Sigma^\rho} \sum_{a \in \Sigma^\eta} 
x_{a,\beta}^* x_{a,\beta}\|,
$$  
we have
$$
\| \langle \xi \mid \xi \rangle_\rho \| \le \| \xi \|_\A.
$$ 
On the other hand,
by the equality
$\langle \xi \mid \xi \rangle_\A
=
\lambda_\rho(
\langle \xi \mid \xi \rangle_\rho)$,
we have
\begin{equation*}
\| \langle \xi \mid \xi \rangle_\A \|
\le
\| \lambda_\rho\|
\| \langle \xi \mid \xi \rangle_\rho \|
=
\| \lambda_\rho(1)\|
\| \langle \xi \mid \xi \rangle_\rho \|
= C_\rho
\| \langle \xi \mid \xi \rangle_\rho \|.
\end{equation*}
Therefore we have
$$
\| \xi \|_\rho 
\le \| \xi \|_\A 
\le
C_\rho^{\frac{1}{2}}
\| \xi \|_\rho
\quad
\text{ and similarly }
\quad
\| \xi \|_\eta 
\le \| \xi \|_\A 
\le
C_\rho^{\frac{1}{2}}
\| \xi \|_\eta.
$$
\end{pf}
\section{Fock  Hilbert $C^*$-quad modules and creation operators}
In this section, 
we will consider
relative tensor products of Hilbert $C^*$-quad modules
and
 introduce Fock space of Hilbert $C^*$-quad modules
which is two-dimensional analogue of Fock space of Hilbert $C^*$-bimodules.
The Hilbert $C^*$-module
$\HK$ is originally a Hilbert $C^*$-right module 
$(\HK, \varphi_\A)$

over $\A$  
with $\A$-valued inner product $\langle \cdot \mid \cdot \rangle_\A$.
It has two other 
 multi structure of Hilbert $C^*$-bimodules.
The Hilbert $C^*$-bimodule
$(\phi_\rho,\HK, \varphi_\rho)$ over $\BR$  
and
the Hilbert $C^*$-bimodule
$(\phi_\eta,\HK, \varphi_\eta)$ over $\BE$.
This situation is written as in the figure:
\begin{equation*}
\begin{CD} 
         @. \BR @.      \\
@.   @V{\phi_\rho}VV @. \\
\BE @>{\phi_\eta}>> \HK @<{\varphi_\eta}<< \BE \\ 
@.   @A{\varphi_\rho}AA @. \\
         @. \BR @.      
\end{CD}
\end{equation*}
There exist faithful
completely positive maps 
$\lambda_\rho: \BR \longrightarrow \A$
and
$\lambda_\eta: \BE \longrightarrow \A$
satisfying
\eqref{eqn:innerprod}
so that the three norms induced by their respect inner proucts
$\langle \cdot \mid \cdot \rangle_\A, 
 \langle \cdot \mid \cdot \rangle_\rho,
 \langle \cdot \mid \cdot \rangle_\eta
 $
 are equivalent to each other.
The Hilbert $C^*$-right module 
$(\HK, \varphi_\A)$
over $\A$  
with multi structure of  Hilbert $C^*$-bimodules
$(\phi_\rho,\HK, \varphi_\rho)$ over $\BR$  
and
$(\phi_\eta,\HK, \varphi_\eta)$ over $\BE$
 is called a {\it Hilbert} $C^*$-{\it quad module over} 
$(\A; \BR, \BE)$.
We will define two kinds of relative tensor products
$$
\HK \otimes_{\eta} \HK, 
\qquad
\HK \otimes_{\rho} \HK
$$
as  Hilbert $C^*$-quad modules over 
$(\A; \BR, \BE)$.
The latter one should be written  
vertically as
$$
\begin{matrix}
\HK \\
\otimes_{\rho}\\
\HK
\end{matrix}
$$
rather than horizontally
$\HK 
\otimes_{\rho}
\HK.
$ 
The first relative tensor product
is defined as
\begin{equation*}
\HK\otimes_{\eta}\HK:= \HK\otimes_{\BE}\HK
\end{equation*}
the relative tensor product as Hilbert $C^*$-modules over $\BE$,
the left $\HK$ is a right $\BE$-module through $\varphi_\eta$
and   
the right $\HK$ is a left $\BE$-module through $\phi_\eta$.
It has a right $\BR$-valued inner product 
and
a right $\BE$-valued inner product defined by
\begin{align*}
\langle \xi\otimes_\eta\zeta \mid \xi'\otimes_\eta\zeta'\rangle_\rho
& := \langle \zeta \mid \phi_\eta(\langle \xi \mid \xi'\rangle_\eta)\zeta'\rangle_\rho, \\
\langle \xi\otimes_\eta\zeta \mid \xi'\otimes_\eta\zeta'\rangle_\eta
& := \langle \zeta \mid \phi_\eta(\langle \xi \mid \xi'\rangle_\eta)\zeta'\rangle_\eta
\end{align*}
 respectively.
It has two right actions, 
$\id \otimes\varphi_\rho$ from $\BR$
and  
 $\id \otimes\varphi_\eta$ from $\BE$.
It also has two left actions,
$\phi_\rho\otimes\id$ from $\BR$
and $\phi_\eta\otimes\id$ from $\BE$.
By these operations 
$\HK\otimes_{\eta}\HK$
is a Hilbert $C^*$-bimodule over $\BR$ and also
is a Hilbert $C^*$-bimodule over $\BE$.
It also has a  right $\A$-valued inner product defined by
\begin{equation*}
\langle \xi\otimes_\eta\zeta \mid \xi'\otimes_\eta\zeta'\rangle_\A
:= \lambda_\eta(\langle \xi\otimes_\eta\zeta \mid \xi'\otimes_\eta\zeta'\rangle_\eta)
(= \lambda_\rho(\langle \xi\otimes_\eta\zeta \mid \xi'\otimes_\eta\zeta'\rangle_\rho))
\end{equation*}
and a right $\A$-action $\id\otimes \varphi_\A$ and a left $\A$-action
$\phi\otimes\id$.
By these structure
$\HK\otimes_{\eta}\HK$
is a Hilbert $C^*$-quad module over $(\A;\BR,\BE)$.
\begin{equation*}
\begin{CD} 
         @. \BR @.      \\
@.   @V{\phi_\rho\otimes\id}VV @. \\
\BE @>{\phi_\eta\otimes\id}>> \HK\otimes_{\eta}\HK @<{\id\otimes\varphi_\eta}<< \BE \\ 
@.   @A{\id\otimes\varphi_\rho}AA @. \\
         @. \BR @.      
\end{CD}
\end{equation*} 
We denote the above operations
$\phi_\rho\otimes\id, \phi_\eta\otimes\id,
\id\otimes\varphi_\rho, 
\id\otimes\varphi_\eta
$
still
by
$\phi_\rho, \phi_\eta,
\varphi_\rho, \varphi_\eta
$
respectively.
Similarly we consider the other  relative 
 tensor product
 defined by
\begin{equation*}
\HK\otimes_{\rho}\HK:= \HK\otimes_{\BR}\HK
\end{equation*}
the relative tensor product as Hilbert $C^*$-modules over $\BR$,
the left $\HK$ is a right $\BR$-module through $\varphi_\rho$
and   
the right $\HK$ is a left $\BR$-module through $\phi_\rho$.
By symmetrically to the above,
$\HK\otimes_{\rho}\HK$
is a Hilbert $C^*$-quad module over $(\A;\BR,\BE)$.
The following lemma is routine.
\begin{lem}
Let
${\cal H}_i = \HK, i=1,2,3$.
The correspondences
\begin{align*}
(\xi_1 \otimes_\eta \xi_2) \otimes_\rho \xi_3 \in ({\cal H}_1 \otimes_\eta  {\cal H}_2)\otimes_\rho {\cal H}_3
& \longrightarrow 
\xi_1 \otimes_\eta ( \xi_2 \otimes_\rho \xi_3)  \in {\cal H}_1 \otimes_\eta ( {\cal H}_2 \otimes_\rho {\cal H}_3), \\
(\xi_1 \otimes_\rho \xi_2) \otimes_\eta \xi_3 \in
({\cal H}_1 \otimes_\rho  {\cal H}_2)\otimes_\eta {\cal H}_3
& \longrightarrow
\xi_1 \otimes_\rho (\xi_2 \otimes_\eta \xi_3)
\in
{\cal H}_1 \otimes_\rho( {\cal H}_2 \otimes_\eta {\cal H}_3)
\end{align*}
yield  isomorphisms of Hilbert $C^*$-quad modules respectively.
\end{lem}
We write the isomorphism class of the former Hilbert $C^*$-quad modules as 
${\cal H}_1 \otimes_\eta  {\cal H}_2\otimes_\rho {\cal H}_3
$
and
that of the latter ones as
$
{\cal H}_1 \otimes_\rho  {\cal H}_2\otimes_\eta {\cal H}_3
$
respectively.

\medskip

We note that the direct sum
$\BE \oplus \BR$ 
has a structure of  a Hilbert $C^*$-quad module by the following operations:
For $b_1 \oplus b_2, b'_1 \oplus b'_2 \in \BE \oplus \BR$
and $ y \in \A$,
set
\begin{align*}
(b_1 \oplus b_2)  \varphi_\A(y) & := b_1\psi_\eta(y) \oplus b_2  \psi_\rho(y) \in \BE \oplus \BR, \\ 
\langle b_1 \oplus b_2  \mid b'_1 \oplus b'_2 \rangle_{\A} 
& := \lambda_\eta(b_1^* b'_1) + \lambda_\rho(b_2^* b'_2) \in \A.
\end{align*}
It is direct to see 
\begin{equation*}
\langle b_1 \oplus b_2  \mid (b'_1 \oplus b'_2) \varphi_\A(y) \rangle_{\A} 
=
\langle b_1 \oplus b_2  \mid b'_1 \oplus b'_2 \rangle_\A\cdot y  
\end{equation*}
so that 
$\BE \oplus \BR$ is a 
Hilbert $C^*$-right module over $\A$.
Its Hilbert $C^*$-bimodule structure 
over $\BR$ and over $\BE$ are defined as follows: 
For $ w \in \BR, z \in \BE$, set 

1. The right $\BR$-action $\varphi_\rho$ and the right $\BE$-action $\varphi_\eta$:
\begin{equation*}
(b_1 \oplus b_2) \varphi_\rho(w)
 =   b_2 w,  \qquad
(b_1 \oplus b_2) \varphi_\eta(z)
 =  b_1 z.
\end{equation*}

2. The left $\BR$-action $\phi_\rho$ and the left $\BE$-action $\phi_\eta$:
\begin{equation*}
\phi_\rho(w) (b_1 \oplus b_2) 
 =  w b_2, \qquad
\phi_\eta(z) (b_1 \oplus b_2) 
 =  z b_1.
\end{equation*}

3. The right $\BR$-valued inner product 
$\langle \cdot | \cdot \rangle_{\rho}$ 
and 
the right $\BE$-valued inner product 
$\langle \cdot | \cdot \rangle_{\eta}$:
\begin{equation*}
\langle b_1 \oplus b_2 \mid b'_1 \oplus b'_2 \rangle_{\rho}
= b_1^* b'_1 \in \BR, \qquad 
\langle b_1 \oplus b_2 \mid b'_1 \oplus b'_2 \rangle_{\eta}
 = b_2^* b'_2 \in \BE.
\end{equation*}
Let us  define the Fock Hilbert $C^*$-quad module
as a two-dimensional analogue of the Fock space of Hilbert $C^*$-bimodules.
Put
$\Gamma_n = \{ (\pi_1,\dots,\pi_n)) \mid \pi_i = \eta,\rho\}, n=1,2,\dots.$
We set
\begin{align*}
F_0(\kappa) & = \BE \oplus \BR, \qquad F_1(\kappa)  = \HK, \\
F_2(\kappa) & = (\HK\otimes_\eta \HK) \oplus (\HK\otimes_\rho \HK), \\
F_3(\kappa) & = (\HK\otimes_\eta \HK \otimes_\eta \HK) 
         \oplus (\HK\otimes_\eta \HK \otimes_\rho \HK) \\
         & \oplus (\HK\otimes_\rho \HK \otimes_\eta \HK) 
         \oplus (\HK\otimes_\rho \HK \otimes_\rho \HK),\\
\cdots      & \cdots \\
F_n(\kappa) & 
= \oplus_{(\pi_1,\cdots, \pi_{n-1}) \in \Gamma_{n-1}}
\HK \otimes_{\pi_1} \HK\otimes_{\pi_2}\cdots \otimes_{\pi_{n-1}} \HK \\
\cdots      & \cdots 
 \end{align*}
as Hilbert $C^*$-quad modules.
We will define the Fock Hilbert $C^*$-quad module
$F_\kappa$ by setting
\begin{equation*}
F_\kappa := \overline{\oplus_{n=0}^\infty F_n(\kappa)}
\end{equation*}
which is the completion of the algebraic direct sum
$\oplus_{n=0}^\infty F_n(\kappa)$
of the Hilbert $C^*$-quad modules
under the norm
$\| \xi \|_\A $
on
$\oplus_{n=0}^\infty F_n(\kappa)$
induced by the $\A$-valued inner product 
$\langle \cdot \mid \cdot \rangle_\A$ 
on $F_n(\kappa)$.
The Hilbert $C^*$-quad modules 
$F_n(\kappa), n=0,1,\dots $
and hence 
the algebraic direct sum
$\oplus_{n=0}^\infty F_n(\kappa)$
have also both $\BR$-valued inner products
and $\BE$-valued inner products.
As in Lemma \ref{lem:norm},
both of the two norms 
$
\| \xi \|_{\rho}
$
and
$
\| \xi \|_{\eta}
$
induced by the inner products
are equivalent to the norm
$\| \xi \|_\A$.

For $\xi\in \HK$
we define operators 
$s_\xi$ and $t_\xi$   from $F_0(\kappa)$ to $F_1(\kappa)$
by
\begin{equation*}
s_\xi(b_1 \oplus b_2) = \xi\varphi_\eta(b_1),
\qquad 
t_\xi(b_1 \oplus b_2) = \xi\varphi_\rho(b_2)
\end{equation*}
for $b_1 \oplus b_2\in \BE \oplus \BR$.
\begin{lem} \hspace{6cm}
\begin{enumerate}
\renewcommand{\labelenumi}{(\roman{enumi})}
\item
$s_\xi$ is a right $\BE$-module map from $F_0(\kappa)$ to $F_1(\kappa)$.  
\item
$t_\xi$ is a right $\BR$-module map from $F_0(\kappa)$ to $F_1(\kappa)$.  
\item
Both the maps
$s_\xi, t_\xi: F_0(\kappa) \longrightarrow F_1(\kappa)$
are right $\A$-module maps.
\end{enumerate}
\end{lem}
\begin{pf}
For $b_1 \oplus b_2\in \BE \oplus \BR$ and $z \in \BE$,
we have
\begin{equation*}
s_\xi((b_1 \oplus b_2)\varphi_\eta(z)) 
 = s_\xi(b_1 z) 
 = \xi\varphi_\eta(b_1 z)
 = (s_\xi(b_1 \oplus b_2) )\varphi_\eta( z). 
\end{equation*}
Hence 
$s_\xi$ is a right $\BE$-module map
and similarly $t_\xi$ 
is a right $\BR$-module map.
For $y \in \A$,
by Lemma 3.6, we have
\begin{equation*}
s_\xi((b_1 \oplus b_2)\varphi_\A(y)) 
 = \xi\varphi_\eta(b_1 \psi_\eta(y)) 
 = (\xi\varphi_\eta(b_1))\varphi_\A(y) 
 = (s_\xi(b_1 \oplus b_2) )\varphi_\A(y). 
\end{equation*}
Hence 
$s_\xi$ and similarly 
$t_\xi$ are right $\A$-module maps.
\end{pf}
 For $\xi \in \HK$, denote by
 $s_\xi^*, t_\xi^* : F_1(\kappa) \longrightarrow F_0(\kappa)$
 the adjoints of 
$s_\xi, t_\xi : F_0(\kappa) \longrightarrow F_1(\kappa)$
 with respect to the right $\A$-valued inner products.
\begin{lem}
For $\xi, \xi'\in \HK$, 
we have
\begin{enumerate}
\renewcommand{\labelenumi}{(\roman{enumi})}
\item
$s_\xi^* \xi' = \langle \xi \mid \xi' \rangle_\eta \oplus 0$  in 
$\BE \oplus \BR$.
\item
$t_\xi^* \xi' = 0 \oplus \langle \xi \mid \xi' \rangle_\rho$  in 
$\BE \oplus \BR$.
\end{enumerate}
\end{lem} 
\begin{pf}
For $b_1 \oplus b_2 \in \BE \oplus \BR$, we have
\begin{align*}
\langle b_1 \oplus b_2 \mid s_\xi^* \xi' \rangle_\A
& = \langle \xi\varphi_\eta(b_1) \mid  \xi' \rangle_\A \\
& = \lambda_\eta(b_1^*\langle \xi \mid  \xi' \rangle_\eta) \\
& = \lambda_\eta(\langle b_1 \oplus b_2 \mid
\langle \xi \mid  \xi' \rangle_\eta \oplus 0 \rangle_\eta) \\
& = \langle  b_1 \oplus b_2 \mid
\langle \xi \mid  \xi' \rangle_\eta \oplus 0 \rangle_\A 
\end{align*}
so that
$s_\xi^* \xi' = \langle \xi \mid \xi' \rangle_\eta \oplus 0$.

(ii) is similar to (i).
\end{pf} 
 For $\xi \in \HK$ and
$\xi_1 \otimes_{\pi_1} \cdots \otimes_{\pi_{n-1}}\xi_n \in F_n(\kappa)$
with
$(\pi_1,\dots,\pi_{n-1})\in \Gamma_{n-1}$, 
set
\begin{align}
s_\xi 
(\xi_1 \otimes_{\pi_1}
\cdots  \otimes_{\pi_{n-1}}\xi_n)
& =
\xi\otimes_\eta \xi_1 \otimes_{\pi_1}
\cdots  \otimes_{\pi_{n-1}}\xi_n, \\
t_\xi 
(\xi_1 \otimes_{\pi_1}
\cdots  \otimes_{\pi_{n-1}}\xi_n)
& =
\xi\otimes_\rho \xi_1 \otimes_{\pi_1}
\cdots  \otimes_{\pi_{n-1}}\xi_n. 
\end{align}
The following lemma is direct.
\begin{lem}
For $\xi\in \HK$ and $n=1,2,\dots $, we have 
\begin{enumerate}
\renewcommand{\labelenumi}{(\roman{enumi})}
\item
$s_\xi$ is a right $\BE$-module map from $F_n(\kappa)$ to $F_{n+1}(\kappa)$.  
\item
$t_\xi$ is a right $\BR$-module map from $F_n(\kappa)$ to $F_{n+1}(\kappa)$.  
\item
Both the maps
$s_\xi, t_\xi: F_n(\kappa) \longrightarrow F_{n+1}(\kappa)$
are right $\A$-module maps.
\end{enumerate}
\end{lem}
Denote by
$s_\xi^*, t_\xi^* : F_{n+1}(\kappa) \longrightarrow F_n(\kappa)$
 the adjoints of 
$s_\xi, t_\xi : F_n(\kappa) \longrightarrow F_{n+1}(\kappa)$
 with respect to the right $\A$-valued inner products.
\begin{lem}
For $\xi\in \HK$
and
$\xi_1 \otimes_{\pi_1}\xi_2 \otimes_{\pi_2}
\cdots  \otimes_{\pi_n}\xi_{n+1} \in F_{n+1}(\kappa)$, we have
\begin{enumerate}
\renewcommand{\labelenumi}{(\roman{enumi})}
\item
$s_\xi^* 
(\xi_1 \otimes_{\pi_1}\xi_2 \otimes_{\pi_2}
\cdots  \otimes_{\pi_n}\xi_{n+1})
 =
\begin{cases}
\phi_\eta(\langle \xi\mid  \xi_1\rangle_\eta)\xi_2 \otimes_{\pi_2}
\cdots  \otimes_{\pi_n}\xi_{n+1} 
   & \text{ if } \pi_1 = \eta, \\
0  & \text{ if } \pi_1 = \rho,
\end{cases}
$
\item
$t_\xi^* 
(\xi_1 \otimes_{\pi_1}\xi_2 \otimes_{\pi_2}
\cdots  \otimes_{\pi_n}\xi_{n+1})
 =
\begin{cases}
\phi_\rho(\langle \xi\mid  \xi_1\rangle_\rho)\xi_2 \otimes_{\pi_2}
\cdots  \otimes_{\pi_n}\xi_{n+1} 
   & \text{ if } \pi_1 = \rho, \\
0  & \text{ if } \pi_1 = \eta,
\end{cases}
$
\end{enumerate}
\end{lem}
\begin{pf}
(i)
Let $\gamma = \eta$ or $\rho$. 
For
$
\zeta_1 \otimes_{\theta_1}
\zeta_2 \otimes_{\theta_2}
\cdots
\otimes_{\theta_{n-1}}\zeta_n
\in F_n(\kappa)$,
we have 
\begin{align*}
& \langle
\zeta_1 \otimes_{\theta_1}
\zeta_2 \otimes_{\theta_2}
\cdots
\otimes_{\theta_{n-1}}\zeta_n
\mid
s_\xi^* 
(\xi_1 \otimes_{\pi_1}\xi_2 \otimes_{\pi_2}
\cdots  \otimes_{\pi_n}\xi_{n+1}) \rangle_\gamma \\
= & \langle
s_\xi(\zeta_1 \otimes_{\theta_1}
\zeta_2 \otimes_{\theta_2}
\cdots
\otimes_{\theta_{n-1}}\zeta_n)
\mid
\xi_1 \otimes_{\pi_1}\xi_2 \otimes_{\pi_2}
\cdots  \otimes_{\pi_n}\xi_{n+1} \rangle_\gamma \\
= & \langle
\xi \otimes_\eta \zeta_1 \otimes_{\theta_1}
\zeta_2 \otimes_{\theta_2}
\cdots
\otimes_{\theta_{n-1}}\zeta_n)
\mid
\xi_1 \otimes_{\pi_1}\xi_2 \otimes_{\pi_2}
\cdots  \otimes_{\pi_n}\xi_{n+1} \rangle_\gamma \\
= & 
{\begin{cases}
\langle
\zeta_1 \otimes_{\theta_1}
\zeta_2 \otimes_{\theta_2}
\cdots
\otimes_{\theta_{n-1}}\zeta_n)
\mid
\phi_\eta(\langle \xi \mid \xi_1\rangle_\eta)\xi_2 \otimes_{\pi_2}
\cdots  \otimes_{\pi_n}\xi_{n+1} \rangle_\gamma 
  & \text{ if } \pi_1 = \eta, \\
0 & \text{ if } \pi_1 = \rho.
\end{cases}} 
\end{align*}
Hence the desired formulae hold
with respect to the inner products
$
\langle  \cdot \mid \cdot \rangle_\eta,
 \langle \cdot \mid \cdot \rangle_\rho
$
and hence to 
the $\A$-valued inner product
$\langle \cdot \mid \cdot \rangle_\A$
because of the equality
$$
\langle \cdot \mid \cdot \rangle_\A
= \lambda_\eta(\langle \cdot \mid \cdot \rangle_\eta)
= \lambda_\rho(\langle \cdot \mid \cdot \rangle_\rho).
$$
\end{pf}
We have shown in the proof of the above lemma that 
the adjoints of 
$s_\xi, t_\xi : F_n(\kappa) \longrightarrow F_{n+1}(\kappa)$
 with respect to the other two inner products 
$
\langle \cdot \mid \cdot \rangle_\eta,
 \langle \cdot \mid \cdot \rangle_\rho
$ 
have the same fom as above.
We denote by
$\bar{\varphi}_\rho, \bar{\varphi}_\eta, \bar{\phi}_\rho, \bar{\phi}_\eta $
the right $\BR$-action, the right $\BE$-action, the left $\BR$-action, the left $\BE$-action
on $F_n(\kappa)$ and hence on $F_\kappa$ respectively.
The left actions $\bar\phi_\rho$ of $\BR$ 
and
 $\bar{\phi}_\eta$ of $\BE$ 
satisfy the following equalities
\begin{align*}
\bar{\phi}_\rho(w)(b_1 \oplus b_2) = w b_2, 
& \qquad  
\bar{\phi}_\eta(z)(b_1 \oplus b_2) = z b_1, \\
\bar{\phi}_\rho(w)(\xi_1 \otimes_{\pi_1}\xi_2 \otimes_{\pi_2}
\cdots  \otimes_{\pi_{n-1}}\xi_n)
& = 
(\phi_\rho(w)\xi_1) \otimes_{\pi_1}\xi_2 \otimes_{\pi_2}
\cdots  \otimes_{\pi_{n-1}}\xi_n),\\
\bar{\phi}_\eta(z)(\xi_1 \otimes_{\pi_1}\xi_2 \otimes_{\pi_2}
\cdots  \otimes_{\pi_{n-1}}\xi_n)
& = 
(\phi_\eta(z)\xi_1) \otimes_{\pi_1}\xi_2 \otimes_{\pi_2}
\cdots  \otimes_{\pi_{n-1}}\xi_n)
\end{align*}
for $w \in \BR, z \in \BE$, $b_1 \oplus b_2 \in \BE \oplus \BR$
and $\xi_1 \otimes_{\pi_1}\xi_2 \otimes_{\pi_2}
\cdots  \otimes_{\pi_{n-1}}\xi_n \in F_n(\kappa)$.
The following lemma is direct.
\begin{lem} 
For $w \in \BR, z \in \BE$, 
we have for $n=0,1,\dots$
\begin{enumerate}
\renewcommand{\labelenumi}{(\roman{enumi})}
\item
$\bar{\phi}_\rho(w)$ 
is a right $\BR$-module map from $F_n(\kappa)$ to $F_{n}(\kappa)$.  
\item
$\bar{\phi}_\eta(z) $
is a right $\BE$-module map from $F_n(\kappa)$ to $F_{n}(\kappa)$.  
\item
Both the maps
$ \bar{\phi}_\rho(w), \bar{\phi}_\eta(z)$
are right $\A$-module maps on $F_\kappa$.
\end{enumerate}
\end{lem}
\begin{lem} For $\xi \in \HK,  w \in \BR, z \in \BE$,
we have
\begin{enumerate}
\renewcommand{\labelenumi}{(\roman{enumi})}
\item
$t_{\xi\varphi_\rho(w)} = t_\xi\bar{\phi}_\rho(w)$
and hence
$t_\xi = \sum_{a \in \Sigma^\eta} 
t_{v_a}\bar{\phi}_\rho(\langle v_a \mid \xi \rangle_\rho)$.
\item
$s_{\xi\varphi_\eta(z)} = s_\xi 
\bar{\phi}_\eta(z)$
and hence
$s_\xi = \sum_{\alpha \in \Sigma^\rho} 
s_{u_\alpha}\bar{\phi}_\eta(\langle u_\alpha \mid \xi \rangle_\eta)$.
\end{enumerate}
\end{lem}
\begin{pf}
(i) 
We have
\begin{align*}
t_{\xi\varphi_\rho(w)}
(\xi_1 \otimes_{\pi_1}\xi_2 \otimes_{\pi_2}
\cdots  \otimes_{\pi_{n-1}}\xi_n)
& = \xi\varphi_\rho(w) \otimes_\rho
\xi_1 \otimes_{\pi_1}\xi_2 \otimes_{\pi_2}
\cdots  \otimes_{\pi_{n-1}}\xi_n \\
& = \xi  \otimes_\rho
(\phi_\rho(w) \xi_1) \otimes_{\pi_1}\xi_2 \otimes_{\pi_2}
\cdots  \otimes_{\pi_{n-1}}\xi_n \\
& = (t_\xi \bar{\phi}_\rho(w)) \xi_1 \otimes_{\pi_1}\xi_2 \otimes_{\pi_2}
\cdots  \otimes_{\pi_{n-1}}\xi_n).
\end{align*}
(ii) is similar to (i).
\end{pf}
By Lemma \ref{lem:4.4}, we have 
$
\varphi_\eta(\psi_\eta(y)) = \varphi_\rho(\psi_\rho(y)) =\varphi_\A(y)
$
for
$y \in \A$,
the above lemma implies the equalities:
\begin{equation}
t_{\xi\varphi_\A(y)} = t_\xi \bar{\phi}_\rho(\psi_\rho(y))
\quad
\text{and}
\quad
s_{\xi\varphi_\A(y)} = s_\xi \bar{\phi}_\eta(\psi_\eta(y))
\quad
\text{ for }
y\in \A.
\end{equation}
The following lemma 
is immediate.
\begin{lem} For $\xi \in \HK, w \in \BR, z \in \BE$,
we have
\begin{enumerate}
\renewcommand{\labelenumi}{(\roman{enumi})}
\item
$\bar{\phi}_\rho(w)s_\xi = s_{\phi_\rho(w)\xi}$ and
$\bar{\phi}_\rho(w)t_\xi = t_{\phi_\rho(w)\xi}$.
\item
$\bar{\phi}_\eta(z) s_\xi  = s_{\phi_\eta(z)\xi}$ and  
$\bar{\phi}_\eta(z) t_\xi  = t_{\phi_\eta(z)\xi}$.
\end{enumerate}
\end{lem}

We set
\begin{equation}
s_\alpha  = s_{u_\alpha} \quad \text{ for } \alpha \in \Sigma^\rho
\quad \text{ and }
\quad 
t_a  = t_{v_a} \quad \text{ for } a \in \Sigma^\eta. 
\end{equation} 
By Lemma 4.7,
 we have for $\xi \in \HK$
\begin{equation}
s_\xi  = \sum_{\alpha \in \Sigma^\rho} 
s_{\alpha} \bar{\phi}_\eta(\langle u_\alpha \mid \xi \rangle_\eta), \qquad
t_\xi  = \sum_{a \in \Sigma^\eta} 
t_{a} \bar{\phi}_\rho(\langle v_a \mid \xi \rangle_\rho). \label{eqn:expan}
\end{equation}
Let $P_n$ be  the projection on $F_\kappa$ 
onto $F_n(\kappa)$ for $ n=0,1,\dots$.
We also define 
two  projections on $F_\kappa$ by
\begin{align*}
P_\rho
& = \text{The projection onto } 
\sum_{n=1}^\infty
\sum_{(\pi_1,\cdots, \pi_{n}) \in \Gamma_{n}} 
\HK \otimes_{\eta}\HK\otimes_{\pi_1}\HK\otimes_{\pi_2}\cdots
\otimes_{\pi_{n}}\HK, \\
P_\eta
& = \text{The projection onto } 
\sum_{n=1}^\infty
\sum_{\pi_1,\cdots, \pi_{n}\in \Gamma_{n}}
\HK \otimes_{\rho}\HK\otimes_{\pi_1}\HK\otimes_{\pi_2}\cdots
\otimes_{\pi_{n}}\HK.
\end{align*}
\begin{lem} Keep the above notations.
$$
\sum_{\alpha \in \Sigma^\rho} s_\alpha s_\alpha^* = P_1 + P_\rho
\quad 
\text{ and }
\quad
\sum_{a \in \Sigma^\eta} t_a t_a^* = P_1 + P_\eta.
$$
Hence
\begin{equation}
\sum_{\alpha \in \Sigma^\rho} s_\alpha s_\alpha^* 
+
\sum_{a \in \Sigma^\eta} t_a t_a^*  + P_0 = 1_{F_\kappa}  + P_1.
\end{equation}
\end{lem}
\begin{pf}
For
$
\xi_1\otimes_{\pi_1}\xi_2\otimes_{\pi_2}
\cdots
\otimes_{\pi_{n-1}}\xi_n
\in F_n(\kappa)
$
with $2 \le n \in {\Bbb N}$,
we have
\begin{align*}
 &s_\alpha s_\alpha^*
  (\xi_1\otimes_{\pi_1}\xi_2\otimes_{\pi_2}
  \cdots
  \otimes_{\pi_{n-1}}\xi_n) \\
=&
{\begin{cases}
u_\alpha\otimes_\eta
\phi_\eta(\langle u_\alpha \mid \xi_1\rangle_\eta)\xi_2\otimes_{\pi_2}
\cdots
\otimes_{\pi_{n-1}}\xi_n)
   & \text{ if } \pi_1 = \eta,\\
0  & \text{ if } \pi_1 = \rho.
\end{cases} }
\end{align*}
As
$
u_\alpha \otimes_\eta
\phi_\eta(\langle u_\alpha \mid \xi_1\rangle_\eta)\xi_2
 = u_\alpha 
\varphi_\eta(\langle u_\alpha \mid \xi_1\rangle_\eta)
\otimes_\eta
\xi_2
$
and
$\sum_{\alpha \in \Sigma^\rho} u_\alpha 
\varphi_\eta(\langle u_\alpha \mid \xi_1\rangle_\eta)
= \xi_1$,
we have
\begin{equation*}
\sum_{\alpha \in \Sigma^\rho} s_\alpha s_\alpha^*
(\xi_1\otimes_{\pi_1}\xi_2\otimes_{\pi_2}
\cdots
\otimes_{\pi_{n-1}}\xi_n)
 =
\begin{cases}
\xi_1\otimes_{\eta}\xi_2\otimes_{\pi_2}
\cdots
\otimes_{\pi_{n-1}}\xi_n
   & \text{ if } \pi_1 = \eta,\\
0  & \text{ if } \pi_1 = \rho.
\end{cases} 
\end{equation*}
Hence we have
\begin{equation*}
\sum_{\alpha \in \Sigma^\rho} s_\alpha s_\alpha^* |_{\oplus_{n=2}^{\infty}F_n(\kappa)}
= P_\rho|_{\oplus_{n=2}^{\infty}F_n(\kappa)}.
\end{equation*}
For $\xi \in F_1(\kappa) = \HK$,
we have
$
s_\alpha s_\alpha^*
\xi
 = s_\alpha(\langle u_\alpha \mid \xi \rangle_\eta \oplus 0)
= u_\alpha \varphi_\eta(\langle u_\alpha \mid \xi \rangle_\eta) 
$
so that
\begin{equation*}
\sum_{\alpha \in \Sigma^\rho} s_\alpha s_\alpha^*
\xi
 = \sum_{\alpha \in \Sigma^\rho}
u_\alpha \varphi_\eta(\langle u_\alpha \mid \xi \rangle_\eta) = \xi.
\end{equation*}
Hence we have
\begin{equation*}
\sum_{\alpha \in \Sigma^\rho} s_\alpha s_\alpha^* |_{F_1(\kappa)}
= 1_{F_1(\kappa)}.
\end{equation*}
As
$ s_\alpha s_\alpha^*(b_1 \oplus b_2) = 0$ for 
$b_1 \oplus b_2 \in \BE \oplus \BR$,
we have
\begin{equation*}
\sum_{\alpha \in \Sigma^\rho} s_\alpha s_\alpha^* |_{F_0(\kappa)}
= 0.
\end{equation*}
Therefore
we conclude that
\begin{equation*}
\sum_{\alpha \in \Sigma^\rho} s_\alpha s_\alpha^*
= P_\rho + P_1
\quad
\text{ and similarly }
\quad
\sum_{a \in \Sigma^\eta} t_a t_a^* = P_\eta + P_1.
\end{equation*}
As 
$P_\eta + P_\rho + P_0 + P_1 = 1_{F_{\kappa}},
$
one has
\begin{equation*}
\sum_{\alpha \in \Sigma^\rho} s_\alpha s_\alpha^* 
+
\sum_{a \in \Sigma^\eta} t_a t_a^*  + P_0 = 1_{F_\kappa}  + P_1.
\end{equation*}
\end{pf}
\begin{lem} 
$s_\zeta^* s_\xi =\bar{\phi}_\eta(\langle \zeta \mid \xi \rangle_\eta)$
and
$t_\zeta^* t_\xi =\bar{\phi}_\rho(\langle \zeta \mid \xi \rangle_\rho)$
for $\zeta, \xi \in \HK$.
\end{lem}
\begin{pf}
The equalities for
$\xi_1\otimes_{\pi_1}\xi_2\otimes_{\pi_2}
\cdots
\otimes_{\pi_{n-1}}\xi_n
\in F_n(\kappa)$
\begin{align*}
s_\zeta^* s_\xi
(\xi_1\otimes_{\pi_1}\xi_2\otimes_{\pi_2}
\cdots
\otimes_{\pi_{n-1}}\xi_n)
 & =
s_\zeta^*(\xi \otimes_\eta \xi_1\otimes_{\pi_1}\xi_2\otimes_{\pi_2}
\cdots
\otimes_{\pi_{n-1}}\xi_n)\\
&=
\phi_\eta(\langle \zeta \mid \xi \rangle_\eta)
\xi_1\otimes_{\pi_1}
\xi_2\otimes_{\pi_2}
\cdots
\otimes_{\pi_{n-1}}\xi_n
\end{align*}
hold so that
$s_\zeta^* s_\xi =\bar{\phi}_\eta(\langle \zeta \mid \xi \rangle_\eta)$
on
$\oplus_{n=1}^\infty F_n(\kappa)$.
As
\begin{align*}
 s_\zeta^* s_\xi(b_1 \oplus b_2) 
& = s_\zeta^*(\xi\varphi_\eta(b_1)) \\
& = \langle \zeta \mid \xi\varphi_\eta(b_1) \rangle_\eta \oplus 0\\
& = \langle \zeta \mid \xi \rangle_\eta b_1\oplus 0 \\
& =\phi_\eta(\langle \zeta \mid \xi \rangle_\eta) (b_1 \oplus b_2), 
\end{align*}
we have
$s_\zeta^* s_\xi =\phi_\eta(\langle \zeta \mid \xi \rangle_\eta)$
on
$F_0(\kappa)$.
Hence
$s_\zeta^* s_\xi =\bar{\phi}_\eta(\langle \zeta \mid \xi \rangle_\eta)$
on
$F_\kappa$
and similarly
$t_\zeta^* t_\xi =\bar{\phi}_\rho(\langle \zeta \mid \xi \rangle_\rho)$.
\end{pf}
As $\bar{\phi}_\rho(y) = \bar{\phi}_\eta(y)$ on 
$F_n(\kappa),n=1,2,3,\dots $ for $y \in \A$,
we write
$\bar{\phi}_\rho(y)( = \bar{\phi}_\eta(y))$ 
as
$\bar{\phi}(y)$ on 
$F_n(\kappa),n=1,2,3,\dots $ for $y \in \A$.
\begin{lem}
 For $\alpha \in \Sigma^\rho, a \in \Sigma^\eta$ and 
$y \in \A$, we have
\begin{equation*}
\begin{aligned}
s_\alpha^* \bar{\phi}(y)s_\alpha & = \bar{\phi}_\eta(\rho_\alpha(y)), \\
s_\alpha s_\alpha^* \bar{\phi}_\eta(y) & =  \bar{\phi}(y)s_\alpha s_\alpha^*,
\end{aligned}
\qquad
\begin{aligned}
t_a^* \bar{\phi}(y) t_a & = \bar{\phi}_\rho(\eta_a(y)),\\
t_a t_a^* \bar{\phi}_\rho(y) & = \bar{\phi}(y)t_a t_a^*.
\end{aligned}
\end{equation*}
\end{lem}
\begin{pf}
The equalities on $F_n(\kappa), n=1,2,\dots$
\begin{align*}
 s_\alpha^* \bar{\phi}(y)s_\alpha 
 (\xi_1\otimes_{\pi_1}
 \cdots
 \otimes_{\pi_{n-1}}\xi_n) 
=&
s_\alpha^* (\phi(y)u_\alpha 
\otimes_\eta \xi_1\otimes_{\pi_1}
\cdots
\otimes_{\pi_{n-1}}\xi_n)\\
=&
\phi_\eta(\langle u_\alpha \mid \phi(y)u_\alpha\rangle_\eta) 
\xi_1\otimes_{\pi_1}
\cdots
\otimes_{\pi_{n-1}}\xi_n)\\
=&
\bar{\phi}_\eta(\rho_\alpha(y)) 
(\xi_1\otimes_{\pi_1}
\cdots
\otimes_{\pi_{n-1}}\xi_n)
\end{align*}
imply 
$s_\alpha^* \bar{\phi}(y)s_\alpha  = \bar{\phi}_\eta(\rho_\alpha(y))
$
on $F_n(\kappa), n=1,2,\dots$.
We have on $F_0(\kappa)$  
\begin{align*}
s_\alpha^* \bar{\phi}(y)s_\alpha (b_1 \oplus b_2)
& = 
\langle u_\alpha \mid  \phi(y) u_\alpha \varphi_\eta(b_1) \rangle_\eta \oplus 0 \\
& = 
\langle u_\alpha \mid  \varphi_\eta(\rho_\alpha(y) b_1 \rangle_\eta \oplus 0 \\
& = 
\langle u_\alpha \mid u_\alpha \rangle_\eta \rho_\alpha(y) b_1  \oplus 0 \\
& =  \rho_\alpha(y) b_1  \oplus 0 \\
& =  \phi_\eta(\rho_\alpha(y)) (b_1  \oplus b_2)
\end{align*}
so that
$
s_\alpha^* \bar{\phi}(y)s_\alpha
= 
\phi_\eta(\rho_\alpha(y))
$
on
$F_0(\kappa)$.
Thus we have
$
s_\alpha^* \bar{\phi}(y)s_\alpha
= 
\bar{\phi}_\eta(\rho_\alpha(y))
$
on
$F_\kappa$.
We similarly have
$
t_a^* \bar{\phi}(y) t_a  = \bar{\phi}_\rho(\eta_a(y)).
$

We also have
\begin{align*}
& s_\alpha s_\alpha^* \bar{\phi}(y)
 (\xi_1\otimes_{\pi_1}\xi_2\otimes_{\pi_2}
 \cdots
 \otimes_{\pi_{n-1}}\xi_n) \\
=&
{\begin{cases}
u_\alpha \varphi_\eta(\langle u_\alpha \mid \phi(y)\xi_1 \rangle_\eta) 
\otimes_{\eta}\xi_2\otimes_{\pi_2}
\cdots
\otimes_{\pi_{n-1}}\xi_n
  & \text{ if } \pi_1 = \eta, \\
0 & \text{ if } \pi_1 \ne \eta.  
\end{cases}}
\end{align*}
Since we have
\begin{align*}
  & u_\alpha \varphi_\eta(\langle u_\alpha \mid \phi(y)\xi_1 \rangle_\eta) \\
= & u_\alpha \varphi_\eta(\langle u_\alpha\varphi_\eta(\rho_\alpha(y*)) \mid \xi_1 \rangle_\eta) =u_\alpha \varphi_\eta(\rho_\alpha(y)\langle u_\alpha \mid \xi_1 \rangle_\eta)
=\phi(y)u_\alpha \varphi_\eta(\langle u_\alpha \mid \xi_1 \rangle_\eta),
\end{align*}
the equality
$s_\alpha s_\alpha^* \bar{\phi}(y) 
=  \bar{\phi}(y)s_\alpha s_\alpha^*$
on $F_n(\kappa), n=1,2,\dots$
holds.
For 
$b_1 \oplus b_2 \in F_0(\kappa)$,
the equality
$s_\alpha s_\alpha^* \phi_\eta(y)(b_1 \oplus b_2) 
=  \bar{\phi}(y)s_\alpha s_\alpha^*(b_1 \oplus b_2) =0$
holds
so that we conclude 
$s_\alpha s_\alpha^* \bar{\phi}_\eta(y) 
=  \bar{\phi}(y)s_\alpha s_\alpha^*$
on $F_\kappa$
and similarly 
$
t_a t_a^* \bar{\phi}_\rho(y) = \bar{\phi}_\rho(y)t_a t_a^*.
$
\end{pf}
\begin{lem}
For
$\omega = (\alpha, b, a, \beta) \in \SK$
and $y \in \A$,
we have
\begin{equation*}
t_a s_\beta t_b^* s_\alpha^* \ \bar{\phi}(y) = \bar{\phi}(y) \ t_a s_\beta t_b^* s_\alpha^*.
\end{equation*}
\end{lem}
\begin{pf}
By the preceding lemma, we have
\begin{align*}
t_a s_\beta  t_b^* s_\alpha^* \bar{\phi}(y) 
& =  t_a s_\beta  t_b^*t_b t_b^* s_\alpha^* \bar{\phi}(y) s_\alpha s_\alpha^* 
   =  t_a s_\beta \bar{\phi}( \eta_b(\rho_\alpha(y)))  t_b^* s_\alpha^*,\\
 \bar{\phi}(y) t_a  s_\beta t_b^* s_\alpha^*
& =  t_a  t_a^* \bar{\phi}(y) t_a s_\beta  s_\beta^* s_\beta t_b^* s_\alpha^*
   =  t_a s_\beta \bar{\phi}( \rho_\beta(\eta_a(y)) ) t_b^* s_\alpha^*. 
\end{align*}
The desired equality holds by \eqref{eqn:kappa}.
\end{pf}
\begin{lem} 
For $\alpha \in \Sigma^\rho,a \in \Sigma^\eta$ and $ y \in \A$, we have
$$
\bar{\phi}_\rho(S_\alpha y S_\alpha^*) P_\rho 
= s_\alpha \bar{\phi}_\rho(y) s_\alpha^*, 
\qquad
\bar{\phi}_\eta(T_a y T_a^*) P_\eta 
= t_a \bar{\phi}_\eta(y) t_a^*.
$$
\end{lem}
\begin{pf}
 For 
$
\xi_1\otimes_{\pi_1}\xi_2\otimes_{\pi_2}
\cdots
\otimes_{\pi_{n-1}}\xi_n
\in F_n(\kappa), n=2,3,\dots,
$ 
we have
\begin{align*}
& \bar{\phi}_\rho(S_\alpha y S_\alpha^*) P_\rho
(\xi_1\otimes_{\pi_1}\xi_2\otimes_{\pi_2}
\cdots
\otimes_{\pi_{n-1}}\xi_n) \\
= &
\begin{cases}
(\phi_\rho(S_\alpha y S_\alpha^*) \xi_1)
\otimes_\eta \xi_2\otimes_{\pi_2}
\cdots
\otimes_{\pi_{n-1}}\xi_n) 
  & \text{ if } \pi_1 = \eta, \\
0 & \text{ if } \pi_1 = \rho.
\end{cases}
\end{align*}
For $\xi_1 = \sum_{\omega \in \SK} 
e_\omega \otimes E_\omega x_\omega$
with $x_\omega \in \A$, 
we have
\begin{equation*}
\phi_\rho(S_\alpha y S_\alpha^*) \xi_1
=
\sum_{\omega \in \SK,t(\omega) = \alpha} 
e_\omega \otimes E_\omega \eta_{r(\omega)} (x_\omega).
\end{equation*}
On the other hand, we have
\begin{align*}
& s_\alpha \bar{\phi}_\rho(y) s_\alpha^*
(\xi_1\otimes_{\pi_1}\xi_2\otimes_{\pi_2}
\cdots
\otimes_{\pi_{n-1}}\xi_n) \\
= &
\begin{cases}
s_\alpha (\phi(y)\phi_\eta(\langle u_\alpha \mid \xi_1\rangle_\eta)
\xi_2 ) \otimes_{\pi_2}
\cdots
\otimes_{\pi_{n-1}}\xi_n) 
  & \text{ if } \pi_1 = \eta, \\
0 & \text{ if } \pi_1 = \rho.
\end{cases}
\end{align*}
As
 $\langle u_\alpha \mid \xi_1\rangle_\eta =
\sum_{\omega \in \SK, \alpha = t(\omega)} T_{r(\omega)} E_\omega x_\omega T_{r(\omega)}^*,
$ 
we have
\begin{align*}
 s_\alpha (\phi(y)(\phi_\eta(\langle u_\alpha \mid \xi_1\rangle_\eta) \xi_2)
& =s_\alpha \phi_\eta(y \langle u_\alpha \mid \xi_1\rangle_\eta) \xi_2  \\
& =\sum_{\omega \in \SK,  t(\omega)=\alpha} 
u_\alpha \otimes_\eta
\phi_\eta(
T_{r(\omega)}\eta_{r(\omega)}(y)x_\omega T_{r(\omega)}^*)\xi_2 \\
& =\sum_{\omega \in \SK, t(\omega)=\alpha}  u_\alpha \varphi_\eta
(
T_{r(\omega)}\eta_{r(\omega)}(y) x_\omega T_{r(\omega)}^*) \otimes_\eta \xi_2  \\
&
=\sum_{\omega \in \SK, t(\omega)=\alpha} 
(e_\omega \otimes 
E_\omega \eta_{r(\omega)}(y) x_\omega)
\otimes_\eta \xi_2.
\end{align*}
Hence we have
$\bar{\phi}_\rho(S_\alpha y S_\alpha^*) P_\rho 
= s_\alpha \bar{\phi}_\rho(y) s_\alpha^*$ 
on $F_n(\kappa) , n=2,3,\dots.$

For $\xi \in \HK$,
we have
$
\bar{\phi}_\rho(y) s_\alpha^* \xi
= \bar{\phi}_\rho(y) (\langle u_\alpha \mid \xi \rangle_\eta \oplus 0)
=0
$
so that
\begin{equation*}
s_\alpha \bar{\phi}_\rho(y) s_\alpha^* \xi
=\bar{\phi}_\rho(S_\alpha y S_\alpha^*) P_\rho \xi
=0.
\end{equation*}
Therefore we have
$\bar{\phi}_\rho(S_\alpha y S_\alpha^*) P_\rho 
= s_\alpha \bar{\phi}_\rho(y) s_\alpha^*$ 
on $F_\kappa.$

The other equality
$
\bar{\phi}_\eta(T_a y T_a^*) P_\eta 
= t_a \bar{\phi}_\eta(y) t_a^*
$ is similarly shown.
\end{pf}
\begin{lem} 
For $w \in \BR, z \in \BE$ and $\alpha \in \Sigma^\rho, a\in \Sigma^\eta$,
 we have
\begin{equation*}
\begin{aligned}
s_\alpha  s_\alpha^* \bar{\phi}_\rho(w) & = \bar{\phi}_\rho(w) s_\alpha  s_\alpha^*,\\
s_\alpha  s_\alpha^* \bar{\phi}_\eta(z) & = \bar{\phi}_\eta(z) s_\alpha  s_\alpha^*,
\end{aligned}
\qquad
\begin{aligned}
 t_a  t_a^* \bar{\phi}_\rho(w) & = \bar{\phi}_\rho(w) t_a  t_a^*, \\
t_a  t_a^* \bar{\phi}_\eta(z) & =\bar{\phi}_\eta(z) t_a  t_a^*
\end{aligned}
\end{equation*}
and hence
\begin{equation*}
\begin{aligned}
P_\rho \bar{\phi}_\rho(w) & = \bar{\phi}_\rho(w) P_\rho,\\
P_\rho \bar{\phi}_\eta(z) & = \bar{\phi}_\eta(z) P_\rho,
\end{aligned} 
\qquad
\begin{aligned}
P_\eta \bar{\phi}_\rho(w) & = \bar{\phi}_\rho(w) P_\eta,\\
P_\eta \bar{\phi}_\eta(z) & = \bar{\phi}_\eta(z) P_\eta.
\end{aligned}
\end{equation*}
\end{lem}
\begin{pf}
 For 
$
\xi_1\otimes_{\pi_1}\xi_2\otimes_{\pi_2}
\cdots
\otimes_{\pi_{n-1}}\xi_n
\in F_n(\kappa),
$ 
we have
\begin{align*}
& \bar{\phi}_\rho(w) s_\alpha  s_\alpha^*
(\xi_1\otimes_{\pi_1}\xi_2\otimes_{\pi_2}
\cdots
\otimes_{\pi_{n-1}}\xi_n) \\
= &
\begin{cases}
\phi_\rho(w) u_\alpha \varphi_\eta( \langle u_\alpha \mid \xi_1\rangle_\eta)
\otimes_\eta \xi_2\otimes_{\pi_2}
\cdots
\otimes_{\pi_{n-1}}\xi_n) 
  & \text{ if } \pi_1 = \eta, \\
0 & \text{ if } \pi_1 = \rho
\end{cases}
\end{align*}
and
\begin{align*}
&  s_\alpha  s_\alpha^* \bar{\phi}_\rho(w)
(\xi_1\otimes_{\pi_1}\xi_2\otimes_{\pi_2}
\cdots
\otimes_{\pi_{n-1}}\xi_n) \\
= &
\begin{cases}
u_\alpha\varphi_\eta( \langle u_\alpha \mid \phi_\rho(w) 
 \xi_1 \rangle_\eta)
\otimes_\eta \xi_2\otimes_{\pi_2}
\cdots
\otimes_{\pi_{n-1}}\xi_n) 
  & \text{ if } \pi_1 = \eta, \\
0 & \text{ if } \pi_1 = \rho.
\end{cases}
\end{align*}
Since
$
\phi_\rho(w) u_\alpha
= u_\alpha\varphi_\eta( \langle u_\alpha \mid \phi_\rho(w) u_\alpha\rangle_\eta),
$
we  have
\begin{align*}
\phi_\rho(w) u_\alpha 
\varphi_\eta( \langle u_\alpha \mid \xi_1\rangle_\eta)
& =
u_\alpha\varphi_\eta( \langle u_\alpha \mid \phi_\rho(w) u_\alpha
\varphi_\eta( \langle u_\alpha \mid \xi_1\rangle_\eta) \rangle_\eta) \\
& =
u_\alpha\varphi_\eta( \langle u_\alpha \mid \phi_\rho(w) 
\sum_{\beta \in \Sigma^\rho}u_\beta
\varphi_\eta( \langle u_\beta \mid \xi_1\rangle_\eta) \rangle_\eta) \\
& =
u_\alpha\varphi_\eta( \langle u_\alpha \mid \phi_\rho(w) 
 \xi_1 \rangle_\eta). 
\end{align*}
Hence
$s_\alpha  s_\alpha^* \bar{\phi}_\rho(w)  = \bar{\phi}_\rho(w) s_\alpha  s_\alpha^*$
holds.
Similarly 
by 
$
\phi_\eta(z) u_\alpha
= u_\alpha\varphi_\eta( \langle u_\alpha \mid \phi_\eta(z) u_\alpha\rangle_\eta),
$
the equality
$\bar{\phi}_\eta(z) s_\alpha  s_\alpha^* = s_\alpha  s_\alpha^* \bar{\phi}_\eta(z)$
holds.
As 
$
P_\rho = \sum_{\alpha \in \Sigma^\rho} s_\alpha s_\alpha^* 
$ 
on $\oplus_{n=2}^\infty F_n(\kappa)$,
we see
that
$P_\rho$ commutes with
$\bar{\phi}_\rho(w) $ and $\bar{\phi}_\eta(z)$.

The other four equalities in the right hand side are similarly shown.
\end{pf}
Let us denote by
${\cal L}_\A(\HK)$
and
${\cal L}_\A(F_\kappa)$
the $C^*$-algebras of all bounded adjointable 
right $\A$-module maps on $\HK$ and on $F_\kappa$
with respect to the right $\A$-valued 
inner products respectively. 
For $L \in {\cal L}_\A(\HK)$,
define $\bar{L} \in {\cal L}_\A(F_\kappa)$
by
\begin{equation*}
\bar{L}(\xi_1\otimes_{\pi_1}\xi_2\otimes_{\pi_2}
\cdots
\otimes_{\pi_{n-1}}\xi_n) 
=(L\xi_1)\otimes_{\pi_1}\xi_2\otimes_{\pi_2}
\cdots
\otimes_{\pi_{n-1}}\xi_n 
\end{equation*}
for $\xi_1\otimes_{\pi_1}\xi_2\otimes_{\pi_2}
\cdots
\otimes_{\pi_{n-1}}\xi_n \in F_\kappa.
$
\begin{lem} 
For $L \in {\cal L}_\A(\HK)$ and
$\alpha \in \Sigma^\rho, a \in \Sigma^\eta$,
 we have
\begin{equation*}
s_\alpha^* \bar{L} s_\alpha 
= \bar{\phi}_\eta(\langle u_\alpha \mid L u_\alpha \rangle_\eta),
\qquad 
t_a^* \bar{L} t_a 
= \bar{\phi}_\rho(\langle v_a \mid L v_a \rangle_\rho).
\end{equation*}
\end{lem}
\begin{pf}
For 
$
\xi_1\otimes_{\pi_1}
 \cdots
 \otimes_{\pi_{n-1}}\xi_n \in F_n(\kappa), n=1,2,\dots$,
 we have
\begin{align*}
 s_\alpha^* \bar{L} s_\alpha 
  (\xi_1\otimes_{\pi_1}
  \cdots
  \otimes_{\pi_{n-1}}\xi_n) 
=&
s_\alpha^* ( L u_\alpha 
\otimes_\eta \xi_1\otimes_{\pi_1}
\cdots
\otimes_{\pi_{n-1}}\xi_n)\\
=&
\bar{\phi}_\eta(\langle u_\alpha \mid L u_\alpha\rangle_\eta) 
\xi_1\otimes_{\pi_1}
\cdots
\otimes_{\pi_{n-1}}\xi_n.
\end{align*}
For $b_1 \oplus b_2 \in \BE \oplus \BR$,
we have
\begin{equation*}
 s_\alpha^* \bar{L} s_\alpha(b_1 \oplus b_2)
= s_\alpha^* L(u_\alpha \varphi_\eta(b_1))
= \langle u_\alpha \mid Lu_\alpha\varphi_\eta(b_1) \rangle_\eta \oplus 0  
= \langle u_\alpha \mid Lu_\alpha \rangle_\eta b_1 \oplus 0  
 \end{equation*}
Since
$\langle u_\alpha \mid Lu_\alpha \rangle_\eta b_1 \oplus 0
=
\bar{\phi}_\eta(\langle u_\alpha \mid Lu_\alpha \rangle_\eta )(b_1 \oplus b_2)
$
we have
$$ 
s_\alpha^* \bar{L} s_\alpha = \bar{\phi}_\eta(\langle u_\alpha \mid Lu_\alpha \rangle_\eta )
\qquad
\text { on } F_\kappa.
$$
The other equality
for
$t_a^* \bar{L} t_a $ is similarly shown.
\end{pf}
\begin{lem} 
For $ w \in \BR, z \in \BE$ 
and 
$\alpha \in \Sigma^\rho, a \in \Sigma^\eta$,
we have 
\begin{equation*}
\begin{aligned}
s_\alpha^* \bar{\phi}_\rho(w)  s_\alpha 
& = 
\bar{\phi}_\eta( \widehat{\rho}_\alpha (w)) \\
s_\alpha^* \bar{\phi}_\eta(z)  s_\alpha 
& = 
\bar{\phi}_\eta( \widehat{\rho}^\eta_\alpha (z)),
\end{aligned}
\qquad
\begin{aligned}
t_a^* \bar{\phi}_\eta(z)  t_a 
& = \bar{\phi}_\rho(\widehat{\eta}_a (z)),\\
t_a^* \bar{\phi}_\rho(w)  t_a 
& = 
\bar{\phi}_\rho( \widehat{\eta}^\rho_a (w)).
\end{aligned}
\end{equation*}
\end{lem}
\begin{pf}
By Lemma 3.12, 
we have
$$
 \widehat{\rho}_\alpha(w) = \langle u_\alpha \mid \phi_\rho(w) u_\alpha \rangle_\eta,
 \qquad
\widehat{\rho}^\eta_\alpha (z) = \langle u_\alpha \mid \phi_\eta(z) u_\alpha \rangle_\eta.
 $$
Hence the preceding lemma implies 
\begin{align*}
s_\alpha^* \bar{\phi}_\rho(w)s_\alpha 
& = \bar{\phi}_\eta(\langle u_\alpha \mid \phi_\rho(w) u_\alpha \rangle_\eta )
=\bar{\phi}_\eta(\widehat{\rho}_\alpha(w)), \\
s_\alpha^* \bar{\phi}_\eta(z)s_\alpha
& =
\bar{\phi}_\eta(\langle u_\alpha \mid \phi_\eta(z) u_\alpha \rangle_\eta)
= 
\bar{\phi}_\eta(\widehat{\rho}^\eta_\alpha(z)).
\end{align*}
The other two equalities in the right hand side are similarly shown.
\end{pf}
\begin{cor} 
For 
$w = \sum_{\alpha \in \Sigma^\rho}S_\alpha w_\alpha S_\alpha^*$
as in \eqref{eqn:w}
and 
$z = \sum_{a \in \Sigma^\eta}T_a z_a T_a^*$
as in \eqref{eqn:z}, we have
\begin{align}
\bar{\phi}_\rho(w)  & = \sum_{\alpha \in \Sigma^\rho}
s_\alpha \bar{\phi}_\eta (w_\alpha) s_\alpha^*  + P_\eta \bar{\phi}_\rho(w) P_\eta
+ P_0 \bar{\phi}_\rho(w) P_0, \label{eqn:barphirho}\\
\bar{\phi}_\eta(z)  & = \sum_{a \in \Sigma^\eta}
t_a \bar{\phi}_\rho(z_a) t_a^* + P_\rho \bar{\phi}_\eta(z) P_\rho
+ P_0 \bar{\phi}_\eta(z) P_0. \label{eqn:barphieta}
\end{align}
\end{cor}
\begin{pf}
As $P_\rho + P_\eta + P_1 + P_0 = 1$ on  $F_\kappa$
and $\sum_{\alpha \in \Sigma^\rho} s_\alpha s_\alpha^* = P_\rho + P_1$ ,
one has 
\begin{equation*}
\bar{\phi}_\rho(w)  
 = \sum_{\alpha \in \Sigma^\rho} s_\alpha s_\alpha^* \bar{\phi}_\rho(w)  s_\alpha s_\alpha^*
+ \bar{\phi}_\rho(w) P_\eta +   \bar{\phi}_\rho(w) P_0.
\end{equation*}
Since
$s_\alpha^* \bar{\phi}_\rho(w)  s_\alpha 
= \bar{\phi}_\eta(\widehat{\rho}_\alpha(w))
= \bar{\phi}_\eta(w_\alpha),
$
we have the equality  \eqref{eqn:barphirho}.
The  equality \eqref{eqn:barphieta} is similarly shown.
\end{pf}
\begin{lem}
\hspace{6cm}
\begin{enumerate}
\renewcommand{\labelenumi}{(\roman{enumi})}
\item
$\bar{\phi}_\rho : \BR \longrightarrow {\cal L}_\A(F_\kappa)$
 is a faithful $*$-homomorphism.
\item 
$\bar{\phi}_\eta : \BE \longrightarrow {\cal L}_\A(F_\kappa)$ 
is a faithful $*$-homomorphism.
\end{enumerate}
\end{lem} 
\begin{pf}
(i)
It is enough to show that 
$\phi_\rho: \BR \longrightarrow {\cal L}_\A(\HK)$ is injective.
For 
$w = \sum_{\alpha \in \Sigma^\rho} S_\alpha w_\alpha S_\alpha^*$
as in \eqref{eqn:w},
suppose that
$\phi_\rho(w) =0$
on $\HK$.
By Lemma 3.12, 
we have
$\widehat{\rho}_\alpha(w) =0$
for all $\alpha \in \Sigma^\rho$
so that
$w_\alpha =0$ for all $\alpha \in \Sigma^\rho$,
which shows $w=0$.
(ii) is similar to (i).
\end{pf}

\section{The $C^*$-algebras associated to the Hilbert $C^*$-quad modules}
In this section, 
we will study the $C^*$-algebras generated by
the operators
$s_\xi, t_\xi$ for $\xi \in \HK$.
For $\xi,\zeta \in F_\kappa$,
denote by 
$\theta_{\xi, \zeta}$
the rank one operator on $F_\kappa$
defined by
\begin{equation*}
\theta_{\xi, \zeta}(\gamma) =
\xi \varphi_\A(\langle \zeta \mid \gamma \rangle_\A)
\qquad \text{ for } \gamma \in F_\kappa.
\end{equation*}
It is immediate to see that
the operators 
$\theta_{\xi, \zeta}$
for $\xi, \zeta \in F_\kappa$
are  $\A$-module maps through $\varphi_\A$.
Let us denote by
${\cal K}_\A(F_\kappa)$
the $C^*$-subalgebra 
of ${\cal L}_\A(F_\kappa)$
generated by the rank one operators
$\theta_{\xi, \zeta}$
for $\xi,\zeta \in F_\kappa$.
Put the projections for $\alpha \in \Sigma^\rho, a \in \Sigma^\eta$
\begin{equation*}
p_\alpha  = S_\alpha S_\alpha^* 
\in \BR \subset F_0(\kappa), \qquad
q_a  = T_a T_a^* 
\in \BE \subset F_0(\kappa).
\end{equation*}
They are regarded as vectors in $F_\kappa$.
\begin{lem}
$
\sum_{\alpha \in \Sigma^\rho} \theta_{p_\alpha,p_\alpha}
+
\sum_{a \in \Sigma^\eta} \theta_{q_a,q_a} 
= P_0 : \text{the projection on }F_\kappa \text{ onto } F_0(\kappa). 
$
\end{lem}
\begin{pf}
For $\alpha \in \Sigma^\rho$ and $b_2 \in \BR$,
we have
\begin{align*}
\theta_{p_\alpha, p_\alpha}(b_2)
& = p_\alpha \varphi_\A(\langle p_\alpha \mid b_2 \rangle_\A) \\
& = p_\alpha \psi_\eta(\lambda_\eta(p_\alpha^*  b_2)) \\
& = p_\alpha \sum_{\alpha',\beta' \in \Sigma^\rho}
        S_{\alpha'} S_{\beta'}^* p_\alpha b_2 S_{\beta'}S_{\alpha'}^* 
 = p_\alpha b_2
\end{align*}
so that
$
\sum_{\alpha \in \Sigma^\rho}\theta_{p_\alpha,p_\alpha}(b_2) = b_2
$
for
$b_2 \in \BR$.
Similarly we have
$\theta_{q_a,q_a}(b_1) = q_a b_1$
for $q_a \in \BE$
so that 
$
\sum_{a \in \Sigma^\eta} \theta_{q_a,q_a}(b_1) = b_1.
$
As 
$\theta_{q_a,q_a}(b_2)=0$ for $b_2 \in \BR$,
$\theta_{p_\alpha,p_\alpha}(b_1)=0$ for $b_1 \in \BE$
and
$$
\theta_{p_\alpha, p_\alpha}(\xi)
=\theta_{q_a,q_a}(\xi)
=0
\quad 
\text{ for }
\quad
\xi \in F_n(\kappa), n=1,2,\dots,
$$ 
the operator
$
\sum_{\alpha \in \Sigma^\rho} \theta_{p_\alpha,p_\alpha}
+
\sum_{a \in \Sigma^\eta} \theta_{q_a,q_a} 
$
is the projection on 
$F_\kappa$
onto
$
F_0(\kappa). 
$
\end{pf}
Put
$\epsilon_\omega := e_\omega \otimes E_\omega \in \HK$ 
for
$\omega \in \SK$.
Then we see 
\begin{equation*}
\langle \epsilon_\omega \mid \epsilon_{\omega'} \rangle_\A
=
\begin{cases}
E_\omega & \text{ if } \omega = \omega',\\
  0      & \text{ if } \omega \ne \omega'.
\end{cases}
\end{equation*}
\begin{lem}
$\{\epsilon_\omega\}_{\omega \in \SK}$ 
forms an orthogonal basis of $\HK$
with respect to the $\A$-valued inner product 
$\langle \cdot \mid \cdot \rangle_\A$
as a right $\A$-module through $\varphi_\A$.
\end{lem}
\begin{pf}
For 
$\xi =\sum_{\omega' \in \SK}e_{\omega'} \otimes E_{\omega'} x_{\omega'}$
with $x_{\omega'} \in \A$,
one has 
\begin{equation*}
\langle \epsilon_\omega \mid \xi \rangle_\A
=\sum_{\omega' \in \SK} \langle \epsilon_\omega \mid
\epsilon_{\omega'} \rangle_\A x_{\omega'}
= E_\omega x_{\omega}
\end{equation*}
so that
\begin{equation}
\xi = 
 \sum_{\omega \in \SK} \epsilon_\omega\varphi_\A( E_\omega x_\omega)
=
\sum_{\omega \in \SK}
\epsilon_\omega \varphi_\A(\langle \epsilon_\omega \mid \xi 
\rangle_\A). \label{eqn:theta}
\end{equation}
\end{pf}
\begin{lem}
$\sum_{\omega \in \SK} \theta_{\epsilon_\omega, \epsilon_\omega} =P_1$
the projection 
on $F_\kappa$ onto $F_1(\kappa)$.
\end{lem}
\begin{pf}
By \eqref{eqn:theta}, we have
$
\xi 
=\sum_{\omega \in \SK} \theta_{\epsilon_\omega, \epsilon_\omega}(\xi)
$
for $\xi \in \HK$.
Since
$\theta_{\epsilon_\omega, \epsilon_\omega}(\xi')=0$
for
$\xi' \in F_n(\kappa)$ with $n \ne 1$,
we have
$\sum_{\omega \in \SK} \theta_{\epsilon_\omega, \epsilon_\omega} =P_1$.
\end{pf}
By the preceding lemmas,  we have
\begin{cor}
$P_0, P_1 \in {\cal K}_\A(F_\kappa)$.
\end{cor}
The $C^*$-subalgebra of ${\cal L}_\A(F_\kappa)$
generated by 
the operators
$s_\xi, t_\xi$ for $\xi \in \HK$
is denoted by
${\cal T}_\HK$ 
and is called the  Toeplitz quad module algebra.

\noindent
{\bf Definition.}
The $C^*$-algebra 
$\OHK$ associated with the Hilbert $C^*$-quad module $\HK$
is defined as the quotient $C^*$-algebra
of ${\cal T}_\HK$
by the ideal
${\cal T}_\HK \cap {\cal K}_\A(F_\kappa)$.

We set the quotients of the operators in $\OHK$ for 
$ \alpha \in \Sigma^\rho, a \in \Sigma^\eta$:
\begin{equation*}
U_\alpha  := [s_\alpha] \in \OHK \qquad 
V_a  := [t_a] \in \OHK.
\end{equation*}
Since
$$
   \sum_{\alpha \in \Sigma^\rho} s_\alpha s_\alpha^*  
+  \sum_{a \in \Sigma^\eta} t_a t_a^* + P_0 =1 + P_1 
$$
 and
$P_0, P_1 \in {\cal K}_\A(F_\kappa)$ 
 we have
\begin{equation}
\sum_{\alpha \in \Sigma^\rho} 
U_\alpha U_\alpha^* +
 \sum_{a \in \Sigma^\eta} V_a V_a^* =1.
\end{equation}
We also have for $w \in \BR$
$$
\bar{\phi}_\rho(w)
=
   \sum_{\alpha \in \Sigma^\rho} \bar{\phi}_\rho(w) s_\alpha s_\alpha^*  
+  \sum_{a \in \Sigma^\eta} \bar{\phi}_\rho(w) t_a t_a^* + \bar{\phi}_\rho(w) P_0  
- \bar{\phi}_\rho(w) P_1. 
$$
As 
$\bar{\phi}_\rho(w)s_\alpha = s_{\phi_\rho(w) u_\alpha}$ and
$\bar{\phi}_\rho(w)t_a = t_{\phi_\rho(w) v_a}$
by Lemma 4.8,
the operators
$\bar{\phi}_\rho(w)$ for $w \in \BR$
and similarly
$\bar{\phi}_\eta(z)$ for $z \in \BE$
belong to
$\TK$ modulo
${\cal K}_\A(F_\kappa)$.
Let
$\Phi_\rho(w), \Phi_\eta(z) \in \OHK$
denote the quotient images of
$\bar{\phi}_\rho(w), \bar{\phi}_\eta(z)$
for $ w \in \BR, z \in \BE$ in the quotient
$\OHK = \TK / \TK \cap {\cal K}_\A(F_\kappa).$ 
The following lemma is clear by 
\eqref{eqn:expan}.
\begin{lem}
The $C^*$-algebra $\OHK$ 
is generated by the
partial isometries 
$U_\alpha, V_a$ 
for $ \alpha \in \Sigma^\rho, a \in \Sigma^\eta$
 and the elements
$\Phi_\rho(w), \Phi_\eta(z)$
for $ w \in \BR, z \in \BE$.
\end{lem}

We will show that
the $C^*$-algebra
$\OHK$ has a universal property subject to the operator relations
inheritated from Lemma 4.14 and Lemma 4.16.
The following proposition is direct.
\begin{prop}
The operators
$\Phi_\rho(w), \Phi_\eta(z), $
$U_\alpha, V_a $
for $w \in \BR, z \in \BE,$
$\alpha\in \Sigma^\rho,a \in \Sigma^\eta
$
satisfy the relations: 
\begin{align*}
\sum_{\alpha \in \Sigma^\rho} U_\alpha U_\alpha^* +
& \sum_{a \in \Sigma^\eta} V_a V_a^* =1, \\ 
U_\alpha U_\alpha^* \Phi_\rho(w)  = \Phi_\rho(w) U_\alpha U_\alpha^*, 
& \qquad
V_a V_a^* \Phi_\rho(w)  = \Phi_\rho(w) V_a V_a^*, \\
U_\alpha U_\alpha^* \Phi_\eta(z)  = \Phi_\eta(z) U_\alpha U_\alpha^*, 
& \qquad
V_a V_a^* \Phi_\eta(z)  = \Phi_\eta(z) V_a V_a^*, \\
\Phi_\rho(\widehat{\rho}_\alpha(w))
=
U_\alpha^* \Phi_\rho(w) U_\alpha, 
&  \qquad
\Phi_\eta(\widehat{\eta}_a(z))
=
V_a^* \Phi_\eta(z) V_a,\\ 
\Phi_\eta(\widehat{\rho}^\eta_\alpha(z))
=
U_\alpha^* \Phi_\eta(z) U_\alpha, 
&  \qquad
\Phi_\rho(\widehat{\eta}^\rho_a(w))
=
V_a^* \Phi_\rho(w) V_a,\\ 
\Phi_\rho(y) & = \Phi_\eta(y) 
\end{align*}
for
$
w \in \BR, z \in \BE, \alpha \in \Sigma^\rho, a \in \Sigma^\eta,
y \in \A.
$
\end{prop}
The ten relations above are called the relations $(\HK)$.
We will henceforth prove that the $C^*$-algebra
$\OHK$ has the universal property subject to the relations 
$(\HK)$.
Let
$\BK$ be the $C^*$-subalgebra of $\OHK$
generated by the operators
$\Phi_\rho(w), \Phi_\eta(z) $
for $w \in \BR, z \in \BE$.
\begin{lem}
Assume that the algebra $\A$ is commutative.
Then 
$\BK$ is commutative
by the relations $(\HK)$.
\begin{pf}
As the algebra $\A$ is commutative,
 the algebras 
$\BR$ and $\BE$ are both commutative
by Lemma 3.8.
Hence it is enough to prove that
$\Phi_\rho(w)$ 
commutes with $\Phi_\eta(z)$  
for $w \in \BR, z \in \BE$.
For $\alpha \in \Sigma^\rho$,
it follows that
\begin{equation*}
\Phi_\eta(z)\Phi_\rho(w)
U_\alpha U_\alpha^*
 = U_\alpha  
\Phi_\eta(\widehat{\rho}^\eta_\alpha(z)) 
\Phi_\rho(\widehat{\rho}_\alpha(w)) U_\alpha^*. 
\end{equation*}
As
$\widehat{\rho}_\alpha(w) \in \A$,
we have
$\Phi_\rho(\widehat{\rho}^\eta_\alpha(w))
=
\Phi_\eta(\widehat{\rho}^\eta_\alpha(w))
$
so that
\begin{align*}
\Phi_\eta(\widehat{\rho}^\eta_\alpha(z)) \Phi_\rho(\widehat{\rho}_\alpha(w))
& =  
\Phi_\eta(\widehat{\rho}^\eta_\alpha(z)\widehat{\rho}_\alpha(w)) \\
\Phi_\rho(\widehat{\rho}_\alpha(w)) \Phi_\eta(\widehat{\rho}^\eta_\alpha(z))
& =
\Phi_\eta(\widehat{\rho}_\alpha(w)\widehat{\rho}^\eta_\alpha(z)).  
\end{align*}
Both the elements
$\widehat{\rho}^\eta_\alpha(z), \widehat{\rho}_\alpha(w)$
belong to the commutative algebra $\BE$,
so that
\begin{equation*}
\Phi_\eta(z)\Phi_\rho(w)
U_\alpha U_\alpha^*
 = U_\alpha  
\Phi_\rho(\widehat{\rho}_\alpha(w)) \Phi_\eta(\widehat{\rho}^\eta_\alpha(z))
 U_\alpha^* 
 =  
 \Phi_\rho(w)  \Phi_\eta(z)U_\alpha U_\alpha^*.
\end{equation*}
Similarly we have 
\begin{equation*}
\Phi_\eta(z)\Phi_\rho(w)
V_a V_a^*
=\Phi_\rho(w)  \Phi_\eta(z) V_a V_a^*
\end{equation*}
for $a \in \Sigma^\eta$.
As
$
\sum_{\alpha \in \Sigma^\rho} U_\alpha U_\alpha^* +
 \sum_{a \in \Sigma^\eta} V_a V_a^* =1,
$
one concludes that
\begin{equation*}
\Phi_\eta(z)\Phi_\rho(w)
=\Phi_\rho(w)  \Phi_\eta(z). 
\end{equation*}
\end{pf}
\end{lem}
Put
$\Sigma^{\rho\cup\eta} = \Sigma^\rho \cup \Sigma^\eta$.
We set for $\gamma \in \Sigma^{\rho\cup\eta}$ and $X \in \BK$
\begin{equation*}
\rho_\gamma^\kappa(X)
= W_\gamma^* X W_\gamma
\quad
\text{ where }
\quad
W_\gamma
 = 
{\begin{cases}
U_\alpha & \text{ if } \gamma = \alpha \in \Sigma^\rho, \\
V_a & \text{ if } \gamma = a \in \Sigma^\eta.
\end{cases}}
\end{equation*}
Since
$
U_\alpha^* \BK U_\alpha \subset \BK
$
and
$
V_a^* \BK V_a \subset \BK,
$
we have
$$
\rho_\gamma^\kappa(\BK) \subset \BK,
$$
so that 
we have a family of endomorphisms
$\rho_\gamma^\kappa, \gamma \in \Sigma^{\rho\cup\eta}
$
on
$\BK$.
In what follows,
we assume that
the algebra $\A$ is commutative,
so that the algebras $\BR, \BE$ and $\BK$ 
are all commutative. 
\begin{lem}
The triplet
$(\BK, \rho^\kappa, \Sigma^{\rho\cup\eta})$
is a $C^*$-symbolic dynamical system.
\end{lem}
\begin{pf}
Since
$\sum_{\gamma \in \Sigma^{\rho\cup\eta}}W_\gamma W_\gamma^*
= \sum_{\alpha \in \Sigma^\rho} U_\alpha U_\alpha^* 
+
\sum_{a \in \Sigma^\eta} V_a V_a^* =1
 $
 and
$W_\gamma W_\gamma^*$
commutes with
$\BK$,
the family
$\rho_\gamma^\kappa, \gamma \in \Sigma^{\rho\cup\eta}
$
yields an endomorphisms on
$\BK$.
We have
\begin{align*}
\sum_{\gamma \in \Sigma^{\rho\cup\eta}}
\rho_\gamma^\kappa(1) 
& =
\sum_{\alpha \in \Sigma^\rho} U_\alpha^* U_\alpha 
+
\sum_{a \in \Sigma^\eta} V_a^* V_a \\
& =
\sum_{\alpha \in \Sigma^\rho} \Phi_\rho(\rho_\alpha(1)) 
+
\sum_{a \in \Sigma^\eta} \Phi_\eta(\eta_a(1)) \\
& \ge \Phi_\rho(1) + \Phi_\eta(1) \ge 2.
\end{align*} 
Hence
$(\BK, \rho^\kappa, \Sigma^{\rho\cup\eta})$
is a $C^*$-symbolic dynamical system.
\end{pf}
For 
$\mu = \mu_1 \cdots \mu_n \in B_n(\Sigma^{\rho\cup\eta})$
where 
$\mu_1,  \dots, \mu_n \in \Sigma^{\rho\cup\eta}$,
denote by
\begin{align*}
|\mu|_\rho & = \text{the number of symbols of  } \Sigma^\rho 
\text{ appearing in the word } \mu_1 \cdots \mu_n, \\
|\mu|_\eta & = \text{the number of symbols of } \Sigma^\eta 
\text{ appearing in the word } \mu_1 \cdots \mu_n.
\end{align*}
Hence
$|\mu|_\rho + |\mu|_\eta =n$.
For $n \in \Zp$,
denote by
${\cal F}_n$
the $C^*$-subalgebra of $\OHK$
generated by 
the operators
$W_{\gamma_1 \cdots \gamma_n} b W_{\gamma'_1 \cdots \gamma'_n}^*
$
for
$b \in \BK$
and
$\gamma_1 \cdots \gamma_n \in B_n(\Sigma^{\rho\cup\eta})$
such that
$ |\gamma_1 \cdots \gamma_n|_\rho = |\gamma'_1 \cdots \gamma'_n|_\rho$ 
and
$ |\gamma_1 \cdots \gamma_n|_\eta = |\gamma'_1 \cdots \gamma'_n|_\eta$.
Since
$
\sum_{\gamma \in \Sigma^{\rho\cup\eta}}
W_\gamma W_\gamma^* 
=1
$
and
$W_\gamma W_\gamma^*$
commutes with
$\BK$,
the equality
\begin{equation*}
W_{\gamma_1 \cdots \gamma_n} b W_{\gamma'_1 \cdots \gamma'_n}^*
=
\sum_{\gamma_{n+1} \in \Sigma^{\rho\cup\eta}}
W_{\gamma_1 \cdots \gamma_n}W_{\gamma_{n+1}} 
\rho^\kappa_{\gamma_{n+1}}(b)W_{\gamma_{n+1}}^* W_{\gamma'_1 \cdots \gamma'_n}^*
\end{equation*}
gives rise to an embedding
$
{\cal F}_n \hookrightarrow {\cal F}_{n+1}, n \in \Zp.
$
Let
$\FHK$
be the 
$C^*$-subalgebra of $\OHK$ 
generated by
$\cup_{n=0}^\infty {\cal F}_n$.

We  define $2$-parameter unitary groups
on $\HK$ by setting
\begin{align*}
\theta^\rho_{r_1}: & u_\alpha\varphi_\eta(z) \in \HK \longrightarrow 
     e^{2{\pi}ir_1}u_\alpha\varphi_\eta(z) \in \HK, \quad \alpha \in \Sigma^\rho, z \in \BE\\
\theta^\eta_{r_2}: & v_a\varphi_\rho(w) \in \HK   \longrightarrow 
     e^{2{\pi}ir_2} v_a\varphi_\rho(w) \in \HK, \quad a \in \Sigma^\eta, w \in \BR
\end{align*}
for
$r_1, r_2 \in {\Bbb R}/{\Bbb Z} = {\Bbb T}. $
They extend on 
$
\HK\otimes_{\pi_1}\cdots \otimes_{\pi_{n-1}}\HK
$ 
for $(\pi_1,\dots,\pi_{n-1}) \in \Gamma_{n-1}$
by
\begin{equation*}
\theta^{\rho\otimes n}_{r_1}
  =
\theta^\rho_{r_1} \otimes_{\pi_1}\cdots 
\otimes_{\pi_{n-1}}\theta^\rho_{r_1}, \qquad
\theta^{\eta\otimes n}_{r_2}
 =
\theta^\eta_{r_2} \otimes_{\pi_1}\cdots 
\otimes_{\pi_{n-1}}\theta^\eta_{r_2},
\end{equation*}
which naturally extend on $F_n(\kappa)$.
We put
\begin{equation*}
u^\rho_{r_1}  = \sum_{n=0}^\infty \theta^{\rho\otimes n}_{r_1},
\qquad
 u^\eta_{r_2}  = \sum_{n=0}^\infty \theta^{\eta\otimes n}_{r_2}
 \qquad
\text{ on } F_\kappa = \oplus_{n=0}^\infty F_n(\kappa)
 \end{equation*}
for
$r_1, r_2 \in {\Bbb R}/{\Bbb Z} = {\Bbb T}$,
where
$\theta^{\rho\otimes n}_{r_1},
 \theta^{\eta\otimes n}_{r_2}
$
for $n=0$ are defined by 
$\theta^{\rho\otimes 0}_{r_1} (b_1 \oplus b_2) = e^{2{\pi}ir_1}b_1 \oplus b_2$,
$\theta^{\eta\otimes 0}_{r_2} (b_1 \oplus b_2) = b_1 \oplus e^{2{\pi}ir_2} b_2$.
We note that 
$u^\rho_{r_1}u^\eta_{r_2} = u^\eta_{r_2} u^\rho_{r_1}$.
Define 
\begin{equation*}
g_{(r_1,r_2)} = Ad(u^\rho_{r_1}u^\eta_{r_2})
\text{ on } F_\kappa
\text{ for } (r_1,r_2) \in {\Bbb T}^2.
\end{equation*}
The following lemma is straightforward.
\begin{lem}
For $ w \in \BR, z \in \BE, \alpha \in \Sigma^\rho, a \in \Sigma^\eta$
and
$(r_1,r_2) \in {\Bbb T}^2$,
we have
\begin{align*}
g_{(r_1,r_2)}(\bar{\phi}_\rho(w)) = \bar{\phi}_\rho(w), & \qquad
g_{(r_1,r_2)}(\bar{\phi}_\eta(z)) = \bar{\phi}_\eta(z), \\
g_{(r_1,r_2)}(s_\alpha)  = e^{2{\pi}ir_1} s_\alpha, & \qquad
g_{(r_1,r_2)}(t_a)  = e^{2{\pi}ir_2} t_a.
\end{align*}
\end{lem}
It is easy to see that
$g_{(r_1,r_2)}({\cal K}_\A(F_\kappa)) = {\cal K}_\A(F_\kappa)$
so that 
$g_{(r_1,r_2)}$ defines an automorphism on $\OHK$ for $(r_1,r_2) \in {\Bbb T}^2$,
which is still denoted by  
$g_{(r_1,r_2)}$.
The automorphisms
$g_{(r_1,r_2)}, (r_1,r_2) \in {\Bbb T}^2$ 
define an action
\begin{equation*}
g : (r_1,r_2) \in {\Bbb T}^2 \longrightarrow 
g_{(r_1,r_2)} \in \Aut(\OHK)
\end{equation*}
of ${\Bbb T}^2$,  called  the gauge action on 
$\OHK$.
Define 
a faithful  conditional expectation 
${\cal E}_{\HK}$
from
$\OHK$ onto the fixed point algebra 
$(\OHK)^g$
by setting
$$
{\cal E}_{\HK}(X) = 
\int_{(r_1,r_2) \in {\Bbb T}^2} g_{(r_1,r_2)}(X) \
dr_1 dr_2, \qquad X \in \OHK.
$$ 
Then the following lemma holds. 
Its proof is routine.
\begin{lem}
The fixed point algebra $(\OHK)^g$ of $\OHK$ under the action
$g$ of ${\Bbb T}^2$
coincides with  $\FHK$.
\end{lem}
We will  prove that the algebra
$\OHK$ has a universal property subject to the relations 
$(\HK)$.
Let us denote by
$\OU$
the universal $C^*$-algebra
generated by the operators
$ w\in \BR, z\in \BE $
and partial isometries
${\tt u}_\alpha,\alpha \in \Sigma^\rho,  
{\tt v}_a, a \in \Sigma^\eta
$
satisfying the following operator relations: 
\begin{align*}
\sum_{\beta \in \Sigma^\rho} {\tt u}_\beta {\tt u}_\beta^* +
& \sum_{b \in \Sigma^\eta} {\tt v}_b {\tt v}_b^* =1, \\ 
{\tt u}_\alpha {\tt u}_\alpha^* w  = w {\tt u}_\alpha {\tt u}_\alpha^*, 
& \qquad
{\tt v}_a {\tt v}_a^* w  = w {\tt v}_a {\tt v}_a^*, \\
{\tt u}_\alpha {\tt u}_\alpha^* z  = z {\tt u}_\alpha {\tt u}_\alpha^*, 
& \qquad
{\tt v}_a {\tt v}_a^* z = z {\tt v}_a {\tt v}_a^*, \\
\widehat{\rho}_\alpha(w)
=
{\tt u}_\alpha^* w {\tt u}_\alpha, 
&  \qquad
\widehat{\eta}_a(z)
=
{\tt v}_a^* z {\tt v}_a,\\ 
\widehat{\rho}^\eta_\alpha(z)
=
{\tt u}_\alpha^* z {\tt u}_\alpha, 
&  \qquad
\widehat{\eta}^\rho_a(w)
=
{\tt v}_a^* w {\tt v}_a,\\ 
\iota_\eta(y) & = \iota_\rho(y) 
\end{align*}
for
$
w \in \BR, z \in \BE,  \alpha \in \Sigma^\rho, a \in \Sigma^\eta,
y \in \A
$
where
$\iota_\eta: \A \hookrightarrow \BE$
and 
$\iota_\rho: \A \hookrightarrow \BR$
are natural embeddings.
The above ten relations of operators are also called the relations $(\HK)$.
Denote by 
$\BU$ the 
$C^*$-subalgebra of $\OU$
generated by the elements
$w \in \BR$
and
$z \in \BE$.
We set
for $\gamma \in \Sigma^\rho\cup\Sigma^\eta$
and $x \in \BU$ 
\begin{equation*}
\rho^{uni}_\gamma(x) 
 = w_\gamma^* x w_\gamma
 \quad
 \text{ where }
\quad
w_\gamma 
 =
{\begin{cases}
{\tt u}_\alpha & \text{ if } \gamma = \alpha \in \Sigma^\rho,\\
{\tt v}_a & \text{ if } \gamma = a \in \Sigma^\eta.
\end{cases}} 
\end{equation*}
Since
$
{\tt u}_\alpha^* \BU {\tt u}_\alpha \subset \BU
$
and
$
{\tt v}_a^* \BU {\tt v}_a \subset \BU,
$
we have
$$
\rho_\gamma^{uni}(\BU) \subset \BU,
$$
so that 
we have a family of endomorphisms
$\rho_\gamma^{uni}, \gamma \in \Sigma^{\rho\cup\eta}
$
on
$\BU$.
Similarly as in the preceding lemma, we have 
\begin{lem}
The triplet
$(\BU, \rho^{uni}, \Sigma^{\rho\cup\eta})$
is a $C^*$-symbolic dynamical system.
\end{lem}
For $n \in \Zp$,
denote by
${\cal F}_n^{uni}$
the $C^*$-subalgebra of $\OU$
generated by 
the operators
$w_{\gamma_1 \cdots \gamma_n} b w_{\gamma'_1 \cdots \gamma'_n}^*
$
for
$b \in \BU$
and
$\gamma_1 \cdots \gamma_n \in B_n(\Sigma^{\rho\cup\eta})$
such that
$ |\gamma_1 \cdots \gamma_n|_\rho = |\gamma'_1 \cdots \gamma'_n|_\rho$ 
and
$ |\gamma_1 \cdots \gamma_n|_\eta = |\gamma'_1 \cdots \gamma'_n|_\eta$.
Since
$
\sum_{\gamma \in \Sigma^{\rho\cup\eta}}
w_\gamma w_\gamma^* 
=1,
$
and
$w_\gamma w_\gamma^*$
commutes with
$\BU$,
the equality
\begin{equation*}
w_{\gamma_1 \cdots \gamma_n} b w_{\gamma'_1 \cdots \gamma'_n}^*
=
\sum_{\gamma_{n+1} \in \Sigma^{\rho\cup\eta}}
w_{\gamma_1 \cdots \gamma_n}w_{\gamma_{n+1}} 
\rho^{uni}_{\gamma_{n+1}}(b)w_{\gamma_{n+1}}^* w_{\gamma'_1 \cdots \gamma'_n}^*
\end{equation*}
gives rise to an embedding
$
{\cal F}_n^{uni} \hookrightarrow {\cal F}_{n+1}^{uni}, n \in \Zp.
$
Let
$\FU$
be the 
$C^*$-subalgebra of $\OU$ 
generated by
$\cup_{n=0}^\infty {\cal F}_n^{uni}$.
By the universality of the algebra
$\OU$
subject to the relations $(\HK)$,
the correspondences for each
$(r_1,r_2) \in {\Bbb T}^2$
\begin{align*}
w\in \OU\longrightarrow w\in \OU, & \qquad 
z\in \OU\longrightarrow z\in \OU,\\
{\tt u}_\alpha\in \OU\longrightarrow e^{2{\pi}ir_1} {\tt u}_\alpha \in \OU, & \qquad 
{\tt v}_a\in \OU\longrightarrow e^{2{\pi}ir_2} {\tt v}_a \in \OU,
\end{align*}
give rise to an automorphism of $\OU$,
which we denote by
$g_{(r_1,r_2)}^{uni}$.
Similarly to the preceding discussions,
$g^{uni}$ yields an action
of ${\Bbb T}^2$
on $\OU$, 
called the gauge action on $\OU$. 
Define similarly to the preceding discussions 
a faithful conditional expectation
${\cal E}_{\HK}^{uni}$
from
$\OU$ onto the fixed point algebra 
$(\OU)^{g^{uni}}$
by setting
$$
{\cal E}_{\HK}^{uni}(X) = 
\int_{(r_1,r_2) \in {\Bbb T}^2} g_{(r_1,r_2)}(X) \
dr_1 dr_2, \qquad X \in \OU.
$$
Similarly to the previous discussions, we have 
\begin{lem}
The fixed point algebra $(\OU)^{g^{uni}}$ of $\OU$ under the gauge action
$g^{uni}$ of ${\Bbb T}^2$
coincides with $\FU$.
\end{lem}
By the universality of the algebra $\OU$
subject to the relations $(\HK)$,
there exists a surjective $*$-homomorphism
$\varPsi: \OU \longrightarrow \OHK$
satisfying
\begin{align*}
\varPsi(w) =\Phi_\rho(w), & \qquad
\varPsi(z) =\Phi_\eta(z),\\
\varPsi({\tt u}_\alpha) = U_\alpha, & \qquad
\varPsi({\tt v}_a) = V_a
\end{align*}
for
$w \in \BR, z \in \BE, \alpha \in \Sigma^\rho, a \in \Sigma^\eta$.
We will prove that 
there exists a $*$-homomorphism
$\pi_\kappa: \BK \longrightarrow \BU$
such that
$ \pi_\kappa(\Phi_\rho(w)) = w, \pi_\kappa(\Phi_\eta(z)) = z$. 
\begin{lem} For $\alpha \in \Sigma^\rho, a \in \Sigma^\eta$, we have
\begin{enumerate}
\renewcommand{\labelenumi}{(\roman{enumi})}
\item 
The correspondence
$
U_\alpha \Phi_\eta(z) U_\alpha^* 
\in U_\alpha \Phi_\eta(\BE) U_\alpha^*
\longrightarrow 
{\tt u}_\alpha z {\tt u}_\alpha^* 
\in {\tt u}_\alpha \BE {\tt u}_\alpha^*
$ for $z \in \BE$
yields a $*$-homomorphism. 
\item
The correspondence
$
V_a \Phi_\rho(w) V_a^* 
\in V_a \Phi_\rho(\BR) V_a^*
\longrightarrow 
{\tt v}_a w {\tt v}_a^* 
\in {\tt v}_a \BR {\tt v}_a^*
$ for $w \in \BR$
yields a $*$-homomorphism. 
\end{enumerate}
\end{lem}
\begin{pf}
(i)
As
$U_\alpha^* U_\alpha = \Phi_\eta(\rho_\alpha(1)) = \Phi_\eta(P_\alpha)
$
and
${\tt u}_\alpha P_\alpha
= {\tt u}_\alpha \widehat{\rho}_\alpha(1) = {\tt u}_\alpha$,
the maps
\begin{align*}
Ad(U_\alpha^*): 
& U_\alpha \Phi_\eta(z) U_\alpha^* 
\in U_\alpha \Phi_\eta(\BE) U_\alpha^*
\longrightarrow 
U_\alpha ^* U_\alpha \Phi_\eta(z) U_\alpha^*U_\alpha
\in \Phi_\eta(P_\alpha \BE P_\alpha ), \\
Ad({\tt u}_\alpha): 
& 
P_\alpha z P_\alpha
\in P_\alpha \BE P_\alpha (\subset \BE) 
\longrightarrow 
{\tt u}_\alpha z {\tt u}_\alpha^* 
\in {\tt u}_\alpha \BE {\tt u}_\alpha^*
\end{align*}
are $*$-homomorphisms.
As
$\Phi_\eta : \BE \longrightarrow \Phi_\eta(\BE)$ is a $*$-isomorphism,
the desired map
$$
U_\alpha \Phi_\eta(z) U_\alpha^* 
\in U_\alpha \Phi_\eta(\BE) U_\alpha^*
\longrightarrow 
{\tt u}_\alpha z {\tt u}_\alpha^* 
\in {\tt u}_\alpha \BE {\tt u}_\alpha^*
$$
 for $z \in \BE$
yields a $*$-homomorphism.

(ii) is similar to (i).
\end{pf}
\begin{lem} \hspace{6cm}
\begin{enumerate}
\renewcommand{\labelenumi}{(\roman{enumi})}
\item 
For $\alpha \in \Sigma^\rho$,
the correspondence:
\begin{equation*}
\Phi_{\gamma_1}(x_{j_1})\Phi_{\gamma_2}(x_{j_2})\cdots
\Phi_{\gamma_n}(x_{j_n})U_\alpha U_\alpha^*
\in \BK U_\alpha U_\alpha^*
\longrightarrow
x_{j_1} x_{j_2}\cdots
x_{j_n} {\tt u}_\alpha {\tt u}_\alpha^*
\in \BU {\tt u}_\alpha {\tt u}_\alpha^*
\end{equation*}
for
$
x_{j_k} \in \BR (\gamma_k = \rho)
$
and
$
x_{j_k} \in \BE (\gamma_k = \eta)
$
gives rise to a $*$-homomorphism
from
$\BK U_\alpha U_\alpha^*$ to 
$\BU {\tt u}_\alpha {\tt u}_\alpha^*$.
\item
 For $a \in \Sigma^\eta$,
the correspondence:
\begin{equation*}
\Phi_{\gamma_1}(x_{j_1})\Phi_{\gamma_2}(x_{j_2})\cdots
\Phi_{\gamma_n}(x_{j_n})V_a V_a^*
\in \BK V_a V_a^*
\longrightarrow
x_{j_1} x_{j_2}\cdots
x_{j_n} {\tt v}_a {\tt v}_a^*
\in \BU {\tt v}_a {\tt v}_a^*
\end{equation*}
for
$
x_{j_k} \in \BR (\gamma_k = \rho)
$
and
$
x_{j_k} \in \BE (\gamma_k = \eta)
$
gives rise to a $*$-homomorphism
from
$\BK V_a V_a^*$ to 
$\BU {\tt v}_a {\tt v}_a^*$.
\end{enumerate}
\end{lem}
\begin{pf}
(i)
Since
$\widehat{\rho}_\alpha(w) \in \A \subset \BE$
for $w \in \BR$,
we see 
$\Phi_\rho(\widehat{\rho}_\alpha(w))
= \Phi_\eta(\widehat{\rho}_\alpha(w))
$
so that
$$
U_\alpha^* \Phi_\rho(w) U_\alpha
= \Phi_\rho(\widehat{\rho}_\alpha(w))
= \Phi_\eta(\widehat{\rho}_\alpha(w)),
\qquad
U_\alpha^* \Phi_\eta(z) U_\alpha
= \Phi_\eta(\widehat{\rho}^\eta_\alpha(z))
$$
for $w \in \BR, z \in \BE$.
For 
$
x_{j_k} \in \BE
$
or
$\in \BR$,
put
$$
\widehat{x}_{j_k} 
: =
\begin{cases}
\widehat{\rho}_\alpha(x_{j_k}) & \text{ if } x_{j_k} \in \BR,\\
\widehat{\rho}^\eta_\alpha(x_{j_k}) & \text{ if } x_{j_k} \in \BE.
\end{cases}
$$
We then have $\widehat{x}_{j_k}  \in \BE$ so that
\begin{align*}
& \Phi_{\gamma_1}(x_{j_1})\Phi_{\gamma_2}(x_{j_2})\cdots
  \Phi_{\gamma_n}(x_{j_n})U_\alpha U_\alpha^* \\
=&  
U_\alpha U_\alpha^* \Phi_{\gamma_1}(x_{j_1})
U_\alpha U_\alpha^* \Phi_{\gamma_2}(x_{j_2})U_\alpha U_\alpha^*
\cdots
U_\alpha U_\alpha^* \Phi_{\gamma_n}(x_{j_n})U_\alpha U_\alpha^* \\
=&  
U_\alpha \Phi_\eta( \widehat{x}_{j_1})
         \Phi_\eta(\widehat{x}_{j_2})
\cdots
         \Phi_\eta( \widehat{x}_{j_n}) U_\alpha^* \\ 
=&  
U_\alpha \Phi_\eta( \widehat{x}_{j_1}\widehat{x}_{j_2}
\cdots \widehat{x}_{j_n}) U_\alpha^*.
\end{align*}
By the preceding lemma,
the correspondence
$$
U_\alpha \Phi_\eta( \widehat{x}_{j_1}\widehat{x}_{j_2}
\cdots \widehat{x}_{j_n}) U_\alpha^* 
\in U_\alpha \Phi_\eta(\BE) U_\alpha^*
\longrightarrow
{\tt u}_\alpha  \widehat{x}_{j_1}\widehat{x}_{j_2}
\cdots \widehat{x}_{j_n}  {\tt u}_\alpha^*
\in {\tt u}_\alpha \BE  {\tt u}_\alpha^*
$$
gives rise to a $*$-homomorphism
from
$
 U_\alpha \Phi_\eta(\BE) U_\alpha^*
$ to 
$ {\tt u}_\alpha \BE  {\tt u}_\alpha^*$.
Since we have
\begin{align*}
{\tt u}_\alpha  \widehat{x}_{j_1}\widehat{x}_{j_2}
\cdots \widehat{x}_{j_n}  {\tt u}_\alpha^*
& =
{\tt u}_\alpha {\tt u}_\alpha^* x_{j_1} 
{\tt u}_\alpha {\tt u}_\alpha^* x_{j_2} {\tt u}_\alpha {\tt u}_\alpha^*
\cdots
{\tt u}_\alpha {\tt u}_\alpha^* x_{j_n} {\tt u}_\alpha {\tt u}_\alpha^* \\
& =
x_{j_1} x_{j_2} 
\cdots
x_{j_n} {\tt u}_\alpha {\tt u}_\alpha^*, 
\end{align*}
we have a desired $*$-homomorphism
from
$\BK U_\alpha U_\alpha^*$ to 
$\BU {\tt u}_\alpha {\tt u}_\alpha^*$.

(ii) is similar to (i).
\end{pf}
The above $*$-homomorphisms of  (i)  and of  (ii) are denoted by 
\begin{align*}
\pi_\alpha : 
& \Phi_{\gamma_1}(x_{j_1})\Phi_{\gamma_2}(x_{j_2})\cdots
\Phi_{\gamma_n}(x_{j_n})U_\alpha U_\alpha^*
\in \BK U_\alpha U_\alpha^*
\longrightarrow
x_{j_1} x_{j_2}\cdots
x_{j_n} {\tt u}_\alpha {\tt u}_\alpha^*
\in \BU {\tt u}_\alpha {\tt u}_\alpha^* \\
\pi_a : 
& \Phi_{\gamma_1}(x_{j_1})\Phi_{\gamma_2}(x_{j_2})\cdots
\Phi_{\gamma_n}(x_{j_n})V_a V_a^*
\in \BK V_a V_a^*
\longrightarrow
x_{j_1} x_{j_2}\cdots
x_{j_n} {\tt v}_a {\tt v}_a^*
\in \BU {\tt v}_a {\tt v}_a^* 
\end{align*}
\begin{lem}
There exists a $*$-homomorphism
\begin{equation*}
\pi_\kappa : \BK \longrightarrow \BU
\end{equation*}
such that
$\pi_\kappa(\Phi_\rho(w)) = w$ for $w \in \BR$
and
$\pi_\kappa(\Phi_\eta(z)) = z$ for $z \in \BE$.
\end{lem}
\begin{pf}
Since
$
\sum_{\alpha \in \Sigma^\rho} U_\alpha U_\alpha^* + 
\sum_{a \in \Sigma^\eta} V_a V_a^* =1
$
and
$
\sum_{\alpha \in \Sigma^\rho} {\tt u}_\alpha {\tt u}_\alpha^* + 
\sum_{a \in \Sigma^\eta} {\tt v}_a {\tt v}_a^* =1
$
by putting
\begin{equation*}
\pi_\kappa(X) := 
\sum_{\alpha \in \Sigma^\rho} \pi_\alpha(X U_\alpha U_\alpha^*) {\tt u}_\alpha {\tt u}_\alpha^* 
+ 
\sum_{a \in \Sigma^\eta} \pi_a(X V_a V_a^*){\tt v}_a {\tt v}_a^*
\end{equation*}
for
$X \in \BK$,
we have a desired $*$-homomorphism
from
$\BK$
to
$\BU$.
\end{pf}
We will consider 
the $*$-homomorphism 
$\varPsi: \OU \longrightarrow \OHK$
again.
By the above discussions we have
\begin{lem}
The restriction 
$\varPsi|_\BU : \BU \longrightarrow \BK$
of 
$\varPsi$ to the subalgebra $\BU$ 
is the inverse of
$\pi_\kappa: \BK \longrightarrow \BU$.
Hence 
$\varPsi|_\BU : \BU \longrightarrow \BK$
is a $*$-isomorphism.
 \end{lem}
Therefore we reach the main result of the paper:
\begin{thm} \label{thm:main}
The $C^*$-algebra $\OHK$
associated with the 
Hilbert $C^*$-quad module $\HK$
is canonically $*$-isomorphic to the universal $C^*$-algebra
$\OU$
generated by the operators
$w\in \BR, z\in \BE $
and partial isometries
${\tt u}_\alpha,\alpha \in \Sigma^\rho,  
{\tt v}_a, a \in \Sigma^\eta
$
satisfying the operator relations: 
\begin{align}
\sum_{\beta \in \Sigma^\rho} {\tt u}_\beta {\tt u}_\beta^* +
& \sum_{b \in \Sigma^\eta} {\tt v}_b {\tt v}_b^* =1, \label{eqn:M1} \\ 
{\tt u}_\alpha {\tt u}_\alpha^* w  = w {\tt u}_\alpha {\tt u}_\alpha^*, 
& \qquad
{\tt v}_a {\tt v}_a^* w  = w {\tt v}_a {\tt v}_a^*, \label{eqn:M2} \\
{\tt u}_\alpha {\tt u}_\alpha^* z  = z {\tt u}_\alpha {\tt u}_\alpha^*, 
& \qquad
{\tt v}_a {\tt v}_a^* z = z {\tt v}_a {\tt v}_a^*, \label{eqn:M3} \\
\widehat{\rho}_\alpha(w)
=
{\tt u}_\alpha^* w {\tt u}_\alpha, 
&  \qquad
\widehat{\eta}_a(z)
=
{\tt v}_a^* z {\tt v}_a, \label{eqn:M4} \\ 
\widehat{\rho}^\eta_\alpha(z)
=
{\tt u}_\alpha^* z {\tt u}_\alpha, 
&  \qquad
\widehat{\eta}^\rho_a(w)
=
{\tt v}_a^* w {\tt v}_a, \label{eqn:M5} \\ 
\iota_\eta(y) & = \iota_\rho(y) \label{eqn:M6} 
\end{align}
for
$
w \in \BR, z \in \BE, \alpha \in \Sigma^\rho, a \in \Sigma^\eta,
y \in \A$
where
$\iota_\rho: \A \hookrightarrow \BR$
and 
$\iota_\eta: \A \hookrightarrow \BE$
are natural embeddings.
\end{thm}
\begin{pf}
The triplets
$(\BU, \rho^{uni}, \Sigma^{\rho\cup\eta})$
and
$(\BK, \rho^\kappa, \Sigma^{\rho\cup\eta})$
are both
the $C^*$-symbolic dynamical systems.
As in the discussions of the proof of
\cite[Lemma 3.2]{MaMZ2010},
the above lemma implies that
the restriction 
$\varPsi|_\FU : \FU \longrightarrow \FHK$
of 
$\varPsi$ to the subalgebra $\FU$ 
is a $*$-isomorphism.
The diagram:
\begin{equation*}
\begin{CD}
\OU @>\varPsi>> \OHK \\
@V{{\cal E}^{uni}_{\HK}}VV @VV{{\cal E}_{\HK}}V \\
\FU @>>{\varPsi|_{\FU}}> \FHK
\end{CD}
\end{equation*}
is commutative.
Since
the conditional expectation
${\cal E}^{uni}_{\HK}$ is faithful
and
the restriction  
$\varPsi|_{\FU}$ is $*$-isomorphic,
one concludes that
$\varPsi: \OU \longrightarrow \OHK$
is a $*$-isomorphism
by a routine argument 
as in \cite[2.9.Proposition]{CK}. 
\end{pf}
The above theorem implies the following:
Suppose that there exist two families of partial isometries
$
\hat{u}_\alpha,\alpha \in \Sigma^\rho,  
\hat{v}_a, a \in \Sigma^\eta
$
in  a unital $C^*$-algebra ${\cal D}$  
and two $*$-homomorphisms
$\pi_\rho: \BR \longrightarrow {\cal D}$,
$\pi_\eta: \BE \longrightarrow {\cal D}$
satisfying the relations: 
\begin{align*}
\sum_{\beta \in \Sigma^\rho} \hat{u}_\beta \hat{u}_\beta^* +
& \sum_{b \in \Sigma^\eta} \hat{v}_b \hat{v}_b^* =1, \label{eqn:rep1}\\ 
\hat{u}_\alpha \hat{u}_\alpha^* \pi_\rho(w)  
= \pi_\rho(w) \hat{u}_\alpha \hat{u}_\alpha^*, 
& \qquad
\hat{v}_a \hat{v}_a^* \pi_\rho(w)  = \pi_\rho(w) \hat{v}_a \hat{v}_a^*, \label{eqn:rep2} \\
\hat{u}_\alpha \hat{u}_\alpha^* \pi_\eta(z)  
= \pi_\eta(z) \hat{u}_\alpha \hat{u}_\alpha^*, 
& \qquad
\hat{v}_a \hat{v}_a^* \pi_\eta(z) = \pi_\eta(z) \hat{v}_a \hat{v}_a^*, \label{eqn:rep3} \\
\pi_\rho(\widehat{\rho}_\alpha(w))
=
\hat{u}_\alpha^* \pi_\rho(w) \hat{u}_\alpha, 
&  \qquad
\pi_\eta(\widehat{\eta}_a(z))
=
\hat{v}_a^* \pi_\rho(z) \hat{v}_a, \label{eqn:rep4} \\ 
\pi_\eta(\widehat{\rho}^\eta_\alpha(z))
=
\hat{u}_\alpha^* \pi_\eta(z) \hat{u}_\alpha, 
&  \qquad
\pi_\rho(\widehat{\eta}^\rho_a(w))
=
\hat{v}_a^* \pi_\rho(w) \hat{v}_a, \label{eqn:rep5}\\ 
\pi_\rho(y) & = \pi_\eta(y)  \label{eqn:rep6}
\end{align*}
for 
$w \in \BR, z \in \BE, \alpha \in \Sigma^\rho, a \in \Sigma^\eta, y \in \A$.
Then there exists a $*$-homomorphism
$\pi$ from $\OHK$ onto 
the $C^*$-algebra generated by
$
\pi_\rho(w),  \pi_\eta(z),
 \hat{u}_\alpha,  \hat{v}_a
$
such that
$\pi(w) = \pi_\rho(w),  \pi(z) = \pi_\eta(z),
\pi(U_\alpha) = \hat{u}_\alpha,  
\pi(V_a) =\hat{v}_a
$
for 
$w \in \BR, z \in \BE, \alpha \in \Sigma^\rho, a \in \Sigma^\eta$.

The notions of condition (I) and irreducibility  for $C^*$-symbolic dynamical systems
have been defined in 
\cite{MaContem}(cf. \cite{MaMZ2010}).
The former guarantees a uniqueness of the resulting $C^*$-algebras under operator
relations among generators. 
The latter does a simplicity of them.

\noindent
{\bf Definition.}
We say that 
$\HK$ satisfies {\it condition\/} (I) 
if the $C^*$-symbolic dynamical system 
$(\BK, \rho^\kappa, \Sigma^{\rho\cup\eta})$
satisfies condition (I). 
 We also say that
$\HK$ is {\it irreducible\/} 
if the $C^*$-symbolic dynamical system 
$(\BK, \rho^\kappa, \Sigma^{\rho\cup\eta})$
is irreducible.
\begin{cor}\label{cor:simplicity}
Suppose that
$\A$ is commutative.
\begin{enumerate}
\renewcommand{\labelenumi}{(\roman{enumi})}
\item 
If
$\HK$ satisfies {\it condition } (I),
then the $C^*$-algebra
$\OHK$ is the unique $C^*$-algebra 
subject to the relations 
$(\HK)$.
\item
If in adittion
$\HK$ is irreducible,
the $C^*$-algebra $\OHK$ is simple.
\end{enumerate}
\end{cor}

In the rest of this section,
we will consider a $C^*$-subalgebra of $\OHK$ generated by 
$U_\alpha, V_a$ 
for $ \alpha \in \Sigma^\rho, a \in \Sigma^\eta$
 and elements
$\Phi_\rho(x), \Phi_\eta(x)$
for $ x \in \A$.
 The subalgebra is denoted by 
 $\OAK$.
Since the operators 
$\bar{\phi}_\rho(x)$ and 
$\bar{\phi}_\eta(x)$ for $x \in \A$
are different only on
$F_0(\kappa)$,
we have
$\Phi_\rho(x)=\Phi_\eta(x)$ for $x \in \A$,
which we denote 
by $\Phi(x)$.
Hence the  following relations hold:
\begin{align}
\sum_{\beta \in \Sigma^\rho} U_\beta U_\beta^* &  +
\sum_{b \in \Sigma^\eta} V_b V_b^* =1,  \\
U_\alpha U_\alpha^* \Phi(x)  =\Phi( x) U_\alpha U_\alpha^*,& \qquad
V_a V_a^* \Phi(x)  = \Phi(x) V_a V_a^*, \\
U_\alpha^* \Phi(x) U_\alpha = \Phi(\rho_\alpha(x)),& \qquad
V_a^* \Phi(x) V_a = \Phi(\eta_a(x))
\end{align}
for $x \in \A, \alpha \in \Sigma^\rho, a\in \Sigma^\eta$. 
 For $\omega =(\alpha,b,a,\beta) \in \SK$,
define the operators in $\OAK$ by
\begin{align}
C_\omega 
= & V_a U_\beta V_b^* U_\alpha^*, \\
P_{r,\omega} =  V_a U_\beta U_\beta^* V_a^*, & \qquad 
P_{s,\omega} =  U_\alpha V_b V_b^* U_\alpha^*.
\end{align}  
\begin{lem}
 For $\omega =(\alpha,b,a,\beta) \in \SK$,
 we have
\begin{enumerate}
\renewcommand{\labelenumi}{(\roman{enumi})}
\item 
$C_\omega$ 
is a partial isometry in $\OHK$ satisfying
\begin{align}
C_\omega \Phi(x) =&  \Phi(x) C_\omega, \qquad x \in \A, \label{eqn:CPhi}\\ 
V_a U_\beta = C_\omega U_\alpha V_b,& \qquad
U_\alpha V_b = C_\omega^* V_a U_\beta.
\end{align}
\item
$P_{r,\omega}, P_{s,\omega}$
are projections satisfying
\begin{equation}
P_{r,\omega} = C_\omega C_\omega^*, \qquad
P_{s,\omega} = C_\omega^* C_\omega. \label{eqn:PC}
\end{equation}
\item 
$C_\omega \not\in \BK$.
\end{enumerate}
\end{lem}
\begin{pf}
(i)
By (5.10) and (5.11), we have for $x \in \A$,
\begin{align*}
C_\omega \Phi(x)
= & V_a U_\beta V_b^* \Phi(\rho_\alpha(x)) U_\alpha^* 
=  V_a U_\beta \Phi(\eta_b(\rho_\alpha (x)) V_b^* U_\alpha^* \\
\Phi(x )  C_\omega
= &V_a  \Phi(\eta_a(x)) U_\beta V_b^* U_\alpha^*
= V_a U_\beta \Phi(\rho_\beta ( \eta_a(x) )V_b^* U_\alpha^*. 
\end{align*}  
By \eqref{eqn:kappa}, we have \eqref{eqn:CPhi}.
We also have
\begin{equation*}
C_\omega U_\alpha V_b
=  V_a U_\beta \Phi(\eta_b(\rho_\alpha(1))) 
=  V_a U_\beta \Phi(\rho_\beta( \eta_a(1))) 
=  V_a U_\beta U_\beta^* V_a^* V_a U_\beta 
=  V_a U_\beta
\end{equation*}  
and similarly
$
C_\omega^* V_a U_\beta 
=  U_\alpha V_b.  
$  

(ii)
The equaltions 
\eqref{eqn:PC}
are straightforward.

(iii)
Suppose that
$C_\omega \in \BK$.
Since
$U_\alpha U_\alpha^*$ commutes with
$\Phi_\rho(w), \Phi_\eta(z)$
for
$w \in \BR, z \in \BE$,
it commutes with
$\BK$ and hence with
$C_\omega$.
Hence we have
\begin{equation*}
V_a U_\beta = C_\omega U_\alpha V_b = U_\alpha U_\alpha^* C_\omega U_\alpha V_b.
\end{equation*}
As
$V_a V_a^* U_\alpha U_\alpha^*=0$,
one has
$V_a U_\beta =0$.
Since
$\Phi(\rho_\beta(\eta_a(y))) = U_\beta^* V_a^* y V_a U_\beta =0$,
one has
$\rho_\beta(\eta_a(y)) =0$ for all $y \in \A$,
a contradiction.
\end{pf}

%
For the two $C^*$-symbolic dynamical systems
$(\A, \rho,\Sigma^\rho)$
and
$(\A, \eta,\Sigma^\eta)$
in the $C^*$-textile dynamical system
$\CTDS$,
we define their  union
$(\A, \rho\cup\eta,\Sigma^{\rho\cup\eta})$
as the following way, where
$
\Sigma^{\rho\cup\eta} = \Sigma^\rho \cup \Sigma^\eta.
$
For
$\gamma \in \Sigma^{\rho\cup\eta}$, 
define an endomorphism
$(\rho\cup\eta)_\gamma$ on $\A$ by setting
\begin{equation*}
 (\rho\cup\eta)_\gamma
 =
 \begin{cases}
 \rho_\gamma & \text{ if } \gamma \in \Sigma^\rho,\\
 \eta_\gamma & \text{ if } \gamma \in \Sigma^\eta.
 \end{cases}
\end{equation*}
It is easy to see that 
the triplet
$(\A, \rho\cup\eta,\Sigma^{\rho\cup\eta})$
is a $C^*$-symbolic dynamical system.
Hence we  have a $C^*$-algebra 
${\cal O}_{\rho\cup\eta}$
from 
$(\A, \rho\cup\eta,\Sigma^{\rho\cup\eta})$.
Denote by
${\tt S}_\alpha, \alpha\in \Sigma^\rho$
and
${\tt S}_a, a\in \Sigma^\eta$
its generating partial isometries satisfying the relations:
\begin{gather}
\sum_{\beta \in \Sigma^\rho}  {\tt S}_\beta {\tt S}_\beta^* +
\sum_{b \in \Sigma^\eta}  {\tt S}_b {\tt S}_b^* =1, \\
{\tt S}_\alpha {\tt S}_\alpha^* x = x {\tt S}_\alpha {\tt S}_\alpha^*,
\qquad  
{\tt S}_a {\tt S}_a^* x = x {\tt S}_a {\tt S}_a^*,\\ 
{\tt S}_\alpha^* x {\tt S}_\alpha = \rho_\alpha(x),   \qquad
{\tt S}_a^* x {\tt S}_a = \eta_a(x)
\end{gather}
for 
$x \in \A, \alpha\in \Sigma^\rho, a \in \Sigma^\eta$.
Then we have
\begin{prop}
\begin{enumerate}
\renewcommand{\labelenumi}{(\roman{enumi})}
\item 
The correspondence:
\begin{equation*}
{\tt S}_\alpha  \longleftrightarrow U_\alpha  
\qquad 
{\tt S}_a  \longleftrightarrow V_a,
\qquad
x \longleftrightarrow \Phi(x)  
\end{equation*} 
for $\alpha \in \Sigma^\rho,  a \in \Sigma^\eta,
x \in \A$
yield a $*$-isomorphism between 
the $C^*$-algebras 
${\cal O}_{\rho\cup\eta}$
and
$\OAK$. 
Therefore we have
\begin{equation}
\OAK \cong {\cal O}_{\rho\cup\eta},
\end{equation}
\item 
Hence the $C^*$-algebra $\OAK$ is realized as the universal $C^*$-algebra
generated by partial isometries $U_\alpha, V_a$ for $\alpha \in \Sigma^\rho, a \in \Sigma^\eta$
and elements $x \in \A$ subject to the relations:
\begin{align*}
\sum_{\beta \in \Sigma^\rho} U_\beta U_\beta^* +
& \sum_{b \in \Sigma^\eta} V_b V_b^* =1, \\ 
U_\alpha U_\alpha^* \Phi(x)  = \Phi(x) U_\alpha U_\alpha^*, & \qquad
V_a V_a^* \Phi(x)  = \Phi(x ) V_a V_a^*,\\
U_\alpha^* \Phi(x) U_\alpha = \Phi(\rho_\alpha(x)), &  \qquad
V_a^* \Phi(x) V_a = \Phi(\eta_a(x)).
\end{align*}
 \end{enumerate}
\end{prop}
\begin{pf}
By the universality of the algebra
$ {\cal O}_{\rho\cup\eta}$
subject to the relation (5.17), (5.18), (5.19),
the correspondences 
\begin{equation*}
{\tt S}_\alpha  \longrightarrow U_\alpha  
\qquad 
{\tt S}_a  \longrightarrow V_a,
\qquad
x \longrightarrow \Phi(x)  
\end{equation*} 
for $\alpha \in \Sigma^\rho,  a \in \Sigma^\eta,
x \in \A$
yield a $*$-homomorphism from 
${\cal O}_{\rho\cup\eta}$
to
$\OAK$
which we denote by
$\Psi$. 
Since the action 
$g_{(r_1,r_2)}, (r_1,r_2)\in {\Bbb T}^2$
preserves $\OAK$,
the automorphisms
$g_{(t, t)}$
for
$t \in {\Bbb T}$
 give rise to an action on $\OAK$
 which we denote by 
 $g_t^\A$.
 It satisfies
 \begin{equation*}
 g^\A_t(\Phi(x)) = \Phi(x),
\qquad
   g^\A_t(U_\alpha) = e^{2\pi i t} U_\alpha,
 \qquad
  g^\A_t(V_a) = e^{2\pi i t} V_a
\end{equation*}
for $\alpha \in \Sigma^\rho,  a \in \Sigma^\eta,
x \in \A$ and $t \in {\Bbb T}$.
Let $\hat{g}$ be the
gauge action on 
${\cal O}_{\rho\cup\eta}$
which satisfies
 \begin{equation*}
 \hat{g}_t(x) = x,
\qquad
  \hat{g}_t({\tt S}_\alpha) = e^{2\pi i t}{\tt S}_\alpha,
 \qquad
  \hat{g}_t({\tt S}_a) = e^{2\pi i t} {\tt S}_a
\end{equation*}
for $\alpha \in \Sigma^\rho,  a \in \Sigma^\eta,
x \in \A$ and $t \in {\Bbb T}$.
Hence we have
\begin{equation*}
\Psi \circ \hat{g}_t = g^\A_t \circ \Psi
\qquad \text{ for }\quad
 t \in {\Bbb T}.
\end{equation*}
Let
$(\OAK)^{g^\A}$ 
be the fixed point algebra of
$\OAK$
under the action
$g^\A$.
Denote by 
${\cal E}^\A : \OAK \longrightarrow  
(\OAK)^{g^\A}$
the conditional expectation defined by the formula:
$$
{\cal E}^\A(X) = \int_{{\Bbb T}} g^\A_t(X) \ dt 
\qquad 
\text{ for }
X \in \OAK.
$$ 
Denote by
${\cal E}^{\rho \cup \eta} : {\cal O}_{\rho\cup\eta} \longrightarrow 
  ({\cal O}_{\rho\cup\eta})^{\hat{g}}
$
the conditional expectation similarly defined to the above by the gauge action
$\hat{g}$.
Since we have
\begin{equation*}
\Psi \circ {\cal E}^{\rho\cup\eta}  = {\cal E}^\A \circ \Psi
\end{equation*}
 and
 $
 {\cal E}^{\rho\cup\eta}
$
is faithful,
by a routine argument as in \cite[2.9.Proposition]{CK},
one concludes that $\Psi$ is injective and hence $*$-isomorphic.
\end{pf}

\section{ Relation between $\OHK$ and $\ORE$}
For a $C^*$-textile dynamical system
$\CTDS$,
the author has introduced a $C^*$-algebra
$\ORE$ in \cite{MaPre2011}.
It is realized as the
universal $C^*$-algebra 
$C^*(x, S_\alpha, T_a ; x \in \A, \alpha \in \Sigma^\rho, 
a \in \Sigma^\eta)$  
generated by 
$x \in \A$ 
and two families of partial isometries 
$S_{\alpha}, \alpha \in \Sigma^\rho$,
$T_a, a \in \Sigma^\eta$
subject to the relations 
\eqref{eqn:S}, \eqref{eqn:T} and \eqref{eqn:STTS}
that are
called the relations $(\rho,\eta;\kappa)$.
In this section we will describe a relationship  between the two algebras
$\OHK$ and $\ORE$.
As  both of the $C^*$-algebras $\ORHO$ and $\OETA$
are naturally regarded as $C^*$-subalgebras of   
$\ORE$,
the algebras  
$\BR$ and $\BE$ 
may be realized as 
the $C^*$-subalgebra of $\ORE$
generated by $S_\alpha x S_\alpha^*$
 for $x \in \A, \alpha \in \Sigma^\rho$
and that of   $\ORE$
generated by $T_a x T_a^*$
 for $x \in \A, a \in \Sigma^\eta$
respectively.
\begin{lem}
Let
$S_{\alpha}, \alpha \in \Sigma^\rho$
and
$T_a, a \in \Sigma^\eta$
be partial isometries 
in the algebra $\ORE$
satisfying the relations 
\eqref{eqn:S}, \eqref{eqn:T} and \eqref{eqn:STTS}.
\begin{enumerate}
\renewcommand{\labelenumi}{(\roman{enumi})}
\item 
For $\alpha \in \Sigma^\rho$
and
$z = \sum_{b \in \Sigma^\eta} T_b z_b T_b^* \in \BE$ 
as in \eqref{eqn:z},
\begin{equation}
S_\alpha^* z S_\alpha
= 
\sum
\begin{Sb}
b,a,\beta\\
(\alpha,b,a,\beta)\in\SK
\end{Sb}
T_b \rho_\beta(z_a) T_b^* 
= \widehat{\rho}_\alpha^\eta(z). 
\end{equation}
\item 
For $a \in \Sigma^\eta$
and
$w = \sum_{\beta \in \Sigma^\rho} S_\beta w_\beta S_\beta^* \in \BR$ 
as in \eqref{eqn:w},
\begin{equation}
T_a^* w T_a=
\sum
\begin{Sb}
\alpha, b,\beta\\
(\alpha,b,a,\beta)\in\SK
\end{Sb}
S_\beta \eta_b(w_\alpha) S_\beta^* 
=\widehat{\eta}_a^\rho(w).
\end{equation}
\end{enumerate}
\end{lem}
\begin{pf}
(i)
By \cite[Lemma 4.2]{MaPre2011} 
the following formulae hold
for $\alpha \in \Sigma^\rho, a \in \Sigma^\eta$,
\begin{equation*}  
T_a^* S_\alpha 
=
\sum_{
\stackrel{b,\beta}{\kappa(\alpha, b) = (a,\beta)}}
 S_\beta \eta_b(\rho_\alpha(1))T_b^*, 
 \qquad  
S_\alpha^* T_a 
=
\sum_{
\stackrel{b,\beta}{\kappa(\alpha, b) = (a,\beta)}}
 T_b \rho_\beta(\eta_a(1))S_\beta^*.
\end{equation*} 
It then follows that
\begin{align*}
S_\alpha^* z S_\alpha
& = 
\sum_{a \in \Sigma^\eta} S_\alpha^* T_a z_a T_a^* S_\alpha \\
& = 
\sum_{a \in \Sigma^\eta}
\sum_{
\stackrel{b,\beta}{\kappa(\alpha, b) = (a,\beta)}}
 T_b \rho_\beta(Q_a)S_\beta^*
z_a 
\sum_{
\stackrel{b',\beta'}{\kappa(\alpha, b') = (a,\beta')}}
 S_{\beta'} \rho_{\beta'}(Q_a)T_{b'}^* \\ 
& = 
\sum
\begin{Sb}
b,a,\beta\\
(\alpha,b,a,\beta)\in\SK
\end{Sb}
T_b \rho_\beta(Q_a z_aQ_a) T_b^* 
= \widehat{\rho}_\alpha^\eta(z).
\end{align*}

(ii) is similar to (i).
\end{pf}
Let us denote by 
$S_1, S_2$ isometries satisfying 
$S_1S_1^* + S_2S_2^* =1$.
The $C^*$-algebra generated by them is the Cuntz algebra ${\cal O}_2$ of order $2$.
In  the tensor product $C^*$-algebra  
$\ORE \otimes {\cal O}_2$,
we put
\begin{align*}
\widehat{u}_\alpha = S_\alpha \otimes S_1 \text{ for } \alpha \in \Sigma^\rho, 
& \qquad
\widehat{v}_a = T_a \otimes S_2 \text{ for } a \in \Sigma^\eta, \\ 
\pi_\rho(w) = w \otimes 1  \text{ for } w \in \BR, 
& \qquad
\pi_\eta(z) = z \otimes 1  \text{ for } z  \in \BE.
\end{align*}
By the above lemma, the following lemma is straightforward.
\begin{lem}
The operators 
$\widehat{u}_\alpha,
\widehat{v}_a
$
for
$ \alpha \in \Sigma^\rho,
a \in \Sigma^\eta
$
and
$\pi_\rho(w),
\pi_\eta(z)
$
for
$ w \in \BR, 
 z  \in \BE
 $
 satisfy the relations $(\HK)$.
\end{lem}
We thus have
\begin{thm}
Assume that
$\CTDS$
satisfies condition (I).
Then the correspondences
\begin{align*}
U_\alpha \in \OHK & \longrightarrow  S_\alpha \otimes S_1 \in \ORE\otimes {\cal O}_2,\\  
V_a \in \OHK        & \longrightarrow  T_a \otimes S_2 \in \ORE\otimes {\cal O}_2, \\ 
w \in \BR            & \longrightarrow  w \otimes 1  \in \ORE\otimes {\cal O}_2, \\ 
z \in \BE             & \longrightarrow  z \otimes 1  \in \ORE\otimes {\cal O}_2     
\end{align*}
 for
$\alpha \in \Sigma^\rho,
 a \in \Sigma^\eta
$ 
give rise to a $*$-isomorphism from $\OHK$ to the $C^*$-subalgebra
of $\ORE\otimes {\cal O}_2$
genarated by the partial isometries 
$S_\alpha \otimes S_1,
T_a \otimes S_2 
$
 for
$\alpha \in \Sigma^\rho,
 a \in \Sigma^\eta
$ 
and the elements
$
w\otimes 1,
z \otimes 1
$
for
$
w \in \BR,
z \in \BE$.
That is
\begin{equation*}
\OHK \cong 
C^*( S_\alpha \otimes S_1, T_a \otimes S_2,
w\otimes 1, z \otimes 1:
 \alpha \in \Sigma^\rho,
 a \in \Sigma^\eta,
w \in \BR,
z \in \BE ).
\end{equation*}
\end{thm}
We will present an example.
Let 
$\boldsymbol\alpha$,
$\boldsymbol\beta$
be automorphisms of a unital commutative 
$C^*$-algebra $\A$.
Put
$\Sigma^\rho = \{ \boldsymbol\alpha \},
\Sigma^\eta = \{ \boldsymbol\beta \}
$
and define
$\rho_{\boldsymbol\alpha} = \boldsymbol\alpha,
\eta_{\boldsymbol\beta} = \boldsymbol\beta.
$
We have two $C^*$-symbolic dynamical systems
$
(\A, \boldsymbol\alpha, \{ \boldsymbol\alpha \}),
(\A, \boldsymbol\beta, \{ \boldsymbol\beta \}).
$
Assume that 
$
\boldsymbol\alpha \circ \boldsymbol\beta
=
 \boldsymbol\beta \circ \boldsymbol\alpha.
$
Put
$
\Sigma^{\boldsymbol\alpha \boldsymbol\beta}
= \{ (\boldsymbol\alpha, \boldsymbol\beta) \},
\Sigma^{\boldsymbol\beta \boldsymbol\alpha}
= \{ (\boldsymbol\beta, \boldsymbol\alpha) \}.
$
The specification
$\kappa:
\Sigma^{\boldsymbol\alpha \boldsymbol\beta}
\longrightarrow
\Sigma^{\boldsymbol\beta \boldsymbol\alpha}
$
is unique and satisfies 
$\kappa(\boldsymbol\alpha, \boldsymbol\beta)
= (\boldsymbol\beta, \boldsymbol\alpha).
$
We have a $C^*$-textile dynamical system
$(\A, \boldsymbol\alpha,  \boldsymbol\beta,
\{ \boldsymbol\alpha \}, \{ \boldsymbol\beta \},
\kappa ).
$
We denote by
${\cal H}_\kappa^{
\boldsymbol\alpha, \boldsymbol\beta}
$ 
the associated Hilbert $C^*$-quad module.
Since 
$\SK = \{
 (\boldsymbol\alpha, \boldsymbol\beta,\boldsymbol\beta,
\boldsymbol\alpha) \}$
a singleton  and 
$E_\omega = \boldsymbol\beta(\boldsymbol\alpha(1)) = 1$
 so that
${\cal H}_\kappa^{
\boldsymbol\alpha, \boldsymbol\beta}
 =  \A.$
As
$\boldsymbol\alpha \circ \boldsymbol\beta
= \boldsymbol\beta\circ \boldsymbol\alpha,
$
they induce an action of ${\Bbb Z}^2$
on $\A$.
By the universality subject to the relations
\eqref{eqn:S}, \eqref{eqn:T} and \eqref{eqn:STTS},
one easily sees
that 
the algebra
$\ORE$ for the $C^*$-textile dynamical system
$(\A, \boldsymbol\alpha,  \boldsymbol\beta,
\{ \boldsymbol\alpha \}, \{ \boldsymbol\beta \},
\kappa )
$
is $*$-isomorphic to
the crossed product
$\A \times_{\boldsymbol\alpha, \boldsymbol\beta} {\Bbb Z}^2$.
Take implementing unitaries $U, V$ in
$\A \times_{\boldsymbol\alpha, \boldsymbol\beta} {\Bbb Z}^2$
for the action such that
\begin{equation*}
\boldsymbol\alpha(x) = U^* x U, \qquad
\boldsymbol\beta(x) = V^* x V \qquad\text{for } x \in \A.
\end{equation*}
Since both $U, V$ are unitaries,
we have
$
\boldsymbol\alpha^{-1}(x) = U x U^*, 
\boldsymbol\beta^{-1}(x) = Vx V^* 
$
for
$x \in \A$
which belong to $\A$.
Hence we have 
\begin{equation*}
{\cal B}_{\boldsymbol\alpha} = C^*(U x U^* \mid x \in \A) = \A, \qquad
{\cal B}_{\boldsymbol\beta} = C^*(V x V^* \mid x \in \A) = \A.
\end{equation*}
For 
$ w \in {\cal B}_{\boldsymbol\alpha}, z \in {\cal B}_{\boldsymbol\beta}$
we have
$$
 w = U \boldsymbol\alpha(w) U^*, \qquad
 z = V \boldsymbol\alpha(z) V^*.
 $$
We will write down the Hilbert $C^*$-quad module structure for 
$
{\cal H}_\kappa^{
\boldsymbol\alpha, \boldsymbol\beta}
 = \A$ 
defined in Section 3. 
For $\xi = x, \xi' = x' \in 
{\cal H}_\kappa^{
\boldsymbol\alpha, \boldsymbol\beta}
 = \A,  y \in \A, 
w \in {\cal B}_{\boldsymbol\alpha} =\A, z \in {\cal B}_{\boldsymbol\beta} =\A$,

0. The right $\A$-module and the right $\A$-valued inner product $\langle \cdot \mid \cdot \rangle_\A$:
\begin{equation*}
\xi \varphi_\A(y) = x y \text{  for } y\in \A, 
\qquad 
\langle \xi \mid \xi' \rangle_\A = x^* x'.
\end{equation*}

1. The right action of ${\cal B}_{\boldsymbol\alpha}$ and    
   the right action of ${\cal B}_{\boldsymbol\beta}$:
\begin{equation*}
\xi \varphi_{\boldsymbol\alpha}(w) = x {\boldsymbol\alpha}(w), 
\qquad 
\xi \varphi_{\boldsymbol\alpha}(z) = x {\boldsymbol\beta}(z).
\end{equation*}

2. The right  ${\cal B}_{\boldsymbol\alpha}$-valued inner product and
    the right  ${\cal B}_{\boldsymbol\beta}$-valued inner product:
\begin{equation*}
\langle \xi \mid \xi' \rangle_{\boldsymbol\alpha}
= U x^* x' U^* = {\boldsymbol\alpha}^{-1}(x^*x'),
\qquad 
\langle \xi \mid \xi' \rangle_{\boldsymbol\beta}
= V x^* x' V^* = {\boldsymbol\beta}^{-1}(x^*x').
\end{equation*}

3. The left action of ${\cal B}_{\boldsymbol\alpha}$ and    
   the left action of ${\cal B}_{\boldsymbol\beta}$:
\begin{equation*}
\phi_{\boldsymbol\alpha}(w) \xi = {\boldsymbol\beta}({\boldsymbol\alpha}(w))x, 
\qquad 
\phi_{\boldsymbol\beta}(z) \xi =  {\boldsymbol\alpha}({\boldsymbol\beta}(z)) x.
\end{equation*}

We then have
\begin{equation*}
\widehat{\rho}_{\boldsymbol\alpha} (w) =\boldsymbol\alpha(w), \qquad
\widehat{\eta}_{\boldsymbol\beta}(z) = \boldsymbol\beta(z)
\qquad \text{ for }
w \in {\cal B}_{\boldsymbol\alpha} = \A, 
z \in {\cal B}_{\boldsymbol\beta} = \A.
\end{equation*}
For the $*$-homomorphisms
$
\widehat{\rho}^\eta_{\boldsymbol\alpha}: 
\B_{\boldsymbol\beta} \longrightarrow \B_{\boldsymbol\beta}
$
and
$
\widehat{\eta}^\rho_{\boldsymbol\beta}: 
\B_{\boldsymbol\alpha} \longrightarrow \B_{\boldsymbol\alpha},
$
we have for $z \in {\cal B}_{\boldsymbol\beta}$,
\begin{equation*}
\widehat{\rho}^\eta_{\boldsymbol\alpha} (z)
 = V \boldsymbol\alpha (z_a) V^* 
= V \boldsymbol\alpha (V^* z V) V^*
 = \boldsymbol\beta^{-1}(\boldsymbol\alpha(\boldsymbol\beta (z))) 
   = \boldsymbol\alpha(z)
\end{equation*}
and similarly
$\widehat{\eta}^\rho_{\boldsymbol\beta}(w) =\boldsymbol\beta(w)$
for
$w\in {\cal B}_{\boldsymbol\alpha}$.
Hence the 
$C^*$-algebra 
${\cal O}_{{\cal H}_\kappa^{\boldsymbol\alpha,\boldsymbol\beta}}
$
 defined by the Hilbert $C^*$-quad module 
${\cal H}_\kappa^{
\boldsymbol\alpha, \boldsymbol\beta}
$
is the universal $C^*$-algebra generated by two isometries 
$u, v $ and elements $x \in \A$ subject to the relations
\begin{align*}
u u^* + &  v v^* = 1, \\
u u^* x = x u u^*, & \qquad v v^* x = x v v^*,\\
\boldsymbol\alpha(x) = u^* x u, & \qquad
\boldsymbol\beta(x ) = v^* x v
\end{align*}
for $ x \in \A$.
By the universality of the algebra 
${\cal O}_{{\cal H}_\kappa^{\boldsymbol\alpha,\boldsymbol\beta}}
$
the correspondence
\begin{align*}
u \in {\cal O}_{{\cal H}_\kappa^{\boldsymbol\alpha,\boldsymbol\beta}} 
& \longrightarrow 
\widehat{u}= U \otimes S_1
\in \A\times_{\boldsymbol\alpha,\boldsymbol\beta}{\Bbb Z}^2\otimes {\cal O}_2, \\
v \in {\cal O}_{{\cal H}_\kappa^{\boldsymbol\alpha,\boldsymbol\beta}} 
& \longrightarrow 
\widehat{v}=V \otimes S_2
\in \A\times_{\boldsymbol\alpha,\boldsymbol\beta}{\Bbb Z}^2\otimes {\cal O}_2, \\
x \in \A 
& \longrightarrow 
\widehat{x}=x \otimes 1
\in \A\times_{\boldsymbol\alpha,\boldsymbol\beta}{\Bbb Z}^2\otimes {\cal O}_2
\end{align*}
gives rise to a $*$-homomorphism.
If in particular, 
the action
$(n,m) \in {\Bbb Z}^2 \longrightarrow
\boldsymbol\alpha^n\circ \boldsymbol\beta^m \in \Aut(\A)
$
is outer,
$(\A, \boldsymbol\alpha,  \boldsymbol\beta,
\{ \boldsymbol\alpha \}, \{ \boldsymbol\beta \},
\kappa )
$
satisfies condition (I)
so that 
the above correspondence
gives rise to a $*$-isomorphism,
that is
\begin{equation*}
{\cal O}_{{\cal H}_\kappa^{\boldsymbol\alpha,\boldsymbol\beta}}
\cong
C^*(U \otimes S_1, V \otimes S_2, x \mid x \in \A) \subset 
\A\times_{\boldsymbol\alpha,\boldsymbol\beta}{\Bbb Z}^2\otimes {\cal O}_2.
\end{equation*}

\section{Textile systems of commuting matrices}
We will study $C^*$-algebras associated with the Hilbert $C^*$-quad modules 
defined by textile systems of commuting matrices (cf. \cite{NaMemoir}).
Let $A$ be an $N \times N$ matrix with entries in nonnegative integers.
We may consider a directed graph
$G_A = (V, E_A)$ with vertex set $V = \{v_1,\dots,v_N \}$ 
and edge set $E_A$ consisting of 
 $A(i,j)$ edges from the vertex
$v_i$ to the vertex $v_j$.
Denote by
$\Sigma^A=E_A$.
Let
$\A$ be the $N$-dimensional commutative $C^*$-algebra ${\Bbb C}^N$
with minimal projections 
$E_1,\dots, E_N$
such that
$$
\A = {\Bbb C}E_1 \oplus \cdots \oplus {\Bbb C}E_N.
$$
We set 
for $\alpha \in \Sigma^A, i,j=1,\dots,N$ 
\begin{equation*}
\widehat{A}(i,\alpha,j) =
\begin{cases}
1 & \text{ if } s(\alpha) =v_i, \, r(\alpha) = v_j,\\
0 & \text{ otherwise}
\end{cases}
\end{equation*}
where $s(\alpha), r(\alpha)$ mean the source vertex , the range vertex of an edge $\alpha$
respectively.
We define an endomorphism $\rho_\alpha^A$
on $\A$  for $\alpha \in \Sigma^A$:
\begin{equation*}
\rho^A_\alpha(E_i) = \sum_{j=1}^N \widehat{A}(i,\alpha,j)E_j,
\qquad i=1,\dots,N.
\end{equation*}
Then we have a  $C^*$-symbolic dynamical system
$(\A,\rho^A,\Sigma^A)$.
Let
$B$ be another 
 $N\times N$ matrix with entries in nonnegative integers such that
\begin{equation}
AB = BA. \label{eqn:abba}
\end{equation}
Consider the associated directed graph $G_B = (V, E_B)$
and the $C^*$-symbolic dynamical system
$(\A,\rho^B,\Sigma^B)$
for $B$.  
Let
$S_{\alpha}, \alpha \in E_A$,
$T_a, a \in E_B$
be the generating partial isometries 
of the associated $C^*$-algebras ${\cal O}_{\rho^A}$
and
${\cal O}_{\rho^B}$
respectively.
They 
satisfy the relations:
\begin{align*}
\sum_{\beta \in E_A}S_{\beta}S_{\beta}^*  =1,\qquad
x S_\alpha S_\alpha^* & =  S_\alpha S_\alpha^*  x,\qquad
S_\alpha^* x S_\alpha  = \rho^A_\alpha(x),  \\
\sum_{b \in E_B} T_{b} T_{b}^*  =1,\qquad
x T_a T_a^*  & =  T_a T_a^*  x,\qquad
T_a^* x T_a  = \rho^B_a(x),
\end{align*}
for all $x \in {\cal A}$ and $\alpha \in E_A, a \in E_B$
respectively.
The
 $C^*$-algebras
${\cal O}_{\rho^A}$
and
 ${\cal O}_{\rho^B}$
are isomorphic to the Cuntz-Krieger algebras
${\cal O}_{\widehat{A}}$
and
 ${\cal O}_{\widehat{B}}$
respectively.
Put subalgebras
\begin{align*}
\B_A & = C^*( S_\alpha E_i S_\alpha^* : \alpha \in E_A, i=1,\dots,N) \subset {\cal O}_{\rho^A}, \\
\B_B & = C^*( T_a E_i T_a^* : a \in E_B, i=1,\dots,N) \subset {\cal O}_{\rho^B}.
\end{align*}
Since
$S_\alpha E_i S_\alpha^* \ne 0$ 
if and only if 
$S_\alpha^* S_\alpha = E_i$,
which is equivalent to
$r(\alpha) = v_i$.
Hence
$\B_A$ 
is of $|E_A|$-dimension,
and similarly
$\B_B$ 
is of $|E_B|$-dimension.
By the identities
\begin{equation*}
E_i = \sum_{\alpha \in E_A} S_\alpha \rho^A_\alpha(E_i) S_\alpha^*
    = \sum_{\alpha \in E_A} \sum_{j=1}^N
    \widehat{A}(i,\alpha,j) S_\alpha E_j S_\alpha^*,
\qquad i=1,\dots, N,
\end{equation*}
$E_i $ belongs to $\B_A$ so that
$\A \subset \B_A$
and similarly
$\A \subset \B_B$.
The equality
\eqref{eqn:abba} implies  that the cardinal numbers of the sets of the pairs of directed edges
\begin{align*}
\Sigma^{AB}_{(i,j)}
& =\{(\alpha,b)\in E_A \times E_B 
\mid s(\alpha) = v_i, r(\alpha) = s(b), r(b) = v_j \}, \\ 
\Sigma^{BA}_{(i,j)}
& =\{(a,\beta)\in E_B \times E_A 
\mid s(a) = v_i, r(a) = s(\beta), r(\beta) = v_j \}
\end{align*}
coincide with each other for each
$v_i$ and $v_j$, 
so that one may take a bijection
$
\kappa: \Sigma^{AB}=\cup_{i,j=1}^N \Sigma^{AB}_{(i,j)}
 \longrightarrow 
\Sigma^{BA}=\cup_{i,j=1}^N \Sigma^{BA}_{(i,j)}$
such that
$$
s(\alpha) = s(a),\quad 
r(b) = r(\beta)
\quad
\text{  if  }
\kappa(\alpha,b) = (a,\beta).
$$
We then have
\begin{lem} 
For $(\alpha,b) \in \Sigma^{AB}, (a,\beta) \in \Sigma^{BA}$
with $\kappa(\alpha, b) = (a,\beta)$, we have
\begin{equation*}
\rho^B_b\circ \rho^A_\alpha(E_i) = 
\rho^A_{\beta}\circ \rho^B_a (E_i), \qquad i=1,\dots,N.
\end{equation*}
\end{lem}
\begin{pf}
We have for $i, k=1,\dots,N$
\begin{equation*}
\rho^B_b\circ \rho^A_\alpha(E_i) E_k
= 
\sum_{j=1}^N \widehat{A}(i,\alpha,j)\rho^B_b(E_j) E_k
 = 
 \sum_{j=1}^N \widehat{A}(i,\alpha,j) \widehat{B}(j,b,k) E_k.
\end{equation*}
Hence
$\rho^B_b\circ \rho^A_\alpha(E_i)E_k = E_k
$
if and only if
$v_i = s(\alpha), r(\alpha) = s(b), r(b) = v_k$.
Similarly
for $(a, \beta) \in \Sigma^{BA}$,
 we have 
$\rho^A_{\beta}\circ \rho^B_a (E_i)E_k = E_k
$
if and only if
$v_i = s(a), r(a) = s(\beta), r(\beta) = v_k$.
The condition
 $\kappa(\alpha,b) = (a, \beta)  $
 implies
 $r(\alpha) = s(b), r(a) = s(\beta), s(\alpha) = s(a), r(b) = r(\beta).$
Therefore 
$$
\rho^B_b\circ \rho^A_\alpha(E_i)E_k
=
\rho^A_{\beta}\circ \rho^B_{a}(E_i)E_k
\qquad \text{ for } i,k =1,\dots,N
$$
and hence
$
\rho^B_b\circ \rho^A_\alpha
=
\rho^A_{\beta}\circ \rho^B_a
$
if 
 $\kappa(\alpha,b) = (a, \beta)$.
\end{pf}
We thus have a $C^*$-textile dynamical system
\begin{equation*}
(\A, \rho^A, \rho^B, E_A, E_B, \kappa).
\end{equation*}
We set
$$
E_\kappa
=\{ (\alpha, b, a,\beta) \in E_A \times E_B \times E_B \times E_A | \kappa(\alpha,b) = (a, \beta) \}.
$$
For
$a \in E_B$ and $\alpha \in E_A$,
we will describe the $*$-homomorphisms
$$
\widehat{\eta}^\rho_a:\B_{\rho^A}(={\cal B}_A) \longrightarrow \B_{\rho^A}(={\cal B}_A)
\quad
\text{ and }
\quad
\widehat{\rho}^\eta_\alpha:\B_{\rho^B}(={\cal B}_B) \longrightarrow \B_{\rho^B}(={\cal B}_B)
$$
which will be denoted by $\widehat{\rho}_a^A$
and
 by $\widehat{\rho}^B_\alpha$
respectively.
We set
\begin{equation*}
E_{AB} = \Sigma^{AB}, \qquad
E_{BA} = \Sigma^{BA}. 
\end{equation*}
For $(a, \beta)\in E_{BA}$
there uniquely exists
$(\alpha, b) \in E_{AB}$
such that
$\kappa(\alpha, b) = (a, \beta)$.
We then define a map for
$a \in E_B$
\begin{equation*}
\kappa_{a}: \beta \in \{ \beta \in E_A | (a, \beta) \in E_{BA} \} 
\longrightarrow 
\alpha \in \{ \alpha \in E_A | \kappa(\alpha, b) = (a, \beta) \text{ for some } b \in E_B\}.
\end{equation*}
Similarly, 
we then define a map for $\alpha \in E_A$
\begin{equation*}
\kappa_{\alpha}: b \in \{ b \in E_B | (\alpha,b) \in E_{AB} \} 
\longrightarrow 
a \in \{ a \in E_B | \kappa(\alpha, b) = (a, \beta) \text{ for some } \beta \in E_A\}.
\end{equation*}
We write the projection
$E_i$ also as $E_{v_i}$.
Hence
$w \in \B_A, z \in \B_B
$
are uniquely written as
$w = \sum_{\alpha \in E_A} w(\alpha)S_\alpha E_{r(\alpha)} S_\alpha^*$,
$z = \sum_{a \in E_B} z(a) T_a E_{r(a)}T_a^*$
for $w(\alpha), z(a) \in {\Bbb C}$.
\begin{lem}
Keep the above notations. 
We have
\begin{equation*}
\widehat{\rho}_{a}^A(w) 
= 
\sum_{\beta\in E_A}
w( \kappa_{a}(\beta)) S_{\beta}E_{r(\beta)} S_{\beta}^*,
\qquad
\widehat{\rho}_\alpha^B(z) = 
\sum_{f \in E_B} z(\kappa_\alpha(b)) T_b E_{r(b)}T_b^*.
\end{equation*}
\end{lem}
\begin{pf}
We  have
\begin{equation*}
\widehat{\rho}_{a}^A(w) 
= \sum
\begin{Sb}
\beta, b, \alpha \\
(\alpha, b, a, \beta) \in E_\kappa
\end{Sb}
S_{\beta} \rho_b^B(w(\alpha)E_{r(\alpha)}) S_{\beta}^*  
= \sum
\begin{Sb}
\beta, b, \alpha \\
(\alpha, b, a, \beta) \in E_\kappa
\end{Sb}
w(\alpha)
S_{\beta} T_b^* E_{r(\alpha)} T_b S_{\beta}^*.
\end{equation*} 
Since
$
T_b^* E_{r(\alpha)} T_b = E_{r(b)} = E_{r(\beta)}
$
and
$\alpha = \kappa_{a}(\beta)$,
we see
\begin{equation*}
\widehat{\rho}_{a}^A(w)
=
\sum_{\beta \in E_A}
w( \kappa_{a}(\beta)) S_{\beta}E_{r(\beta)} S_{\beta}^*. 
\end{equation*}
We  similarly have
\begin{equation*}
\widehat{\rho}_{\alpha}^B(z) 
= \sum
\begin{Sb}
b, a, \beta\\
(\alpha, b, a, \beta) \in E_\kappa
\end{Sb}
T_b \rho_{\beta}^A(z(a)E_{r(a)}) T_b^*  
= \sum
\begin{Sb}
b, a, \beta \\
(\alpha, b, a, \beta) \in E_\kappa
\end{Sb}
z(a)
T_b S_{\beta}^* E_{r(\beta)} S_{\beta}^* T_b^*.
\end{equation*} 
Since
$
S_{\beta}^* E_{r(a)} S_{\beta}^* =  E_{r(\beta)} = E_{r(b)} 
$
and
$a = \kappa_\alpha (b)$,
we see
\begin{equation*}
\widehat{\rho}_{\alpha}^B(z)
=
\sum_{b\in E_B}
w( \kappa_\alpha(b)) T_b E_{r(b)} T_b^*. 
\end{equation*}
\end{pf}
For 
$\omega =(\alpha, b,a,\beta) \in E_\kappa$,
put 
$v(\omega) = r(b) (=r(\beta)) \in V$
and
$$
E_\omega = E_{v(\omega)} 
(= \rho^B_b(\rho^A_\alpha(1)) = \rho^A_{\beta}(\rho^B_{a}(1))).
$$ 
Then the Hilbert $C^*$-quad module
${\cal H}_\kappa$
is written as ${\cal H}_\kappa^{A,B}$ and regarded as 
$$
{\cal H}_\kappa^{A,B}
 = \sum_{\omega \in E_\kappa}{}^\oplus {\Bbb C}E_\omega.
$$  
For 
$
\xi = \sum_{\omega \in \EK}{}^\oplus \xi(\omega)E_{v(\omega)},
\xi' = \sum_{\omega \in \EK}{}^\oplus \xi'(\omega)E_{v(\omega)} \in {\cal H}_\kappa^{A,B}
$
and
$
y = \sum_{v \in V} y(v) E_v \in \A$
with
$\xi(\omega), \xi'(\omega)
$ and $y(v) \in {\Bbb C}$, we see

0. The right $\A$-module structure and the right $\A$-valued inner product
are written as follows:
\begin{equation*}
\xi \varphi_\A(y)  = \sum_{\omega \in \EK}{}^\oplus 
\xi(\omega)E_{v(\omega)} y({v(\omega)}),
\qquad
\langle \xi | \xi' \rangle_\A
 = \sum_{\omega \in \EK}\overline{\xi(\omega)} \xi'(\omega) E_{v(\omega)}.
\end{equation*}

Furthermore for
$w = \sum_{\alpha \in E_A} w(\alpha)S_\alpha E_{r(\alpha)} S_\alpha^*$,
$z = \sum_{a \in E_B} z(a) T_a E_{r(a)}T_a^*$
with
$w(\alpha), z(a) \in {\Bbb C}$,

1. The 
right $\B_A$-action $\varphi_A$ 
and the
right $\B_B$-action $\varphi_B$
are written as follows:
\begin{align*}
\xi \varphi_A(w) 
& = \sum_{\omega \in \EK}{}^\oplus 
    \xi(\omega)E_{v(\omega)}\varphi_A(w) 
 = \sum_{\omega \in \EK}{}^\oplus 
    \xi(\omega) w(b(\omega)) E_{v(\omega)}, \\
\xi \varphi_B(z) 
& = \sum_{\omega \in \EK}{}^\oplus 
    \xi(\omega)E_{v(\omega)}\varphi_B(z) 
 = \sum_{\omega \in \EK}{}^\oplus 
    \xi(\omega) z(r(\omega)) E_{v(\omega)}.
\end{align*}

2. The 
left $\B_A$-action $\phi_A$ 
and the
left $\B_B$-action $\phi_B$
are written as follows:
\begin{align*}
 \phi_A(w) \xi 
& = \sum_{\omega \in \EK}{}^\oplus 
    \xi(\omega)E_{v(\omega)}\phi_A(w) \\
& =\sum_{\omega \in \EK}{}^\oplus 
   \xi(\omega)E_{v(\omega)} 
   \rho^B_{r(\omega)}(w(t(\omega ))E_{r(t(\omega ))}) \\
& =\sum_{\omega \in \EK}{}^\oplus 
   \xi(\omega) w(t(\omega))E_{v(\omega)} \rho^B_{r(\omega)}(E_{r(t(\omega))}).
\end{align*}
As 
$\rho^B_{r(\omega)}(E_{r(t(\omega))}) =  E_{v(\omega)}$,
we have
\begin{equation*}
\phi_A(w) \xi  = \sum_{\omega \in \EK}{}^\oplus 
   \xi(\omega) w(t(\omega))E_{v(\omega)}.
\end{equation*}
We also have
\begin{align*}
 \phi_B(z) \xi 
& = \sum_{\omega \in \EK}{}^\oplus 
    \xi(\omega)E_{v(\omega)}\phi_B(z) \\
& =\sum_{\omega \in \EK}{}^\oplus 
   \xi(\omega)E_{v(\omega)} 
   \rho^A_{b(\omega)}(z(l(\omega))E_{r(l(\omega))}) \\
& =\sum_{\omega \in \EK}{}^\oplus 
   \xi(\omega) z(l(\omega))E_{v(\omega)} \rho^A_{b(\omega)}(E_{r(l(\omega))}).
\end{align*}
As 
$\rho^A_{b(\omega)}(E_{r(l(\omega))}) =  E_{v(\omega)}$,
we have
\begin{equation*}
\phi_B(z) \xi  = \sum_{\omega \in \EK}{}^\oplus 
   \xi(\omega) z(l(\omega))E_{v(\omega)}.
\end{equation*}

3. The 
right $\B_A$-valued inner product 
$\langle \cdot | \cdot \rangle_A$
and the
right $\B_B$-valued inner product 
$\langle \cdot | \cdot \rangle_B$
are written  as follows:
\begin{align*}
\langle \xi | \xi' \rangle_A
& =
\sum_{\omega \in \EK} S_{b(\omega)} \overline{\xi(\omega)}
        E_{v(\omega)} \xi'(\omega) S_{b(\omega)}^*, \\
\langle \xi | \xi' \rangle_B
& =
\sum_{\omega \in \EK} T_{r(\omega)} \overline{\xi(\omega)}
        E_{v(\omega)} \xi'(\omega) T_{r(\omega)}^*.   
\end{align*}
As
$ E_{v(\omega)}  = S_{b(\omega)}^* S_{b(\omega)} 
= T_{r(\omega)}^* T_{r(\omega)},
$
we have
\begin{equation*}
\langle \xi | \xi' \rangle_A
 =
\sum_{\omega \in \EK} \overline{\xi(\omega)}
         \xi'(\omega) S_{b(\omega)} S_{b(\omega)}^*, 
\qquad
\langle \xi | \xi' \rangle_B
 =
\sum_{\omega \in \EK}  \overline{\xi(\omega)}
         \xi'(\omega) T_{r(\omega)} T_{r(\omega)}^*.   
\end{equation*}

4. The positive maps
$\lambda_A:\B_A \longrightarrow \A$
and
$\lambda_B:\B_B \longrightarrow \A$
are written as follows:
\begin{equation*}
\lambda_A(w)  = \sum_{\alpha \in E_A} w(\alpha) E_{r(\alpha)}, 
\qquad 
\lambda_B(z)  = \sum_{a \in E_B} z(a) E_{r(a)}.
\end{equation*}
Hence we have
\begin{equation*}
\lambda_A(\langle \xi | \xi' \rangle_A)
 =
\langle \xi | \xi' \rangle_\A,   
\qquad
\lambda_B(\langle \xi | \xi' \rangle_B)
 =
\langle \xi | \xi' \rangle_\A.   
\end{equation*}
Put
$p_\alpha = S_\alpha E_{r(\alpha)} S_\alpha^*$ for $\alpha \in E_A$
and
$q_a = T_a E_{r(a)} T_a^*$ for $a \in E_B$.
Hence
\begin{equation*}
\B_A = \sum_{\alpha \in E_A} {\Bbb C}p_\alpha, \qquad
\B_B = \sum_{a \in E_B} {\Bbb C}q_a.
\end{equation*}
By Lemma 7.2,
we have
\begin{align}
\widehat{\rho}^A_a(w)
& = \sum_{\beta \in E_A}w(\kappa_a(\beta)) p_\beta \quad
\text{ for }
w = \sum_{\alpha \in E_A} w(\alpha)p_\alpha \in \B_A, \label{eqn:rhoA}\\
\widehat{\rho}^B_\alpha(z)
& = \sum_{b \in E_B}z(\kappa_\alpha(b)) q_b \qquad
\text{ for } 
z = \sum_{a \in E_B} z(a) q_a \in \B_B \label{eqn:rhoB}. 
\end{align}
We define
$
\kappa_B :
 E_A \times E_B \times E_A \longrightarrow \{0,1\}
$
and
$
\kappa_A :
 E_B \times E_A \times E_B \longrightarrow \{0,1\}
$
by 
\begin{equation*}
 \kappa_B(\alpha,a,\beta) 
=
{
\begin{cases}
1 & \text{ if } \kappa_a(\beta) =\alpha, \\
0 & \text{ otherwise }
\end{cases}
}
\quad
\text{ and }
\quad
 \kappa_A(a, \alpha, b) 
=
{
\begin{cases}
1 & \text{ if } \kappa_\alpha(b) =a, \\
0 & \text{ otherwise. }
\end{cases}
} 
\end{equation*}
We identify the
 $C^*$-algebra
${\cal O}_{{\cal H}^{A,B}_\kappa}$
with the universal $C^*$-algebra subject to the relation
$({\cal H}^{A,B}_\kappa)$,
and denote  the generating partial isometries
by
${\tt u}_\alpha, \alpha \in E_A$
and
${\tt v}_a, a \in E_B$.
Since 
$\widehat{\rho}^A_a(w) = {\tt v}_a^* w {\tt v}_a$
and
$\widehat{\rho}^B_\alpha(z) = {\tt u}_\alpha^* z {\tt u}_\alpha$
with \eqref{eqn:rhoA} and  \eqref{eqn:rhoB},
we have
\begin{lem}
\begin{equation*}
{\tt v}_a^* p_\alpha {\tt v}_a = 
\sum_{\beta \in E_A} \kappa_B(\alpha,a,\beta) p_\beta,
\qquad
{\tt u}_\alpha^* q_a {\tt u}_\alpha = 
\sum_{b \in E_B}
\kappa_A(a,\alpha, b)  q_b.
\end{equation*}
\end{lem}
Define 
$|E_A| \times |E_A|$-matrix $A^E=[A^E(\alpha,\beta)]_{\alpha,\beta\in E_A}$ 
and
$|E_B| \times |E_B|$-matrix $B^E=[B^E(a,b)]_{a,b\in E_B}$ 
by 
\begin{equation*}
A^E(\alpha,\beta)
=\begin{cases}
1 & \text{ if } r(\alpha) = s(\beta),\\
0 & \text{ if } r(\alpha) \ne s(\beta),
\end{cases}
\qquad
B^E(a,b)
=\begin{cases}
1 & \text{ if } r(a) = s(b),\\
0 & \text{ if } r(a) \ne s(b)
\end{cases}
\end{equation*}
respectively.
We then have
for $\delta \in E_A$,
\begin{equation*}
E_{r(\delta)}
= \sum_{\beta \in E_A}S_\beta\rho^A_\beta(E_{r(\delta)})S_\beta^* 
 = \sum_{\beta\in E_A}S_\beta 
\widehat{A}(r(\delta), \beta, r(\beta)) E_{r(\beta)} S_\beta^* 
 = \sum_{\beta\in E_A} 
A^E(\delta, \beta) p_\beta  
\end{equation*}
and similarly
for
$d \in E_B$
\begin{equation*}
E_{r(d)}
= \sum_{b\in E_B} 
B^E(d, b) q_b.
\end{equation*}
Since
we know that
\begin{equation*}
\widehat{\rho}^B_\alpha(p_\delta)
= 
{\begin{cases}
E_{r(\delta)} & \text{ if } \delta = \alpha,\\
0  & \text{ otherwise ,} 
\end{cases}
}
\qquad
\widehat{\rho}^A_a(q_d)
= 
{\begin{cases}
E_{r(d)} & \text{ if } d = a, \\
0  & \text{ otherwise } 
\end{cases}
}
\end{equation*}
and
$
\sum_{\delta \in E_A} p_\delta 
= 
\sum_{d \in E_B} q_d =1,
$
we have
\begin{equation*}
{\tt u}_\alpha^* {\tt u}_\alpha
= \sum_{\beta\in E_A} 
A^E(\alpha, \beta) p_\beta,
\qquad
{\tt v}_a^* {\tt v}_a 
= \sum_{b \in E_B}
B^E(a,b) q_b.
\end{equation*}
Therefore we have
\begin{prop} \label{prop:Prop1}
The $C^*$-algebra ${\cal O}_{{\cal H}_\kappa^{A,B}}$
is $*$-isomorphic to the universal $C^*$-algebra
generated by 
two families of projections
$\{ p_\alpha \}_{\alpha \in E_A}$,
$\{ q_a \}_{a \in E_B}$
and
two families of partial isometries
$\{ {\tt u}_\alpha \}_{\alpha \in E_A}$,
$\{ {\tt v}_a \}_{a \in E_B}$
subject to the  relations: 
\begin{align}
\sum_{\beta \in E_A} p_\beta 
=
\sum_{b \in E_B} q_b 
 =
& \sum_{\beta \in E_A} {\tt u}_\beta {\tt u}_\beta^* 
+
 \sum_{b \in E_B} {\tt v}_b {\tt v}_b^* =1, \label{eqn:P1} \\ 
{\tt u}_\alpha {\tt u}_\alpha^* p_\alpha  =  {\tt u}_\alpha {\tt u}_\alpha^*, 
& \qquad
{\tt v}_a {\tt v}_a^* q_a  =  {\tt v}_a {\tt v}_a^*, \label{eqn:P2} \\
{\tt u}_\alpha {\tt u}_\alpha^* q_a  = q_a {\tt u}_\alpha {\tt u}_\alpha^*, 
& \qquad
{\tt v}_a {\tt v}_a^* p_\alpha = p_\alpha {\tt v}_a {\tt v}_a^*, \label{eqn:P3}\\
{\tt u}_\alpha^* {\tt u}_\alpha
= \sum_{\beta \in E_A} A^E(\alpha,\beta)p_\beta,
&  \qquad
{\tt v}_a^*  {\tt v}_a
=
\sum_{b \in E_B} B^E(a,b)q_b, \label{eqn:P4}\\ 
{\tt u}_\alpha^* q_a {\tt u}_\alpha = 
\sum_{b \in E_B}
\kappa_A(a,\alpha, b) q_b,
& \qquad
{\tt v}_a^* p_\alpha {\tt v}_a = 
\sum_{\beta \in E_A} \kappa_B(\alpha,a,\beta) p_\beta \label{eqn:P5}
\end{align}
for $\alpha \in E_A$ and $a \in E_B$ where
\begin{align*}
\kappa_A(a,\alpha, b)
& =
{
\begin{cases}
1 & \text{ if } \kappa(\alpha,b)=(a,\beta) \text{ for some } \beta \in E_A,\\
0 & \text{ otherwise,} 
\end{cases}
} \\
\kappa_B(\alpha,a,\beta)
& =
{\begin{cases}
1 & \text{ if } \kappa(\alpha,b) = (a,\beta) \text{ for some } b \in E_B,\\
0 & \text{ otherwise. } 
\end{cases}
}
\end{align*}
\end{prop}
\begin{pf}
The relations 
\eqref{eqn:P1}, \eqref{eqn:P2} and \eqref{eqn:P4}
 imply 
\begin{equation*}
{\tt u}_\alpha^* p_\delta {\tt u}_\alpha
= 
\begin{cases}
\sum_{\beta \in E_A} A^E(\alpha,\beta)p_\beta
& \text{ if } \delta = \alpha,\\
0 
& \text{ otherwise}
\end{cases}
\end{equation*}
and
\begin{equation*}
{\tt v}_a^* q_d {\tt v}_a
=
\begin{cases}
\sum_{b \in E_B}B^E(a,b)q_b & \text{ if } d=a,    \\
0                                  & \text{otherwise.}
\end{cases}   
\end{equation*}
The above two relations are equivalent to 
\eqref{eqn:M3} and \eqref{eqn:M4}.
Since the two $C^*$-algebras
$\B_A$ and $\B_B$
are generated by 
the projections
$\{ p_\alpha \}_{\alpha \in E_A}$
and
$\{ q_a \}_{a \in E_B}$
respectively,
we see that the relations 
\eqref{eqn:P1}, \eqref{eqn:P2}, \eqref{eqn:P3}, \eqref{eqn:P4} and \eqref{eqn:P5}
are equivalent to
the relations $({\cal H}_\kappa^{A,B})$.
\end{pf}
We will further study 
the above operator relations.
\begin{lem}
$p_\alpha$ commutes with $q_a$
for all $\alpha \in E_A$, $a \in E_B$.
\end{lem}
\begin{pf}
For $\alpha \in E_A$,  $a \in E_B$,
by \eqref{eqn:P1}, we have
\begin{equation*}
p_\alpha =
\sum_{\beta \in E_A}{\tt u}_\beta {\tt u}_\beta^* p_\alpha
+
\sum_{b \in E_B} {\tt v}_b {\tt v}_b^* p_\alpha
\end{equation*}
By \eqref{eqn:P2}, we have
\begin{equation*}
{\tt u}_\beta {\tt u}_\beta^* p_\alpha =
\begin{cases}
{\tt u}_\alpha {\tt u}_\alpha^*  & \text{ if } \alpha = \beta,\\
0 & \text{ if } \alpha \ne \beta,
\end{cases}
\qquad
q_a {\tt v}_b {\tt v}_b^*  =
\begin{cases}
{\tt v}_a {\tt v}_a^*  & \text{ if } a = b,\\
0 & \text{ if } a \ne b
\end{cases}
\end{equation*}
so that
$
p_\alpha =
{\tt u}_\alpha {\tt u}_\alpha^* 
+
\sum_{b \in E_B} {\tt v}_b {\tt v}_b^* p_\alpha
$
and hence
$
q_a p_\alpha =
q_a {\tt u}_\alpha {\tt u}_\alpha^* 
+
{\tt v}_a {\tt v}_a^* p_\alpha.
$
Since
$q_a$ commutes with
${\tt u}_\alpha {\tt u}_\alpha^*$
and
$p_\alpha$ commutes with
${\tt v}_a {\tt v}_a^*$,
we have
\begin{equation}
q_a p_\alpha =
 {\tt u}_\alpha {\tt u}_\alpha^* q_a 
+
p_\alpha {\tt v}_a {\tt v}_a^*, \label{eqn:pq}
\end{equation}
which is symmetrically 
equal to
$p_\alpha q_a$.
\end{pf}
Put for
$\alpha,\beta \in E_A$
and
$a,b \in E_B$
\begin{equation*}
\kappa_{AB}(\alpha,b) 
 = 
{
\begin{cases}
1 & \text{ if } r(\alpha) = s(b), \\
0 & \text{ otherwise, }
\end{cases}
} 
\qquad
\kappa_{BA}(a,\beta) 
= 
{
\begin{cases}
1 & \text{ if } r(a) = s(\beta), \\
0 & \text{ otherwise. }
\end{cases}
}
\end{equation*}
\begin{lem} For $\alpha \in E_A, a \in E_B$, we have
\begin{equation*}
{\tt u}_\alpha^* {\tt u}_\alpha 
=
\sum_{b \in E_B}\kappa_{AB}(\alpha,b) q_b, \qquad
{\tt v}_a^* {\tt v}_a 
=
\sum_{\beta \in E_A}\kappa_{BA}(a,\beta) p_\beta.
\end{equation*}
\end{lem}
\begin{pf}
By \eqref{eqn:P1} and \eqref{eqn:P5}
and the equality
$
\sum_{a \in E_B}
\kappa_A(a,\alpha, b) 
=
\kappa_{AB}(\alpha,b)$,
 we have
\begin{equation*}
{\tt u}_\alpha^* {\tt u}_\alpha
=
\sum_{a \in E_B}
{\tt u}_\alpha^* q_a {\tt u}_\alpha 
=
\sum_{a \in E_B}
\sum_{b \in E_B}
\kappa_A(a,\alpha, b) q_b
=
\sum_{b \in E_B}\kappa_{AB}(\alpha,b) q_b.
\end{equation*}
The equality for ${\tt v}_a^* {\tt v}_a$ is similarly shown.
\end{pf}
\begin{lem}
For $\alpha \in E_A$,  $a \in E_B$,
if
$r(\alpha) = r(a)$,
then 
${\tt u}_\alpha^* {\tt u}_\alpha = {\tt v}_a^* {\tt v}_a$.
\end{lem}
\begin{pf}
The condition
$r(\alpha) = r(a)$
implies
$
\sum_{b \in E_B}\kappa_{AB}(\alpha,b) q_b
=
\sum_{b \in E_B}B^E(a,b) q_b.
$
By the equality for ${\tt u}_\alpha^*{\tt u}_\alpha$ 
in the preceding lemma 
and 
the equality for 
${\tt v}_a^* {\tt v}_a$
in
\eqref{eqn:P4},
we have
$
{\tt u}_\alpha^* {\tt u}_\alpha
=
{\tt v}_a^* {\tt v}_a.
$
\end{pf}
Put
\begin{equation*}
\Omega_\kappa
=\{  (\alpha,a) \in   E_A\times E_B |  s(\alpha) = s(a),   
\kappa(\alpha,b) =(a,\beta) \text{ for some } \beta \in E_A, b\in E_B \} 
\end{equation*}
and
$e_{\alpha,a} = p_\alpha q_a$
for $(\alpha, a ) \in \Omega_\kappa$. 
\begin{lem}
\begin{enumerate}
\renewcommand{\labelenumi}{(\roman{enumi})}
\item
$\sum_{(\alpha,a) \in \Omega_\kappa} e_{(\alpha,a)} =1.$
\item
The $C^*$-subalgebra $\BK$ of $\OHK$
generated by the subalgebras $\BR$ and $\BE$ 
is $*$-isomorphic to 
the direct sum
$\sum_{(\alpha,a) \in \Omega_\kappa}{}^\oplus {\Bbb C}e_{(\alpha,a)}$.
It is the $C^*$-algebra of all complex valued continuous functions 
on $\Omega_\kappa$. 
\end{enumerate}
\end{lem}
\begin{pf}
By the equality \eqref{eqn:pq},
one sees that
$p_\alpha q_a \ne 0$
if and only if
$
{\tt u}_\alpha^* q_a {\tt u}_\alpha \ne 0
$
or
$
{\tt v}_a^* p_\alpha {\tt v}_a
\ne 0.
$
The latter condition is equivalent to
the condition that
there exists $b \in E_B$ such that $\kappa_A(a,\alpha,b) \ne 0$
or
there exists $\beta \in E_A$ such that $\kappa_B(\alpha,a,\beta) \ne 0$,
which is also equivalent to the condition that
there exist $b \in E_B$ and $\beta \in E_A$ 
such that $\kappa(\alpha,b) =(a,\beta) $.
Hence we have
$p_\alpha q_a \ne 0$
if and only if
$(\alpha,a) \in \Omega_\kappa$.
Therefore we have
$$
\sum_{(\alpha,a) \in \Omega_\kappa} e_{(\alpha,a)} =
(\sum_{\alpha \in E_A} p_\alpha )\cdot (\sum_{a \in E_B} q_a) =1.
$$
As 
$p_\alpha q_a \cdot p_{\alpha'} q_{a'} =0$
if 
$\alpha \ne \alpha'$ or
$a \ne a'$,
the $C^*$-algebra 
$\BK$
is $*$-isomorphic to
$ \sum_{(\alpha,a) \in \Omega_\kappa}{}^\oplus {\Bbb C}e_{(\alpha,a)}$.
\end{pf}
We define 
two $|\Omega_\kappa| \times |\Omega_\kappa|$-matrcies
$A_\kappa$ and $B_\kappa$  with entries in $\{0,1\}$
by 
\begin{align*}
A_\kappa((\alpha,a),(\delta,b))
& = 
{\begin{cases}
1 & \text{ if there exists } \beta\in E_A  \text{ such that } \kappa(\alpha,b) = (a,\beta),\\  
0 & \text{ otherwise}
\end{cases}
}
\\
\intertext{ for
$(\alpha,a),(\delta,b) \in \Omega_\kappa$, and} 
B_\kappa((\alpha,a), (\beta,d))           
& = 
{\begin{cases}
1 & \text{ if there exists } b \in E_B  \text{ such that } \kappa(\alpha,b) = (a,\beta),\\  
0 & \text{ otherwise}
\end{cases}
}
 \end{align*}
for 
$(\alpha,a), (\beta,d) \in \Omega_\kappa$
respectively.
They represent the concatenations of edges as in the following figures respectively:
\begin{equation*}
{\begin{CD}
\circ @>\alpha>> \circ @>\delta>> \\
@V{a}VV  @V{b}VV @. \\ 
\circ @>\beta>> \circ @.
\end{CD}
}
\qquad 
\quad
\text{ and }
 \quad
 \qquad
{
\begin{CD}
\circ @>\alpha>> \circ \\
@V{a}VV   @V{b}VV \\
\circ @>\beta >> \circ  \\
@V{d}VV   @. 
\end{CD}
}
\end{equation*}
\begin{prop} \label{prop:Prop2}
The $C^*$-algebra ${\cal O}_{{\cal H}_\kappa^{A,B}}$
is $*$-isomorphic to the universal $C^*$-algebra
generated by 
a family 
$\{ e_{(\alpha,a)} \}_{(\alpha,a) \in\Omega_\kappa }$
of projections and
two families of partial isometries
$\{ {\tt u}_\alpha \}_{\alpha \in E_A}$,
$\{ {\tt v}_a \}_{a \in E_B}$
subject to the  relations: 
\begin{align}
\sum_{(\alpha,a) \in\Omega_\kappa } e_{(\alpha,a)} 
 =
& \sum_{\beta \in E_A} {\tt u}_\beta {\tt u}_\beta^* 
+
 \sum_{b \in E_B} {\tt v}_b {\tt v}_b^* =1, \label{eqn:P1'}\\ 
{\tt u}_\alpha {\tt u}_\alpha^*  = 
& \sum_{a \in E_B} {\tt u}_\alpha {\tt u}_\alpha^* e_{(\alpha,a)}
=\sum_{a \in E_B}  e_{(\alpha,a)} {\tt u}_\alpha {\tt u}_\alpha^*, \label{eqn:P2'}\\
{\tt v}_a {\tt v}_a^*  = 
&  \sum_{ \alpha \in E_A} {\tt v}_a {\tt v}_a^* e_{(\alpha,a)} 
=  \sum_{ \alpha \in E_A} e_{(\alpha,a)} {\tt v}_a {\tt v}_a^* , \label{eqn:P3'}\\
{\tt u}_\alpha^* e_{(\alpha,a)} {\tt u}_\alpha 
=
&
\sum_{(\delta,b) \in\Omega_\kappa }A_\kappa((\alpha,a),(\delta,b)) e_{(\delta,b)}, \label{eqn:P4'}\\
{\tt v}_a^*e_{(\alpha,a)}{\tt v}_a
=
&
\sum_{(\beta,d) \in\Omega_\kappa }B_\kappa((\alpha,a),(\beta,d)) e_{(\beta,d)} \label{eqn:P5'}
\end{align}
for $\alpha \in E_A$ and $a \in E_B$.
\end{prop}
\begin{pf}
Let ${\tt u}_\alpha, \alpha\in E_A$ and ${\tt v}_a, a \in E_B$
be the partial isometries as in Proposition \ref{prop:Prop1}.
The equalities \eqref{eqn:P2'} and  \eqref{eqn:P3'}
follow from the equalities \eqref{eqn:P1} and \eqref{eqn:P2}, \eqref{eqn:P3}.
As
${\tt u}_\alpha^* p_\alpha {\tt u}_\alpha = {\tt u}_\alpha {\tt u}_\alpha^*$
by \eqref{eqn:P2},
we have by \eqref{eqn:P5}
\begin{align*}
{\tt u}_\alpha^* e_{(\alpha,a)} {\tt u}_\alpha 
 = {\tt u}_\alpha^* q_a {\tt u}_\alpha 
& = \sum_{b \in E_B}\kappa_A(a,\alpha,b) q_b \\
& = \sum_{b \in E_B}\kappa_A(a,\alpha,b) \sum_{\delta \in E_A} p_\delta q_b \\
& =
\sum_{(\delta,b) \in\Omega_\kappa }A_\kappa((\alpha,a),(\delta,b)) e_{(\delta,b)}.
\end{align*}
The equality \eqref{eqn:P5'} is similarly shown.
Hence the equalities \eqref{eqn:P1'},\dots, \eqref{eqn:P5'}
follow from the equalities \eqref{eqn:P1},\dots, \eqref{eqn:P5}.
Conversely, from the projections
$e_{(\alpha,a)}, (\alpha,a) \in \Omega_\kappa$
by putting 
$$
p_\alpha =
\sum_{a \in E_B} e_{(\alpha,a)}, \qquad
q_a =
\sum_{\alpha \in E_A} e_{(\alpha,a)}
$$
the equalities \eqref{eqn:P1},\dots, \eqref{eqn:P5}
follow from the equalities \eqref{eqn:P1'},\dots, \eqref{eqn:P5'}.
\end{pf}
We then see the following theorem:
\begin{thm}\label{thm:CK} 
The $C^*$-algebra ${\cal O}_{{\cal H}_\kappa^{A,B}}$
associated with the 
Hilbert $C^*$-quad module ${\cal H}_\kappa^{A,B}$
defined by commuting matrices
$A, B$
and a specification $\kappa$
is 
generated by partial isometries
$
S_{(\alpha,a)}, T_{(\alpha,a)}
$ 
for
$(\alpha,a)\in \Omega_\kappa$
satisfying the  relations: 
\begin{align*}
\sum_{(\delta,b) \in\Omega_\kappa} 
& S_{(\delta,b)}S_{(\delta,b)}^*
+
\sum_{(\beta,d) \in\Omega_\kappa} 
T_{(\beta,d)}T_{(\beta,d)}^* = 1, \\
S_{(\alpha,a)}^*S_{(\alpha,a)}
=
&
\sum_{(\delta,b) \in\Omega_\kappa}
A_\kappa((\alpha,a),(\delta,b)) (
S_{(\delta,b)}S_{(\delta,b)}^* + T_{(\delta,b)} T_{(\delta,b)}^*),\\
T_{(\alpha,a)}^* T_{(\alpha,a)}
=
&
\sum_{(\beta,d) \in\Omega_\kappa} B_\kappa((\alpha,a),(\beta,d)) 
(
S_{(\beta,d)}S_{(\beta, d)}^* + T_{(\beta,d)} T_{(\beta,d)}^*)
\end{align*}
for $(\alpha,a) \in  \Omega_\kappa$.
\end{thm}
\begin{pf}
The algebra
${\cal O}_{{\cal H}_\kappa^{A,B}}$
is generated by 
${\tt u}_\alpha, \alpha \in E_A, {\tt v}_a, a \in E_A$
and
$e_{(\alpha,a)}, (\alpha,a) \in \Omega_\kappa$
as in the preceding proposition.
For
$(\alpha,a) \in \Omega_\kappa$,
put
\begin{equation}
S_{(\alpha,a)} = e_{(\alpha,a)} {\tt u}_\alpha, \qquad
T_{(\alpha,a)} = e_{(\alpha,a)} {\tt v}_a.
 \end{equation}
Denote by
$C^*(S_{(\alpha,a)}, T_{(\alpha,a)} : (\alpha,a) \in \Omega_\kappa)$
the $C^*$-subalgebra 
of
${\cal O}_{{\cal H}_\kappa^{A,B}}$
generated by elements
$S_{(\alpha,a)},  T_{(\alpha,a)}, (\alpha,a) \in \Omega_\kappa$.
We have
\begin{equation*}
S_{(\alpha,a)}^*S_{(\alpha,a)}
= {\tt u}_\alpha^* e_{(\alpha,a)} {\tt u}_\alpha
=
\sum_{(\delta,b) \in\Omega_\kappa}
A_\kappa((\alpha,a),(\delta,b)) 
e_{(\delta,b)}.
\end{equation*}
As
$e_{(\alpha,a)} {\tt u}_\beta  = 0 
$
for
$\beta \ne \alpha$,
and
$
e_{(\alpha,a)} {\tt v}_b  = 0
$
for
$ b \ne a,$
we have
\begin{align*}
e_{(\alpha,a)} 
& =
\sum_{\beta \in E_A} e_{(\alpha,a)} {\tt u}_\beta {\tt u}_\beta^* e_{(\alpha,a)} 
+
\sum_{b \in E_B} e_{(\alpha,a)} {\tt v}_b {\tt v}_b^* e_{(\alpha,a)} \\
& =
 e_{(\alpha,a)}  {\tt u}_\alpha {\tt u}_\alpha^* e_{(\alpha,a)} 
+
e_{(\alpha,a)} {\tt v}_a {\tt v}_a^* e_{(\alpha,a)} \\
& =
 S_{(\alpha,a)} S_{(\alpha,a)} ^*
 +
T_{(\alpha,a)} T_{(\alpha,a)} ^*
\end{align*}
so that
$e_{(\alpha,a)} $ belongs to the algebra
$C^*(S_{(\alpha,a)}, T_{(\alpha,a)} : (\alpha,a) \in \Omega_\kappa)$
and the equality
\begin{align*}
S_{(\alpha,a)}^*S_{(\alpha,a)}
=
&
\sum_{(\delta,b) \in\Omega_\kappa}
A_\kappa((\alpha,a),(\delta,b)) (
S_{(\delta,b)}S_{(\delta,b)}^* + T_{(\delta,b)} T_{(\delta,b)}^*)\\
\intertext{holds. Similarly we have}
T_{(\alpha,a)}^* T_{(\alpha,a)}
=
&
\sum_{(\beta,d) \in\Omega_\kappa}B_\kappa((\alpha,a),(\beta,d)) 
(
S_{(\beta,d)}S_{(\beta, d)}^* + T_{(\beta,d)} T_{(\beta,d)}^*).
\end{align*}
As 
$\sum_{(\alpha,a) \in \Omega_\kappa}e_{(\alpha,a)} = 1$
and $e_{(\alpha,a)} {\tt u}_\beta =0$ for $\beta \ne \alpha$,
we have
\begin{equation*}
{\tt u}_\alpha 
=  \sum\begin{Sb}
a \in \Sigma^\eta\\
(\alpha,a) \in \Omega_\kappa
\end{Sb}
e_{(\alpha,a)} {\tt u}_\alpha
=  \sum_{
(\alpha,a) \in \Omega_\kappa}
S_{(\alpha,a)} 
\end{equation*}
so that
${\tt u}_\alpha $ and similarly 
${\tt v}_a$ belong to the algebra
$C^*(S_{(\alpha,a)}, T_{(\alpha,a)} : (\alpha,a) \in \Omega_\kappa)$.
Therefore the $C^*$-algebra generated by 
$e_{(\alpha,a)}, {\tt u}_\alpha, {\tt v}_a $
coincides
with the subalgebra
$C^*(S_{(\alpha,a)}, T_{(\alpha,a)} : (\alpha,a) \in \Omega_\kappa)$.
\end{pf}
Put $n = |\Omega_\kappa|$.
Define a 
$2 n \times 2 n$-matrix
$H_\kappa$
with entries in $\{ 0,1 \}$
by the block matrix
\begin{equation*}
H_\kappa
=
\begin{bmatrix}
A_\kappa & A_\kappa \\
B_\kappa & B_\kappa.
\end{bmatrix}.
\end{equation*}
Denote by $I_n$ and $I_{2n}$ the identity matrices of size $n$
and of size $2n$  respectively.
\begin{lem}
\hspace{7cm}
\begin{enumerate}
\renewcommand{\labelenumi}{(\roman{enumi})}
\item
$
{\Bbb Z}^{2n} / (H_\kappa - I_{2n}){\Bbb Z}^{2n}
\cong 
{\Bbb Z}^{n} / (A_\kappa + B_\kappa - I_n){\Bbb Z}^n.
$
\item
$\Ker(H_\kappa - I_{2n})$ in ${\Bbb Z}^{2n}
\cong
\Ker(A_\kappa + B_\kappa - I_n)
$
in 
${\Bbb Z}^n.$
\end{enumerate}
\end{lem}
\begin{pf}
(i)
Put a $2n \times 2n$ block matrix
$
\widehat{H}_\kappa =
\begin{bmatrix}
A_\kappa  & I_n \\
B_\kappa  & 0 
\end{bmatrix}.
$
Then we easily see
\begin{align*}
{\Bbb Z}^{2n} / (H_\kappa - I_{2n}){\Bbb Z}^{2n}
& \cong 
{\Bbb Z}^{2n} / (\widehat{H}_\kappa - I_{2n}){\Bbb Z}^{2n}.
\end{align*}
Define a map
\begin{equation*}
\Psi:
(x_1,\dots,x_n,y_1,\dots,y_n) 
\in {\Bbb Z}^{2n} \longrightarrow
(x_1+y_1,\dots,x_n +y_n) \in {\Bbb Z}^n
\end{equation*}
which is a surjective homomorphism of abelian groups 
from 
${\Bbb Z}^{2n}$ to
${\Bbb Z}^n$.
Since we know
$\Psi((\widehat{H}_\kappa - I_{2n}){\Bbb Z}^{2n} )
=(A_\kappa + B_\kappa - I_n){\Bbb Z}^n,
$
the homomorphism
$\Psi:{\Bbb Z}^{2n} \longrightarrow {\Bbb Z}^n$
induces an isomorphism
from
$
{\Bbb Z}^{2n} / (\widehat{H}_\kappa - I_{2n}){\Bbb Z}^{2n}
$
to
$ 
{\Bbb Z}^{n} / (A_\kappa + B_\kappa - I_n){\Bbb Z}^n.
$
Therefore 
$
{\Bbb Z}^{2n} / (H_\kappa - I_{2n}){\Bbb Z}^{2n}
$
is isomorphic to
$ 
{\Bbb Z}^{n} / (A_\kappa + B_\kappa - I_n){\Bbb Z}^n.
$

(ii)
The groups 
$
\Ker(H_\kappa - I_{2n})\text{ in } {\Bbb Z}^{2n}
$
and
$
\Ker(A_\kappa + B_\kappa - I_n)\text{ in } {\Bbb Z}^{n}
$
are the torsion free part of 
$
{\Bbb Z}^{2n} / (H_\kappa - I_{2n}){\Bbb Z}^{2n}
$
and
that of
$ 
{\Bbb Z}^{n} / (A_\kappa + B_\kappa - I_n){\Bbb Z}^n
$
respectively,
so that they are isomorphic to each other.
\end{pf}
Therefore we reach the following theorem.
\begin{thm}
The $C^*$-algebra ${\cal O}_{{\cal H}_\kappa^{A,B}}$
associated with the 
Hilbert $C^*$-quad module ${\cal H}_\kappa^{A,B}$
defined by  commuting matrices
$A, B$
and a specification $\kappa$
is isomorphic to
the Cuntz-Krieger algebra 
${\cal O}_{H^{A,B}_\kappa}$
for the matrix 
$H^{A,B}_\kappa$.
Its K-groups 
$K_*({\cal O}_{H^{A,B}_\kappa})$
are computed as 
\begin{align*}
K_0({\cal O}_{H^{A,B}_\kappa}) & 
={\Bbb Z}^{n} / (A_\kappa + B_\kappa - I_n){\Bbb Z}^n, \\
K_1({\cal O}_{H^{A,B}_\kappa}) & 
= \Ker(A_\kappa + B_\kappa - I_n) \text{ in } {\Bbb Z}^n,
\end{align*}
where
$n = |\Omega_\kappa|$.
\end{thm}
We will finally present a concrete example.
For $1 < N, M \in {\Bbb N}$,
let $A$ and $B$
be the $1 \times 1$ matrices 
$[N]$
and
$[M]$
respectively.
The directed graph $G_A$
associated to the matrix
$A =[N]$
is a graph consists of a vertex denoted by $v$
 with $N$-self directed loops denoted by $E_A$.
Similarly 
the directed graph $G_B$
consists of the vertex $v$ with $M$-self directed loops
denoted by $E_B$.
We fix a specification
$\kappa: E_A \times E_B \longrightarrow E_B \times E_A$
defined by exchanging
$\kappa(\alpha,a) = (a,\alpha)$
for
$(\alpha,a) \in E_A \times E_B$.
Hence 
$\Omega_\kappa = E_A \times E_B$
so that 
$|\Omega_\kappa | = | E_A | \times | E_B| = N \times M$.
We then know
$\kappa_A((\alpha,a),(\delta,b)) = 1 $
if and only if 
$b=a$.
And
$\kappa_B((\alpha,a),(\beta,d)) = 1 $
if and only if 
$\beta=\alpha$ as in the following figures
respectively.
\begin{equation*}
{\begin{CD}
\circ @>\alpha>> \circ @>\delta>> \\
@V{a}VV  @V{a=b}VV @. 
\end{CD}
}
\qquad 
\text{and}
\qquad 
\begin{CD}
\circ @>\alpha>>  \\
@V{a}VV   @. \\
\circ @>\alpha=\beta >>  \\
@V{d}VV   @. 
\end{CD}
\end{equation*}
In particular,
for the case 
$N=2$ and $M=3$,
we write
$E_A = \{1,2\}, E_B =\{1,2,3\}$
and 
$\Omega_\kappa$ as
$$
\Omega_\kappa =E_A \times E_B
=
\{ (1,1), (1,2), (1,3), (2,1), (2,2), (2,3)\}.
$$
The $6\times 6$ matrices 
$A_\kappa$ and $B_\kappa$ 
are written along the above ordered basis in order as:
\begin{equation*}
A_\kappa =
\begin{bmatrix}
1 & 0 & 0 & 1 & 0 & 0 \\
0 & 1 & 0 & 0 & 1 & 0 \\
0 & 0 & 1 & 0 & 0 & 1 \\
1 & 0 & 0 & 1 & 0 & 0 \\
0 & 1 & 0 & 0 & 1 & 0 \\
0 & 0 & 1 & 0 & 0 & 1 
\end{bmatrix}
\qquad
\quad
\text{ and }
\qquad
\quad
B_\kappa =
\begin{bmatrix}
1 & 1 & 1 & 0 & 0 & 0 \\
1 & 1 & 1 & 0 & 0 & 0 \\
1 & 1 & 1 & 0 & 0 & 0 \\
0 & 0 & 0 & 1 & 1 & 1 \\
0 & 0 & 0 & 1 & 1 & 1 \\
0 & 0 & 0 & 1 & 1 & 1 
\end{bmatrix}
\end{equation*}
respectively so that we have
\begin{equation*}
A_\kappa + B_\kappa - I 
=
\begin{bmatrix}
1 & 1 & 1 & 1 & 0 & 0 \\
1 & 1 & 1 & 0 & 1 & 0 \\
1 & 1 & 1 & 0 & 0 & 1 \\
1 & 0 & 0 & 1 & 1 & 1 \\
0 & 1 & 0 & 1 & 1 & 1 \\
0 & 0 & 1 & 1 & 1 & 1 
\end{bmatrix}.
\end{equation*}
It is easy to see that 
\begin{equation*}
{\Bbb Z}^6/ (A_\kappa + B_\kappa - I){\Bbb Z}^6
\cong {\Bbb Z}/8{\Bbb Z},
\qquad
\Ker(A_\kappa + B_\kappa - I) \text{ in }{\Bbb Z}^6
\cong \{ 0 \}.
\end{equation*}
Therefore the $C^*$-algebra 
${\cal O}_{H^{A,B}_\kappa}$
for $A=[2], B=[3]$
and $\kappa=$exchange
is a Cuntz-Krieger algebra stably isomorphic 
to
the Cuntz algebra ${\cal O}_9$,
whereas the $C^*$-algebra 
${\cal O}_{[2], [3]}^\kappa$
considered in \cite{MaPre2011}
is isomorphic ${\cal O}_2 \otimes {\cal O}_3$
which is isomorphic to ${\cal O}_2$.
We will further study 
these $C^*$-algebras 
${\cal O}_{H^{A,B}_\kappa}$
in a forthcoming paper.

\end{document}